\documentclass[english]{smfart}
\usepackage{smfenum,bull,graphicx,xspace,upref,amssymb,url,smfthm}
\usepackage[all,dvips,cmtip]{xy}

\title[Gauss-Manin systems and~Frobenius~structures~(I)]
{Gauss-Manin systems, Brieskorn lattices and~Frobenius~structures~(I)}

\alttitle{Syst\`emes de Gauss-Manin, r\'eseaux de Brieskorn et structures de Frobenius~(I)}

\author[A.~Douai]{Antoine Douai}
\address{UMR 6621 du CNRS\\ Laboratoire J.A. Dieudonn\'e\\
Universit\'e de Nice\\
Parc Valrose\\
06108 Nice cedex 2\\
France}
\email{douai@math.unice.fr}

\author[C.~Sabbah]{Claude Sabbah}
\address{UMR 7640 du CNRS\\
Centre de Math{\'e}matiques\\
{\'E}cole polytechnique\\
\hbox{F--91128} Palaiseau cedex\\
France}
\email{sabbah@math.polytechnique.fr}
\urladdr{http://www.math.polytechnique.fr/cmat/sabbah/sabbah.html}

\makeatletter
\def\mainmatter{\renewcommand{\baselinestretch}{1.1}\normalfont}
\def\backmatter{\renewcommand{\baselinestretch}{1}\normalfont}

\numberwithin{equation}{section}
\NumberTheoremsAs{equation}
\def\thesubsection{\thesection.\alph{subsection}}

\newdir^{ (}{{}*!/-5pt/\dir^{(}}
\newdir_{ (}{{}*!/-5pt/\dir_{(}}

\def\@tempa{english}
\ifx\@tempa\smf@language
   \def\remsname{Remarks}%
   \def\exemsname{Examples}%
   \def\rgname{rk}%
\else
   \def\remsname{Remarques}%
   \def\exemsname{Exemples}%
   \def\rgname{rg}%
\fi

\newcommand{\RedefinitSymbole}[1]{%
\expandafter\let\csname old\string#1\endcsname=#1
\let#1=\relax
\newcommand{#1}{\csname old\string#1\endcsname\,}%
}
\RedefinitSymbole{\forall} \RedefinitSymbole{\exists}

\newenvironment{theoreme}{\begin{theo}}{\end{theo}}
\newenvironment{proposition}{\begin{prop}}{\end{prop}}
\newenvironment{lemme}{\begin{lemm}}{\end{lemm}}
\newenvironment{corollaire}{\begin{coro}}{\end{coro}}

\newenvironment{definition}{\begin{defi}}{\end{defi}}

\newenvironment{remarque}{\begin{rema}}{\end{rema}}

\newenvironment{exemple}{\begin{exem}}{\end{exem}}
\newenvironment{exemples}{\begin{enonce}[remark]{\exemsname}}{\end{enonce}}

\newenvironment{Corollaire}{\skippointrait\begin{coro}}{\end{coro}}

\newenvironment{Remarques}{\skippointrait\begin{enonce}[remark]{\remsname}}{\end{enonce}}

\newenvironment{Exemples}{\skippointrait\begin{enonce}[remark]{\exemsname}}{\end{enonce}}

\let\wh\widehat
\let\wt\widetilde
\let\ov\overline
\let\emptyset\varnothing
\let\ldots\dots

\let\moins\smallsetminus
\let\epsilon\varepsilon
\let\epsilong\varepsilong
\let\nablaf\bigtriangledown
\let\dag\dagger

\newcommand{\bbullet}{{\scriptscriptstyle\bullet}}

\newcommand{\ooplus}{\mathop\oplus\limits}
\newcommand{\ootimes}{\mathop\otimes\limits}

\newcommand{\module}[1]{\left\vert#1\right\vert}

\DeclareMathOperator{\coker}{Coker}

\DeclareMathOperator{\id}{Id}

\newcommand{\isom}{\stackrel{\sim}{\longrightarrow}}

\renewcommand{\ker}{\mathop{\rm Ker}\nolimits}

\newcommand{\rg}{\mathop{\mathrm{\rgname}}\nolimits}

\newcommand{\cHom}{\mathop{\cH om}\nolimits}

\newcommand{\cExt}{\mathop{\cE xt}\nolimits}

\DeclareMathOperator{\Supp}{Supp}

\DeclareMathOperator{\Char}{Char}
\DeclareMathOperator{\DR}{DR}
\DeclareMathOperator{\ord}{ord}
\DeclareMathOperator{\Mod}{Mod}
\DeclareMathOperator{\SP}{SP}
\DeclareMathOperator{\Sp}{Sp}

\DeclareMathOperator{\image}{image}
\DeclareMathOperator{\cHolreg}{\cH\!\mathit{olreg}}
\DeclareMathOperator{\cHol}{\cH\!\mathit{ol}}

\newcommand{\loccit}{\emph{loc\ptbl cit}}

\def\cf{\textit{cf.}\kern.3em}
\def\eg{\textit{e.g.},\ }
\def\ie{\textit{i.e.},\ }
\def\resp{\textit{resp.}\kern.3em}

\newcommand{\T}{\S\kern .15em }
\newcommand{\ptbl}{.\kern .15em }

\newcommand{\lcr}{[\![}
\newcommand{\rcr}{]\!]}

\def\to{\mathchoice{\longrightarrow}{\rightarrow}{\rightarrow}{\rightarrow}}
\def\mto{\mathchoice{\longmapsto}{\mapsto}{\mapsto}{\mapsto}}
\def\hto{\mathrel{\lhook\joinrel\to}}

\let\iff\ssi

\def\MRE#1{\mathchoice{\xrightarrow{\textstyle~~#1~}}{\stackrel{#1}{\longrightarrow}}{}{}}

\def\cB{\mathcal{B}}

\def\cD{\mathcal{D}}
\def\cE{\mathcal{E}}
\def\cF{\mathcal{F}}
\def\cG{\mathcal{G}}
\def\cH{\mathcal{H}}
\def\cI{\mathcal{I}}
\def\cJ{\mathcal{J}}
\def\cK{\mathcal{K}}
\def\cL{\mathcal{L}}
\def\cM{\mathcal{M}}
\def\cN{\mathcal{N}}
\def\cO{\mathcal{O}}

\def\cU{\mathcal{U}}

\def\cZ{\mathcal{Z}}

\def\varepsilong{\boldsymbol{\varepsilon}}
\let\epsilong\varepsilong

\def\bR{\boldsymbol{R}}

\def\bbA{\mathbb{A}}
\def\CC{\mathbb{C}}
\def\bbD{\mathbb{D}}
\def\NN{\mathbb{N}}
\def\PP{\mathbb{P}}
\def\QQ{\mathbb{Q}}
\def\RR{\mathbb{R}}
\def\ZZ{\mathbb{Z}}

\DeclareMathAlphabet{\mathcalmaigre}{U}{eus}{m}{n}
\def\ccF{\mathcalmaigre{F}}
\def\ccP{\mathcalmaigre{P}}

\let\leq\leqslant
\let\geq\geqslant

\newcommand{\cDXt}{\cD_X[t]\langle\partial_t\rangle}
\newcommand{\cDYt}{\cD_Y[t]\langle\partial_t\rangle}
\newcommand{\cDXD}{\cD_{D\times X}}
\newcommand{\cDXDth}{\cD_{D\times X/D}\langle\!\langle\theta\rangle\!\rangle \langle\theta^{-1}\rangle}
\newcommand{\cDXT}{\cD_X[\tau]\langle\partial_\tau\rangle}
\newcommand{\Afu}{\mathbb{A}^{\!1}}
\newcommand{\Afuan}{\mathbb{A}^{\!1\an}}
\newcommand{\AfN}{\mathbb{A}^{\!N}}
\newcommand{\Afuh}{\wh{\mathbb{A}}^{\!1}}
\newcommand{\Afuhan}{\wh{\mathbb{A}}^{\!1\an}}
\newcommand{\an}{\mathrm{an}}
\newcommand{\coh}{\mathrm{coh}}
\newcommand{\pgood}{p\text{\rm-good}}
\newcommand{\cbbullet}{{\raisebox{1pt}{$\bbullet$}}}

\newcommand{\gE}{\mathfrak{E}}
\newcommand{\gr}{\mathrm{gr}}

\def\NC{{\rm(NC)}\xspace}

\def\sfrac#1#2{{#1}/{#2}}

\def\defin{\mathrel{:=}}

\def\numero{n${}^\circ$\kern2pt}
\makeatother

\begin{document}
\frontmatter

\dedicatory{\`A Fr\'ed\'eric Pham, avec amiti\'e}

\begin{abstract}
We associate to any convenient nondegenerate Laurent polynomial $f$ on the complex torus $(\mathbb{C}^*)^n$ a canonical Frobenius-Saito structure on the base space of its universal unfolding. According to the method of K.~Saito (primitive forms) and of M.~Saito (good basis of the Gauss-Manin system), the main problem, which is solved in this article, is the analysis of the Gauss-Manin system of $f$ (or its universal unfolding) and of the corresponding Hodge theory.
\end{abstract}

\begin{altabstract}
Nous associons \`a tout polyn\^ome de Laurent commode et non d\'eg\'en\'er\'e $f$ par rapport \`a son poly\`edre de Newton sur le tore complexe $(\mathbb{C}^*)^n$ une structure de Frobenius-Saito canonique sur la base de son d\'eploiement universel. En suivant la m\'ethode de K.~Saito (formes primitives) et de M.~Saito (bonnes bases du syst\`eme de Gauss-Manin), le probl\`eme principal, qui est r\'esolu dans cet article, consiste en l'analyse du syst\`eme de Gauss-Manin de $f$ (ou de son d\'eploiement universel) et de la th\'eorie de Hodge correspondante.
\end{altabstract}


\altkeywords{Syst\`eme de Gauss-Manin, r\'eseau de Brieskorn, vari\'et\'e de Frobenius}

\maketitle
\mainmatter

\section*{Introduction}
It is now well-known that it is possible to associate with any germ of isolated hypersurface singularity a canonical Frobenius-Saito structure\footnote{Throughout of this article, we use the words ``Frobenius structure'' or ``Saito structure'' in an equivalent way; some authors also call such a structure a ``conformal Frobenius structure'', as we always assume the existence of a Euler vector field.}
(also called a flat structure by K.~Saito) on the (germ of the) base space of its miniversal unfolding. The theory of primitive forms \cite{KSaito83b} allows such a construction, together with Hodge theory of the Brieskorn lattice, as developed by M.~Saito \cite{MSaito89,MSaito91} (see also \cite{Hertling00} for a presentation of the theory). Such a construction is local, but also canonical.

In \cite{Sabbah96c}, the second author has indicated how to extend such a construction to polynomials on $\CC^n$ which are convenient and nondegenerate with respect to their Newton polyhedron, a notion defined by A.G.~Kouchnirenko in \cite{Kouchnirenko76}. The reason for extending the previous work of K.~Saito and M.~Saito is that, in some aspect of mirror symmetry theory, the Landau-Ginzburg potential often takes the form of the universal unfolding of a function on a smooth complex affine variety and is not necessarily reducible to the unfolding of a germ of complex isolated singularity.

A motivation to extend furthermore the construction to convenient and nondegenerate Laurent polynomials, which is the main result of this article, comes from the computation of S.~Barannikov \cite{Barannikov00} for the particular case of the Laurent polynomial $f(u_1,\dots,u_n)=u_1+\cdots+u_n+1/u_1\cdots u_n$: he associates to such a Laurent polynomial a Frobenius structure on the germ $(\CC^{n+1},0)$, which is shown to be equivalent to that coming from the quantum cohomology of $\PP^n(\CC)$ (defined \eg in \cite{Manin96}). One finds in \cite{H-V00} other such polynomials\footnote{The second author is grateful to R.~Kauffman for pointing him this article.}. In the second part \cite{D-S02b} of this article, we will extend the computation of S.~Barannikov to the case of the linear form $w_0u_0+\cdots+w_nu_n$ restricted to the torus $\prod_iu_i^{w_i}=1$, where $(w_0,\dots,w_n)$ is any sequence of positive integers such that $\gcd(w_0,\dots,w_n)=1$; in such a case, the canonical Frobenius structure is completely determined by its initial data, that we compute explicitly.

In this article, we give a detailed presentation of the construction of the Frobenius structure attached to a regular tame function on an affine manifold (when a primitive homogeneous form exists, see below). We first give the main properties of the Gauss-Manin system and the Brieskorn lattice of any unfolding of such a function. Essentially, we show that
all expected results (if we refer to the local situation) hold: holonomicity, regularity, duality... Such results were only sketched in \cite{Sabbah96c}. In Appendix~\ref{app:A}, we recall with some details basic results written by B.~Abdel-Gadir in \cite{AbdelGadir97a,AbdelGadir97b}, as these papers are hardly available.

Once these basic results are obtained, the main tools that one has to develop in order to mimic K.+M.~Saito's construction of the Frobenius structure are:
\begin{itemize}
\item
the Hodge theory for the Brieskorn lattice of the tame function $f$ on a smooth complex affine variety,
\item
the existence of a primitive and homogeneous section of the Brieskorn lattice.
\end{itemize}
The first one was achieved in \cite{Sabbah96a,Sabbah96b}. In Appendix \ref{app:B}, we recall the analogue, in the affine situation, of the method of M.~Saito to obtain a solution of Birkhoff's problem for the Brieskorn lattice. This slightly simplifies the approach in \cite{Sabbah96b}. Moreover, we make the link with the construction of good bases given in \cite{Douai99}.

We are then mainly interested in finding a primitive homogeneous section corresponding to the minimal exponent. Natural candidates are volume forms. In the singularity case, the minimal exponent has multiplicity one \cite{MSaito91}, so such a primitive homogeneous section is essentially unique, if it exists. In fact, any volume form gives rise to such a primitive homogeneous section. We show that a similar result holds for convenient nondegenerate Laurent polynomials on $(\CC^*)^n$. It was asserted in \cite{Sabbah96b} that the same result holds for convenient nondegenerate polynomials on $\CC^n$. However, let us emphasize that a supplementary assumption is understated there: \emph{the linear forms defining the Newton boundary have nonnegative coefficients}. Without this assumption, the volume form may not be homogeneous, the minimal exponent may not have multiplicity one, and we cannot assert the existence of any primitive homogeneous section.

 A basic tool in the convenient nondegenerate case, which goes back to \cite{K-V85,MSaito88}, consists in the identification between the Hodge exponents at infinity (\ie the Hodge spectrum) and the Newton exponents. A simple proof of this result is given in \T\ref{sec:laurent}. It also applies to the case of convenient nondegenerate polynomials on $\CC^n$ which was treated in \cite{Sabbah96b}.

\section{Partial Fourier transform of regular holonomic $\cD$-modules with~lattices} \label{sec:Fourier}

\subsection{Equivalences of categories}
Let $D\subset \Afuan=\PP^1\moins\{\infty\}$ be a nonempty open disc with coordinate $t$ and let $X$ be a complex manifold. Let $\Sigma\subset D\times X$ be a closed analytic hypersurface. We assume throughout this article that the restriction to $\Sigma$ of the projection $p:D\times X\to X$ is \emph{finite} (\ie proper with finite fibres). On the other hand, we simply denote by $\infty$ the divisor $\{\infty\}\times X\subset\PP^1\times X$.

We denote by $\cD_{D\times X}$ (\resp $\cD_{\PP^1\times X}$) the sheaf of holomorphic differential operators on $D\times X$ (\resp on $\PP^1\times X$). Given a holonomic $\cD_{D\times X}$-module $\cM$, we say that the \emph{singular locus} of $\cM$ is contained in $\Sigma$ if $\cM_{|D\times X\moins\Sigma}$ is $\cO$-locally free of finite rank.

Let $\varphi:X'\to X$ be a holomorphic map of complex manifolds. Then $\varphi^{+0}\cM\defin \cO_{D\times X'}\otimes_{\varphi^{-1}\cO_{D\times X}}\varphi^{-1}\cM$ is naturally equipped with a $\cD_{D\times X'}$-module structure, relative to which it is holonomic and with singular locus contained in $\Sigma'=X'\times_X\Sigma$.

\begin{proposition}[Malgrange]\label{prop:extension}
Let $\cM$ be a holonomic $\cD_{D\times X}$-module with singular locus contained in $\Sigma$. There exists a unique (up to isomorphism) $\cD_{\PP^1\times X}$-module $\wt\cM$ satisfying the following properties:
\begin{itemize}
\item
$\wt\cM$ is holonomic, with singular locus contained in $\Sigma\cup\infty$,
\item
$\wt\cM$ has regular singularities along $\infty$,
\item
$\wt\cM$ is localized along $\infty$, \ie $\wt\cM=\cO_{\PP^1\times X}(*\infty) \otimes_{\cO_{\PP^1\times X}}\wt\cM$,
\item
$\wt\cM_{|D\times X}=\cM$.
\end{itemize}
Any morphism $\cM_1\to\cM_2$ can be lifted in a unique way to a morphism $\wt\cM_1\to\wt\cM_2$. If $\bbD$ denotes the duality functor on the category left holonomic $\cD$-modules, we have $\wt{\bbD\cM}=(\bbD\wt\cM)(*\infty)$. Last, the $\wt{\kern5pt\cbbullet\kern5pt}$ functor commutes with base change $\varphi^{+0}$ with respect to $\varphi:X'\to X$.
\end{proposition}

\begin{proof}
See Appendix, \T\ref{proof:extension}.
\end{proof}

We consider the following categories:
\begin{itemize}
\item
$\cHolreg_\Sigma(\cD_{D\times X})$: the objects are regular holonomic $\cD_{D\times X}$-modules with singular locus contained in $\Sigma$, and the morphisms are morphisms of $\cD$-modules; the full subcategory of modules which satisfy the supplementary condition \NC (see below) will be denoted by $\cHolreg_{\Sigma,\NC}(\cD_{D\times X})$.

\item
$\cHolreg_\Sigma(\cD_{\Afuan\times X})$: the objects are regular holonomic $\cD_{\Afuan\times X}$-modules with singular locus contained in $\Sigma$, and the morphisms are morphisms of $\cD$-modules;

\item
$\cHolreg_{\Sigma\cup\infty}\big(\cD_{\PP^1\times X}(*\infty)\big)$: the objects are regular holonomic $\cD_{\PP^1\times X}(*\infty)$-modules, with singular locus contained in $\Sigma\cup\infty$ (recall, \cf \T\ref{subsec:algebr}, that these are holonomic $\cD_{\PP^1\times X}$-modules), and the morphisms are morphisms of $\cD$-modules;

\item
$\cHolreg_\Sigma(\cDXt)$: the objects are $\cDXt$-modules which are holonomic and regular even at infinity, with singular locus contained in $\Sigma$ (\cf Def\ptbl\ref{def:singetreg}), and the morphisms are morphisms of $\cDXt$-modules.
\end{itemize}

We have natural functors between these categories, described in the diagram below, where
\begin{itemize}
\item
$\text{restr.}$ is the usual restriction,
\item
$\text{an}$ is the ``analytization functor'' $\cD_{\Afuan\times X}\otimes_{\cDXt}$,
\item
$p$ denotes the natural projection to $X$.
\end{itemize}
\begin{figure}[htb]
\[
\xymatrix@C=-1cm@R=1.5cm{
&\cHolreg_{\Sigma\cup\infty}\big(\cD_{\PP^1\times X}(*\infty)\big)\ar[dl]_-{\text{restr.\;}}\ar[rd]^-{p_*}&\\
\cHolreg_\Sigma(\cD_{\Afuan\times X})\ar[rd]_-{\text{restr.}}&&\cHolreg_\Sigma(\cDXt)\ar[ll]_(.6){\text{an}}\\
&\cHolreg_\Sigma(\cD_{D\times X})\ar[uu]_(.3){\wr}|!{[ur];[ul]}\hole
}
\]
\par\smallskip
\centerline{Analytization diagram}
\end{figure}

\begin{Corollaire}[\cite{AbdelGadir97a}]\label{cor:equivcat}
\begin{enumerate}
\item\label{cor:equivcat1}
Any loop in this diagram is a functor isomorphic to the identity of the category which is the origin of the loop.
\item\label{cor:equivcat2}
This diagram commutes with base change, \ie if $f:X'\to X$ is any holomorphic map, and if $\cH^kf^+$ denotes the $k$-th inverse image functor between the corresponding categories on $X$ and $X'$, the diagrams on $X$ and $X'$ correspond each other through $\cH^kf^+$.
\item\label{cor:equivcat3}
This diagram commutes with proper direct image with respect to $X$, \ie with the functors $\cH^kf_+$ is $f:X\to X'$ is a proper holomorphic map (\eg a closed embedding).
\item\label{cor:equivcat4}
This diagram commutes with the natural duality functor of each category.
\end{enumerate}
\end{Corollaire}

\begin{proof}
See Appendix, \S\T\ref{subsec:algebr} and \ref{subsec:hol}.
\end{proof}

\begin{Remarques}\label{rem:tildepairs}
\begin{enumerate}
\item\label{rem:tildepairs1}
In the appendix, the previous statements are proved without the regularity assumption along $\Sigma$. That regularity along $\Sigma$ is preserved by the functors $\cH^kf^+$, $\cH^kf_+$ or duality is now a well-known result (see, \eg \cite{Mebkhout89}).

\item
\label{rem:tildepairs2}
Let $\cM$ be as in Proposition~\ref{prop:extension} and let $\cM_0$ be a $\cO_{D\times X}$-coherent submodule of $\cM$. Let us define $\wt\cM_0$ as a coherent $\cO_{\PP^1\times X}(*\infty)$-submodule of $\wt \cM$: let $\cM'=\cD_{D\times X}\cdot\cM_0$; then $\cM'$ satisfies the same properties as $\cM$ does, and moreover $\cM'$ and $\cM_0$ coincide out of $\Sigma$ (``out of'' means ``after tensoring with meromorphic functions having poles along $\Sigma$''); then glue $\cM_0$ with $\wt\cM'_{|\PP^1\times X\moins\Sigma}$ to get $\wt\cM_0$. By Remark~\ref{rem:lattice}, Cor\ptbl\ref{cor:equivcat}\eqref{cor:equivcat1} holds for such pairs $(\cM,\cM_0)$.
\end{enumerate}
\end{Remarques}

\begin{corollaire}\label{cor:locfreeextension}
Let $\cM$ be as in Proposition~\ref{prop:extension}. Assume moreover that there exists a $\cO_{D\times X}$-submodule $\cM_0\subset\cM$ which is free of finite rank and such that $\cM=\cO_{D\times X}(*\Sigma)\otimes_{\cO_{D\times X}}\cM_0$. Then $M_0\defin p_*\wt\cM_0$ is locally free of finite rank as a $\cO_X[t]$-module and $M$ is locally free of finite rank as a $\cO_X[t](*\Sigma)$-module.
\end{corollaire}

\begin{proof}
See Appendix, \T\ref{proof:locfreeextension}.
\end{proof}

\subsubsection*{The noncharacteristic assumption}
We say that $p$ is noncharacteristic with respect to $M$ or $\cM$ if the following condition \NC is satisfied:

\smallskip
\noindent\NC \hfill\emph{The fibres of $p$ are noncharacteristic with respect to $\Char\cM$}.\hfill\mbox{}

\smallskip
\noindent
The characteristic variety $\Char\cM$ is equal to a union of the zero section $T^*_{D\times X}(D\times X)$ and of conormal spaces $T^*_{Z_i}(D\times X)$, where $Z_i$ is a irreducible closed analytic subset of $D\times X$ contained in $\Sigma$. Among the $Z_i$ are the irreducible components of $\Sigma$.

In geometrical terms, Condition \NC means that, for any $Z_i$, any limit at $(t,x)$ of tangent hyperplanes to $Z_i$ is transversal to $p^{-1}(x)$. In other words, for any $x\in X$, the conormal space $T^*_{Z_i}(D\times X)$ cuts $D\times T^*_xX$ at most along $D\times\{0\}$. This implies that the restriction to $T^*_{Z_i}(D\times X)$ of the natural projection $T^*(D\times X)\to (T^*D)\times X$ is finite. This also implies that the restriction of $p$ to $Z_i$ is finite. In other words, the only components of the characteristic variety of $\cM$ other than the zero section are the $T^*_{\Sigma_i}(D\times X)$, where $\Sigma_i$ are the irreducible components of $\Sigma$.

\begin{figure}[htb]
\centerline{\includegraphics[scale=1]{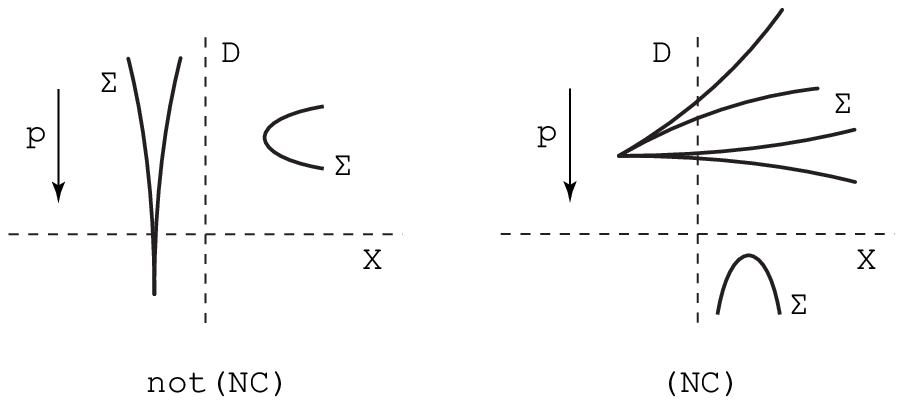}}
\caption{}
\end{figure}

We write $(T^*D)\times X=D\times X\times \CC^\vee$. The following is classical (noticing that the fibres of $p$ are noncharacteristic with respect to $\wt\cM$ along $\infty$).

\begin{lemme}\label{lem:NC}
If $\cM$ is holonomic on $D\times X$ with singular locus contained in $\Sigma$ and satisfies \NC, then $\cM$ is regular along $\Sigma$, $\cM$ is $\cD_{D\times X/X}$-coherent with relative characteristic variety contained in the union of the zero section of $(T^*D)\times X$ and $\Sigma\times\CC^\vee$. Moreover, $M=p_*\wt\cM$ is $\cO_X[t]\langle\partial_t\rangle$-coherent.\hfill\qed
\end{lemme}

\subsubsection*{Equivalence of categories with lattices} When working with lattices, we will always assume that the noncharacteristic assumption \NC is satisfied.

\begin{definition}[Lattices]\label{def:lattices}
Let $\cM,\wt M, M$ be as above and satisfying the noncharacteristic assumption \NC. A \emph{lattice} $\cM_0$ in $\cM$ is a $\cO_{D\times X}$-coherent submodule of $\cM$ which generates $\cM$ as a $\cD_{D\times X/X}$-module.
We similarly have the notion of a lattice in a $\cD_{\Afuan\times X}$-module and in a $\cDXt$-module. In a $\cD_{\PP^1\times X}(*\infty)$-module, a \emph{lattice} $\wt \cM_0$ is a coherent $\cO_{\PP^1\times X}(*\infty)$-module which generates $\wt\cM$ as a $\cD_{\PP^1\times X/X}(*\infty)$-module.

We say that the lattice $\cM_0$ (\resp $\wt\cM_0$, $M_0$) has \emph{Poincar\'e rank one} if, for any vector field $\xi$ on $X$, one has $\xi\cM_0\subset\cM_0+\partial_t\cM_0$ (\resp ...).
\end{definition}

Notice that, if $\cM_0$ is a lattice in $\cM$ (\resp $\wt\cM_0$ is a lattice in $\wt\cM$), then $\cM_0=\cM$ (\resp $\wt\cM_0=\wt\cM$) out of $\Sigma$ (``out of'' is as in Remark~\ref{rem:tildepairs}\eqref{rem:tildepairs2}).

\begin{corollaire}\label{cor:equivcatlat}
Same statement as that of Cor\ptbl\ref{cor:equivcat}\eqref{cor:equivcat1} for the categories of pairs $(\cM,\cM_0)$ of a $\cD$-module with a lattice.
\end{corollaire}

\begin{proof}
See Appendix, Remark \ref{rem:lattice}.
\end{proof}
\subsection{Partial Fourier transform}
Consider the isomorphism of algebras $$\cDXt\isom\cDXT$$ which is the identity on $\cD_X$ and such that $t\mto -\partial_\tau$ and $\partial_t\mto \tau$. Any $\cDXt$-module $M$ becomes a $\cDXT$-module \emph{via} this isomorphism, that we denote $\wh M$ and that we call the \emph{partial Fourier-Laplace transform of $M$ in the $t$-direction}.

Clearly, if $M$ is $\cDXt$-coherent, then $\wh M$ is $\cDXT$-coherent and conversely. Using the homological characterization of holonomic modules (\cf Appendix, \T\ref{subsec:algebr}), one gets:

\begin{proposition}\label{prop:Fourierhol}
$M$ is holonomic iff $\wh M$ is so.\hfill\qed
\end{proposition}

Denote by $\bbD M$ the left $\cDXt$-module associated to the right $\cDXt$-module $\cExt_{\cDXt}^{\dim X+1}(M,\cDXt)$. It is also holonomic. Notice that we clearly have
\[
\cExt_{\cDXt}^{\dim X+1}(M,\cDXt)\!\wh{\phantom{\cM}}=\cExt_{\cDXT}^{\dim X+1}(\wh M,\cDXT).
\]
In local coordinates, the left-right transformation is expressed using the transposition of operators $P\mto{}^t\!P$.

Denote by $P\mto\ov{P}$ the involution of the algebra $\cDXT$ which is the identity on $\cD_X$ and such that $\ov\tau=-\tau$ and $\ov\partial_\tau=-\partial_\tau$. One has ${}^t\wh P=\ov{\wh{{}^t\!P}}$. One deduces that
\begin{equation}\label{eq:dualFourier}
\wh{\bbD M}=\ov{\bbD\wh M}.
\end{equation}

It follows from \cite{Kashiwara78} that
$$G\defin\cD_X[\tau,\tau^{-1}]\langle\partial_\tau\rangle\ootimes_{\cDXT}\wh M=\cO_X[\tau,\tau^{-1}]\ootimes_{\cO_X[\tau]}\wh M$$ is $\cDXT$-holonomic if $\wh M$ is so. Put $\theta=\tau^{-1}$ and identify the action of $t=-\partial_\tau$ on $G$ with that of $\theta^2\partial_\theta$. It also follows from \cite{Kashiwara78} that $G$ is holonomic as a $\cD_X[\theta]\langle\partial_\theta\rangle$-module.

\begin{remarque}[Inverting $\partial_t$]\label{rem:invertdt}
Consider the $\cDXt$-module $$M[\partial_t^{-1}]\defin\cD_X[t]\langle\partial_t,\partial_t^{-1}\rangle\ootimes_{\cDXt}M.$$ This is a holonomic $\cDXt$-module, being the inverse partial Fourier transform of $G$. The kernel and the cokernel of the natural morphism $M\to M[\partial_t^{-1}]$ are thus holonomic. Moreover, they take the form $p^+N',p^+N''$ for some holonomic $\cD_X$-modules $N',N''$ (this follows from Kashiwara's equivalence applied to holonomic $\cDXT$-modules supported on $\tau=0$, after partial Fourier transform).

If $M$ is an object of $\cHolreg_\Sigma(\cDXt)$, then so is $p^+N'$, and this implies that $N'$ is $\cO_X$-locally free of finite rank (\ie a vector bundle with a flat connection).

If moreover $p$ satisfies Condition \NC with respect to $M$, it follows from Cor\ptbl\ref{cor:Fourierreg} below that $N''$ is $\cO_X$-locally free of finite rank, from what we also deduce that $M$ is $\cO_X$-flat and that $M[\partial_t^{-1}]$ also belongs to $\cHolreg_\Sigma(\cDXt)$. Let us sketch a proof of the local $\cO_X$-freeness of $N''$: choose any local coordinate $x$ on $X$, and consider, as in \T\ref{subsubsec:nearby}, the vanishing cycles functor; then $\phi_x\wh M^{\an}=0$ (\cf Appendix \T\ref{proof:Fourierreg}, last part of the proof of Theorem \ref{th:Fourierreg}) and, by Corollary \ref{cor:Fourierreg}, we also have $\phi_x\wh M[\tau^{-1}]^{\an}=0$; as the functor $\phi_x$ is exact on holonomic $\cD_{\Afuhan\times X}$-modules (see, \eg \cite{M-S86}), we conclude that $\phi_xN''=0$; as this vanishing result holds for any local coordinate on $X$, this implies that $N''$ is a vector bundle with a flat connection.
\end{remarque}

\Subsection{Regularity of the partial Fourier transform}

\begin{theoreme}\label{th:Fourierreg}
Let $M$ be a holonomic $\cDXt$-module with singular locus contained in $\Sigma$. Assume that $M$ is regular, including along $\infty$. Then
\begin{enumerate}
\item\label{th:Fourierreg1}
the partial Fourier transform $\wh M^{\an}$ is regular on $\Afuhan\times X$ (but usually not at~$\wh\infty$),
\item\label{th:Fourierreg2}
under the noncharacteristic assumption \NC, $\wh M^{\an}$ is smooth out of \hbox{$\{0\}\times X\subset\Afuhan\times X$}.
\end{enumerate}
\end{theoreme}

\begin{proof}
See Appendix, \T\ref{proof:Fourierreg}.
\end{proof}

\begin{corollaire}\label{cor:Fourierreg}
Under the same conditions and with the noncharacteristic assumption \NC, the $\cO_X[\tau,\tau^{-1}]$-module $G$ is coherent.
\end{corollaire}

\begin{proof}
Denote by $\wh p:\wh \PP^1\times X\to X$ the projection. Under the noncharacteristic assumption, $\wh p^*G$ is $\cO_{\wh\PP^1\times X}(*\wh 0\cup\wh \infty)$-coherent, as a consequence of \ref{th:Fourierreg}\eqref{th:Fourierreg2} and \cite{Kashiwara78}. Therefore, $G$ is $\cO_X[\tau,\tau^{-1}]$-coherent.
\end{proof}

Starting from a regular holonomic $\cD_{D\times X}$-module $\cM$ with singular locus contained in $\Sigma$ and satisfying \NC, we have considered the following objects:
\begin{itemize}
\item
its extension $\wt\cM$ as a holonomic $\cD_{\PP^1\times X}(*\infty)$-module,
\item
the algebraization $M=p_*\wt\cM$, which is a regular holonomic $\cDXt$-module,
\item
the Fourier transform $\wh M$ of $M$, which is a holonomic $\cDXT$-module,
\item
the localization $G=\wh M[\tau^{-1}]$ of $\wh M$ along $\tau=0$, which is a holonomic $\cDXT$-module and also a coherent $\cO_X[\tau,\tau^{-1}]$-module with a flat connection, regular along $\tau=0$. We denote by $\rg G$ the generic rank of $G^\an$ as a $\cO_{\Afuhan\times X}$-module.
\end{itemize}

\begin{definition}\label{def:cF}
Denote by $\ccF$ the composite functor $\cM\mto G$ from $\cHolreg_{\Sigma,\NC}(\cD_{D\times X})$ to the category of coherent $\cO_X[\tau,\tau^{-1}]$-modules with a flat connection, regular along $\tau=0$.
\end{definition}
\begin{remarque}\label{rem:locfree}
Under the condition \NC, $G$ is not far from being locally $\cO_X[\tau,\tau^{-1}]$-free. More precisely, locally on $X$, there exists $N$ such that $G\oplus\cO_X[\tau,\tau^{-1}]^N$ is $\cO_X[\tau,\tau^{-1}]$-locally free. Indeed, let $D'$ be a neighbourhood of $\tau=\infty$ in $\PP^1$ with coordinate $\theta=\tau^{-1}$ and let $\cG=G^\an$ on $D'\times X$. As $G$ is regular along $\theta=\infty$ and has pole along $\theta=0$ only, we have $G=\wh p_*\wt\cG$. On the other hand, since $\cG$ is a $\cO_{D'\times X}[\theta^{-1}]$-module with flat connection, it is locally stably free (\cf~\cite{Malgrange95}), so that, locally on $X$, there exists $N$ such that $\cG\oplus\cO_{D'\times X}[\theta^{-1}]^N$ is free. Equip $\cO_{D'\times X}[\theta^{-1}]^N$ with the trivial connection $d$. We clearly have
\[
[\cG\oplus(\cO_{D'\times X}[\theta^{-1}]^N,d)]^\sim= \wt{\cG}\oplus \wh p^*(\cO_{X}[\theta,\theta^{-1}]^N,d).
\]
Apply then Corollary~\ref{cor:locfreeextension}.
\end{remarque}

\begin{Exemples}
\begin{enumerate}\renewcommand{\theenumi}{\alph{enumi}}
\item
The condition that $p:\Sigma\to X$ is finite insures that $\wh M^{\an}$ has regular singularities. If $p$ is only quasi-finite, this may not hold. If for instance $X$ is a disc with coordinate $x$, consider the Fourier transform $\wh M$ of the regular holonomic $\cDXt$-module of which $(xt-1)^{1/2}$ is a solution; then $\wh M$ contains the line $x=0$ in its singular locus and is irregular along this line: indeed, define $M$ as the quotient of $\cDXt$ by the left ideal generated by the operators
\[
(xt-1)\partial_x-\tfrac12 t,\quad (xt-1)\partial_t-\tfrac12 x,\quad x\partial_x-t\partial_t~;
\]
the Fourier transform $\wh M$ is the quotient of $\cDXT$ by the left ideal generated by
\[
(x\partial_\tau+1)\partial_x-\tfrac12 \partial_\tau,\quad (x\partial_\tau+1)\tau+\tfrac12 x,\quad x\partial_x+\partial_\tau\tau,
\]
which contains the operator $-x^2\partial_x+\tau+\tfrac12 x$ responsible for the irregularity along $x=0$ when $\tau\neq0$.

A more precise relationship between the characteristic variety of $M$ and the irregularity of $\wh M$ along $x=0$ is given in \cite[Chap\ptbl3, Prop\ptbl4.5.7]{Sabbah97}.

\item
Consider the system satisfied by the function $(t^2-x^3)^{1/2}$, \ie the quotient $M$ of $\cDXt$ by the left ideal generated by the operators
\[
P_1=-t\partial_t-\tfrac23x\partial_x+1,\quad P_2=-t\partial_x-\tfrac32x^2\partial_t.
\]
The singular locus $\Sigma$ is defined by $t^2=x^3$, so $p:\Sigma\to X$ is finite, but $\{x=0\}$ is characteristic with respect to $\Sigma$, hence it is characteristic with respect to $M$ (as $T^*_\Sigma(\Afuan\times X)$ is a component of the characteristic variety of $M$).

The Fourier transform is the quotient of $\cDXT$ by the left ideal generated by
\[
\wh P_1=\tau\partial_\tau-\tfrac23x\partial_x+2,\quad \wh P_2=\partial_\tau\partial_x-\tfrac32\tau x^2,
\]
which also contains
\[
\partial_x\wh P_1-\tau\wh P_2=-\tfrac23x\partial_x^2+\tfrac43\partial_x+\tfrac32\tau^2x^2, \quad\partial_\tau\wh P_1+\tfrac23 x\wh P_2=\tau^2\partial_\tau+3\partial_\tau-\tau^2x^2.
\]
One checks that, near a point where $\tau\neq0$, the characteristic variety contains the component $x=0,\wt\tau=0$, where $(\tau,x,\wt\tau,\xi)$ are the symplectic coordinates on \hbox{$T^*(\Afuan\times X)$}. Therefore, the singular locus of $\wh M$ contains $\{x=0\}$ and $G=\cO_X[\tau,\tau^{-1}]\otimes_{\cO_X[\tau]}\wh M$ is not $\cO_X[\tau,\tau^{-1}]$-coherent.
\end{enumerate}
\end{Exemples}

\subsection{Partial microlocalization}
Assume that we are given $\cM$ as in Corollary \ref{cor:equivcatlat}, with $\cM$ in $\cHolreg_\Sigma(\cD_{D\times X})$. In general, we cannot recover $M$ or $G$ in an ``algebraic'' way from $\cM$.
Put (abusing notation) $\cO_X\lcr\theta\rcr=\varprojlim_n\cO_X[\theta]/\theta^n$. We will show in this paragraph that it is possible to recover $G^\wedge\defin \cO_X\lcr\theta\rcr\otimes_{\cO_X[\theta]}G$ from the partially microlocalized object $\cM^\mu$.

Throughout the remaining of this section, we assume that the noncharacteristic assumption \NC is satisfied for $\cM$.

We work with the sheaf of \emph{formal} relative microdifferential operators on $D\times X$, that we denote by $\cE_{D\times X/X}$ (see, \eg \cite{Pham79,Schapira85}). Denote by $\theta$ the variable corresponding to $\partial_t^{-1}$ in $\cE_{D\times X/X}$. We implicitly restrict this sheaf to the image of the section $dt$ of $(T^*D)\times X$, that we identify with $D\times X$, so that we view $\cE_{D\times X/X}$ as a sheaf on $D\times X$. Local sections of $\cE_{D\times X/X}$ are formal series $\sum_{j\geq j_0} a_j(t,x)\theta^j$, where the $a_j$ are holomorphic functions defined on a common open set of $D\times X$. The commutation rule is given by
\[
\theta\cdot a(t,x)=\sum_{k=0}^\infty(-1)^k\Big(\frac{\partial}{\partial t}\Big)^k(a)\cdot \theta^{k+1}.
\]
One identifies $\cD_{D\times X/X}$ with the subring of $\cE_{D\times X/X}$ consisting of polynomials $\sum_{j\leq0}a_j(t,x)\theta^j$ in $\theta$ \emph{via} $\partial_t\mto\theta^{-1}$.

If $\cM$ is a coherent $\cD_{D\times X/X}$-module, we denote by
\begin{equation}\label{eq:Mmu}
\cM^\mu=\cE_{D\times X/X}\ootimes_{\cD_{D\times X/X}}\cM
\end{equation}
the corresponding (formal) partial microlocal module. Notice that $\cM^\mu$ comes equipped with a natural action of differential operators on $X$. In order to take this action into account, consider the sheaf
\begin{equation}\label{eq:cDXDth}
\cDXDth\defin\cD_{D\times X}\ootimes_{\cD_{D\times X/X}}\cE_{D\times X/X}
\end{equation}
of formal series $\sum_{j\geq j_0}P_j\theta^{j}$ where $P_j$ are sections of $\cD_{D\times X/D}$ defined on a fixed open set and of degree bounded independently of $j$. It contains $\cE_{D\times X/X}$ as a subsheaf, and we have $\cM^\mu=\cDXDth\otimes_{\cDXD}\cM$.

\begin{prop}\label{prop:muG}
Let $\cM$ be an object of $\cHolreg_\Sigma(\cD_{D\times X})$, let $M=p_*\wt \cM$ and $G=M[\partial_t^{-1}]$. Under the noncharacteristic assumption \NC,
\begin{enumerate}
\item\label{prop:muG1}
$p_*\cM^\mu$ is $\cO_X\lcr\theta\rcr[\theta^{-1}]$-coherent; the natural action of $t$ and of the vector fields on $X$ defines on it a flat connection, such that $t$ acts as $\theta^2\partial_\theta$;
\item\label{prop:muG2}
the natural morphism
\[
G^\wedge\defin\cO_X\lcr\theta\rcr[\theta^{-1}]\ootimes_{\cO_X[\theta,\theta^{-1}]}G\to p_*\cM^\mu
\]
is an isomorphism of $\cO_X\lcr\theta\rcr[\theta^{-1}]$-modules, compatible with the $\cD_X\lcr\theta\rcr[\theta^{-1}]\langle\partial_\theta\rangle$-module structure.
\end{enumerate}
\end{prop}

We have denoted by $\cD_X\lcr\theta\rcr[\theta^{-1}]\langle\partial_\theta\rangle=\cD_X\lcr\theta\rcr[\theta^{-1}]\langle\theta^2\partial_\theta\rangle$ the sheaf $\cO_X\lcr\theta\rcr[\theta^{-1}]\otimes_{\cO_X[\theta^{-1}]}\cD_X[\theta^{-1}]\langle t\rangle$, with the identification $t=\theta^2\partial_\theta$.

\begin{proof}
For the first point, use that $\cM^\mu$ is supported on $\Sigma$ and apply the standard preparation theorem.

For the second point, denote by $\cE_{\Afu\times X/X}$ the sheaf on $X$ having formal series $\sum_{j\geq j_0}a_j(t,x)\theta^j$ as sections, where the $a_j$ are sections of $\cO_X[t]$ on a common open set. It contains $\cO_X[t]\langle\partial_t\rangle$ as a subsheaf. We then have a morphism
\[
G^\wedge=\cO_X\lcr\theta\rcr[\theta^{-1}]\ootimes_{\cO_X[\theta^{-1}]}M
\to\cE_{\Afu\times X/X}\ootimes_{\cO_X[t]\langle\partial_t\rangle}M,
\]
hence a well defined analytization map $G^\wedge\to p_*\cM^\mu$. Both objects are equipped with a flat connection and the previous morphism is obviously compatible with the connections. The kernel and cokernel of it are therefore of the same kind. It follows from \cite{Malgrange95} that all these objects are locally stably $\cO_X\lcr\theta\rcr[\theta^{-1}]$-free. Hence, after adding a suitable power of $(\cO_X\lcr\theta\rcr[\theta^{-1}],d)$ to both objects, we may assume that, locally on $X$, they are $\cO_X\lcr\theta\rcr[\theta^{-1}]$-free, as well as the kernel and the cokernel of the analytization map. As a consequence, in order to show that the kernel and cokernel are zero, it is enough to show this fact after restriction to any $x^o\in X$.

Choose $x^o\in X$. By flatness of $\cE_{D\times X/X}$ over $\cD_{D\times X/X}$, and as $p$ is finite, the restriction to $x^o$ of $p_*\cM^\mu$ is $p_*$ of the microlocalized module associated to $M^o$. An analogous result also holds for $G$, hence for $G^\wedge$. Last, we know that the morphism $G^\wedge\to p_*\cM^\mu$ is an isomorphism when $X$ is reduced to a point (see, \eg \cite[Prop\ptbl2.3]{Sabbah96b}). Hence, $G^\wedge\to p_*\cM^\mu$ is an isomorphism after restriction to~$x^o$.
\end{proof}

\subsubsection*{Relation with microlocalization}
We still assume that $\cM$ satisfies \NC. Let $(c,x^o)\in\Sigma\subset D\times X$ and let $\{\eta_1,\dots,\eta_r\}$ be the inverse image of $(c,x^o;1)\in T^*_cD\times\{x^o\}$ in $\Char\cM\cap T^*_{(c,x^o)}(D\times X)$. Denote by $\cE_{D\times X,\eta}$ the germ at $(c,x^o;\eta)\in T^*_{(c,x^o)}(D\times X)$ of the sheaf of (formal) microdifferential operators on $D\times X$. Then by \NC and the Preparation Theorem for microdifferential operators,
\[
\cM^\mu_{(c,x^o;\eta_i)}\defin \cE_{D\times X,(c,x^o;\eta_i)}\ootimes_{\cD_{D\times X,(c,x^o)}}\cM_{(c,x^o)}
\]
has finite type over $\cO_{X,x^o}\lcr\theta\rcr[\theta^{-1}]$. Notice that, denoting on the right by $\cM^\mu$ the partial microlocalized module, we also have
\[
\cM^\mu_{(c,x^o;\eta_i)}= \cE_{D\times X,(c,x^o;\eta_i)} \ootimes_{\cDXDth_{x^o}} \cM^\mu_{(c,x^o;1)}.
\]
The natural morphism of $\cO_{X,x^o}\lcr\theta\rcr[\theta^{-1}]$-modules
\begin{equation}\label{eq:mumu}
\cM^\mu_{(c,x^o;1)}\to \ooplus_{i=1}^r\cM^\mu_{(c,x^o;\eta_i)}
\end{equation}
is compatible with the connections. It is moreover known that both $\cO_{X,x^o}\lcr\theta\rcr[\theta^{-1}]$-modules have the same rank. Arguing as in Proposition~\ref{prop:muG}\eqref{prop:muG2}, we conclude that this morphism is an isomorphism. This justifies the ambiguous notation $\cM^\mu$.

\subsection{Partial microlocalization and partial Fourier transform with a lattice}\label{subsec:FTlattice}

Recall that we assume \NC. Let $\cM_0$ be a lattice in $\cM$. Denote by $\cE_{D\times X/X}(0)$ the subsheaf of $\cE_{D\times X/X}$ having sections $\sum_{j\geq0}a_j(t,x)\theta^j$, and put $\cM_0^\mu=\image\big[\cE_{D\times X/X}(0)\otimes_{\cO_{D\times X}}\cM_0\rightarrow\cM^\mu\big]$. Then $\cM_0^\mu$ is $\cE_{D\times X/X}(0)$-coherent and supported on $\Sigma$, therefore $p_*\cM_0^\mu$ is $\cO_X\lcr\theta\rcr$-coherent, by the Preparation Theorem for microdifferential operators. Notice that, by definition, $p_*\cM^\mu=\cO_X\lcr\theta\rcr[\theta^{-1}]\otimes_{\cO_X\lcr\theta\rcr}p_*\cM_0^\mu$.

\begin{proposition}\label{prop:muG0}
Assume that $p_*\cM_0^\mu$ is $\cO_X\lcr\theta\rcr$-locally free of rank equal to $\rg G$. Then,
\begin{enumerate}
\item\label{prop:muG01}
$G$ is $\cO_X[\theta,\theta^{-1}]$-locally free (of rank $\rg G$);
\item\label{prop:muG02}
there exists a unique $\cO_X[\theta]$-locally free submodule $G_0\subset G$ of rank $\rg G$ such that $G=\cO_X[\theta,\theta^{-1}]\otimes_{\cO_X[\theta]}G_0$ and, under the isomorphism of Proposition \ref{prop:muG}\eqref{prop:muG2}, we have $\cO_X\lcr\theta\rcr\otimes_{\cO_X[\theta]}G_0\isom p_*\cM_0^\mu$.
\end{enumerate}
Moreover, the construction of $(G,G_0)$ is compatible with base change.
\end{proposition}

\begin{proof}
We know, by Corollary~\ref{cor:Fourierreg}, that $G$ is $\cO_X[\theta,\theta^{-1}]$-coherent. Denote by $\CC$ the chart of $\wh\PP^1$ with coordinate $\theta$, and consider $G^\an$ on $\CC\times X$, which is $\cO_{\CC\times X}[\theta^{-1}]$-coherent, and $\cO_{\CC^*\times X}$-locally free of rank $\rg G$ when restricted to $\theta\neq0$, as it has a flat connection (see also Theorem~\ref{th:Fourierreg}\eqref{th:Fourierreg2}).

It follows from \cite[Prop\ptbl1.2]{Malgrange95} that there is a bijective correspondence between lattices of $G^\an$ and lattices of $G^\wedge$, where a lattice means here a coherent $\cO_{\CC\times X}$-module (\resp a coherent $\cO_X\lcr\theta\rcr$-module) which generates $G^\an$ or $G^\wedge$ when inverting $\theta$. Moreover, in the analytic case, such a lattice coincides with $G^\an$ out of $\theta=0$.

As we have seen above, $p_*\cM_0^\mu$ is a lattice of $p_*\cM^\mu$. Hence, from Proposition \ref{prop:muG}\eqref{prop:muG2}, we obtain a lattice $G_0^\wedge$ of $G^\wedge$, and consequently a lattice $G_0^\an$ of $G^\an$, such that $\cO_X\lcr\theta\rcr\otimes_{\cO_{\CC\times X}}G_0^\an=G_0^\wedge$.

As $G$ has regular singularities along $\theta=\infty$, we have $G=\wh p_*\wt{G^\an}$. We therefore \emph{define} $G_0$ as $\wh p_*\wt{G_0^\an}$. In particular, $G_0[\theta^{-1}]=G$.

Assume now that $p_*\cM_0^\mu$ is $\cO_X\lcr\theta\rcr$-locally free of rank $\rg G$. Then $G_0^\an$ is also $\cO_{\CC\times X}$-locally free of rank $\rg G$. Cover $X$ by open sets $U$ for which there exists a disc $D'$ centered at $\theta=0$ such that $G_0^\an$ is $\cO_{D'\times U}$ when restricted to $D'\times U$. Apply then Corollary~\ref{cor:locfreeextension} to obtain that $G_0$ is $\cO_X[\theta]$-locally free and that $G$ is $\cO_X[\theta,\theta^{-1}]$-locally on any such $U$.

Let us now consider the base change. We will assume that $i:X'\hto X$ is a closed inclusion. The case of a projection is easy. Notice first that, even without any assumption of local freeness, the functor $\ccF$ (\cf Definition~\ref{def:cF}) is compatible with base change $i^*$ (also denoted by $i^{+0}$ in order to take account of the $\partial_t$ or $\partial_\tau$ action). Consider now the composite mapping
\[
\cE_{D\times X/X}(0)\ootimes_{\cO_{D\times X}}\cM_0\to\hspace*{-5mm}\to\cM_0^\mu\hto\cE_{D\times X/X}\ootimes_{\cD_{D\times X/X}}\cM.
\]
As $\cE_{D\times X/X}(0)$ is $\cO_{D\times X}$-flat and $i^*$ is right exact, we have a surjective morphism $i^*(\cM_0^\mu)\to(i^*\cM_0)^\mu$ and, as $p$ is finite on $\Sigma$, a surjective morphism $i^*(p_*\cM_0^\mu)=p_*i^*(\cM_0^\mu)\to p_*(i^*\cM_0)^\mu$. This morphism becomes an isomorphism after inverting $\theta$ and, by assumption, the left-hand term is locally free of rank $\rg G$ over $\cO_{X'}\lcr\theta\rcr$. It is therefore an isomorphism.
\end{proof}

\begin{Remarques}\label{rem:mulattice}
\begin{enumerate}
\item\label{rem:mulattice1}
If $X$ is reduced to a point, then one may also construct $G_0$ from $M_0$ in the following way. Denote by $M'_0$ the image of $M_0$ in $M[\partial_t^{-1}]$ by the natural morphism $M\to M[\partial_t^{-1}]=G$. Then $G_0$ is the $\CC[\theta]$-submodule of $G$ generated by $M'_0$ (indeed, this submodule satisfies the required properties for $G_0$, see \cite[Prop\ptbl2.1 and 2.3]{Sabbah96b}).
\item\label{rem:mulattice2}
If $\cM_0$ is a lattice in $\cM$, one defines microlocal lattices $\cM^\mu_{0,(c,x^o;\eta_i)}$, and \eqref{eq:mumu} induces an isomorphism of the corresponding lattices.
\item\label{rem:mulattice3}
By the very definition of $G_0$, we have $G_0/\theta G_0=p_*\cM_0^\mu/\theta p_*\cM_0^\mu$.
\end{enumerate}
\end{Remarques}

\Subsection{Behaviour with respect to duality}

\begin{proposition}
The functor $\ccF$ (\cf Definition~\ref{def:cF}) commutes with the duality functor of each category, up to conjugation, \ie $\ccF(\bbD\cM)=\ov{\ccF(\cM)}^*$.
\end{proposition}

\begin{proof}
One considers each individual functor entering in $\ccF$. For $\cM\mto\wt\cM$, the commutation with duality has been seen in Proposition~\ref{prop:extension}. For the algebraization $p_*$, this is stated in Corollary~\ref{cor:equivcat}\eqref{cor:equivcat4}. For the Fourier transform, this is \eqref{eq:dualFourier}. For the localization and the equivalence ``localized $\cD$-modules $\iff$ meromorphic connections'', apply first Lemma~\ref{lem:dualconnexions} to $\cG$ defined in Remark~\ref{rem:locfree}, and then use that $G=\wh p_*\wt\cG$.
\end{proof}

Recall that, if $G$ is any coherent $\cO_X[\tau,\tau^{-1}]$-module, we denote by $\ov G$ the $\cO_X$-module $G$ with an action $\cdot$ of $\CC[\tau,\tau^{-1}]$ defined by $\tau\cdot g=-\tau g$. If we denote by $\ov g$ the element $g$ when viewed in $\ov G$, we will write $\tau\ov g=\ov{-\tau g}=-\ov{\tau g}$. If $G$ is equipped with a connection $\nabla$, then $\ov G$ is equipped with the connection $\ov\nabla$ such that $\ov\nabla_{\partial_\tau}\ov g=-\ov{\nabla_{\partial_\tau} g}$ and $\ov\nabla_\xi g=\ov{\nabla_\xi g}$ if $\xi$ is a vector field on $X$ (this is compatible with the notation used in \eqref{eq:dualFourier}). Notice that $\tau\ov\nabla_{\partial_\tau}\ov g=\ov{\tau\nabla_{\partial_\tau}g}$.

A \emph{sesquilinear pairing} $S$ on $G$ is a $\cO_X[\tau,\tau^{-1}]$-linear morphism
$$S:G\ootimes_{\cO_X[\tau,\tau^{-1}]}\ov G\to\cO_X[\tau,\tau^{-1}]$$
compatible with the connections. We say that $S$ is \emph{nondegenerate} if it induces an isomorphism (compatible with the connections)
$\ov G^*\isom G$, with $G^*\defin\cHom_{\cO_X[\tau,\tau^{-1}]}(G,\cO_X[\tau,\tau^{-1}])$.
We say that $S$ is \emph{$w$-Hermitian} ($w\in\ZZ$) if $S(g'',\ov g')=(-1)^w\ov{S(g',\ov g'')}$.

\begin{corollaire}\label{cor:cP}
Let $\cM$ be an object of $\cHolreg_{\Sigma,\NC}(\cD_{D\times X})$. If $\ccP:\bbD\cM\to \cM$ is a morphism, then the morphism $\ccF(\ccP):\ov{G}^*\to G$ induces a sesquilinear pairing $S$ on $G=\ccF(\cM)$. If $\ker\ccP$ and $\coker\ccP$ are $\cO_{D\times X}$-locally free of finite rank, then $S$ is nondegenerate. If the adjoint $\bbD\ccP:\bbD\cM\to \bbD\bbD\cM=\cM$ is equal to $(-1)^w\ccP$, then $S$ is $w$-Hermitian.\hfill\qed
\end{corollaire}

Similarly, for any $(c,x^o)\in \Sigma$, $\ccP$ induces a morphism $\ccP^\mu_{(c,x^o)}:\ov{\bbD\cM^\mu_{(c,x^o)}}\to\cM^\mu_{(c,x^o)}$, where the duality is taken as $\cDXDth_{(c,x^o)}$-modules (see \eqref{eq:cDXDth}): indeed, by flatness of $\cDXDth_{(c,x^o)}$ over $\cD_{D\times X,(c,x^o)}$, and using that left-right transformations commute only up to conjugation, we have a canonical isomorphism $\ov{\bbD\cM^\mu_{(c,x^o)}}\isom(\bbD\cM)^\mu_{(c,x^o)}$. Using Remark~\ref{rem:dualconnexionsfiltr}\eqref{rem:dualconnexionsfiltr1}, we get a sesquilinear pairing
\[
S^\mu_{(c,x^o)}:\cM^\mu_{(c,x^o)}\ootimes_{\cO_{X,x^o}\lcr\theta\rcr[\theta^{-1}]}\ov{\cM^\mu_{(c,x^o)}}\to \cO_{X,x^o}\lcr\theta\rcr[\theta^{-1}].
\]
The analogue of Corollary~\ref{cor:cP} holds for $S^\mu$.

\medskip
Let $G_0$ be a coherent $\cO_X[\tau^{-1}]$-submodule of $G$ such that
$$G=\cO_X[\tau,\tau^{-1}]\ootimes_{\cO_X[\tau^{-1}]}G_0.$$
Consider similarly $\ov G_0\defin\ov{G_0}\subset \ov G$ and $G_0^*\subset G^*$ (we identify $G_0^*$ with the sheaf of homomorphisms $\lambda:G\to \cO_X[\tau,\tau^{-1}]$ such that $\lambda(G_0)\subset\cO_X[\tau^{-1}]$).

\begin{definition}\label{def:Hemitianw}
The nondegenerate sesquilinear pairing $S$ is said to be \emph{Hermitian of weight $w$} on $(G,G_0)$ if $S$ is $w$-Hermitian and if the isomorphism $\ov G^*\to G$ induced by $S$ sends $\ov G^*_0$ onto $\tau^{w}G_0$, in other words if $S$ induces a perfect pairing
$$G_0\ootimes_{\cO_X[\tau^{-1}]}\ov G_0\to\tau^{-w}\cO_X[\tau^{-1}].$$
\end{definition}

\begin{proposition}[A microdifferential criterion]\label{prop:microcrit}
Let $\cM$ be an object of $\cHolreg_{\Sigma,\NC}(\cDXD)$ with a lattice $\cM_0$. Let $\ccP:\bbD\cM\to\cM$. Assume that
\begin{enumerate}
\item\label{prop:microcrit1}
$\cM_0$ has Poincar\'e rank one (\cf Definition~\ref{def:lattices}),
\item\label{prop:microcrit2}
$\ccP$ satisfies all assumptions of Corollary~\ref{cor:cP},
\item\label{prop:microcrit3}
for some $x^o\in X$, $p_*\cM^\mu_{0,x^o}$ is $\cO_{X,x^o}\lcr\theta\rcr$-free of rank $\rg G$,
\item\label{prop:microcrit4}
$p_*S^\mu_{x^o}\big(\ov{p_*\cM^\mu_{0,x^o}}^*\big)=\theta^{-w}p_*\cM^\mu_{0,x^o}$ (\ie $p_*S^\mu_{x^o}$ is Hermitian of weight $w$ on $(\cM^\mu_{x^o},\cM^\mu_{0,x^o})$).
\end{enumerate}
Then $S$ is Hermitian of weight $w$ on $(G,G_0)$ in some neighbourhood of $x^o$.
\end{proposition}

\begin{proof}
It will have three steps.
\begin{enumerate}
\item
Consider the filtration $F_\bbullet$ by the order of differential operators on $\cDXD$ and extend it on $\cDXDth$, so that $F_{-1}=0$ and $F_0=\cE_{D\times X/X}(0)$ (this is the filtration considered in Remark~\ref{rem:dualconnexionsfiltr}\eqref{rem:dualconnexionsfiltr1}). The lattice $\cM_0$ induces a good filtration $F_k\cM=F_k\cDXD\cdot\cM_0$. Similarly, $p_*\cM^\mu_0$ induces a good filtration of $p_*\cM^\mu$ as a $\cDXDth$-module. At the level of Rees modules (\cf\T\ref{subsec:dualmerom}), we have
\begin{equation}\label{eq:RF}
R_Fp_*\cM^\mu=R_F\cDXDth\ootimes_{R_F\cDXD}R_F\cM.
\end{equation}

As $\cM_0$ has Poincar\'e rank one, we have $F_kp_*\cM^\mu=\theta^{-k}p_*\cM^\mu_0$ for any \hbox{$k\geq0$}. By Assumption~\ref{prop:microcrit}\eqref{prop:microcrit3}, $R_Fp_*\cM^\mu$ is free as a $R_F\cO_{X,x^o}\lcr\theta\rcr[\theta^{-1}]$-module. Therefore, we may apply the analogue of Lemma~\ref{lem:dualconnexionsfiltr} to conclude that the dual complex of $R_Fp_*\cM^\mu_{x^o}$ has cohomology in degree $\dim X+1$ only, and that the filtered module $\bbD p_*\cM^\mu_{x^o}$ is identified with $(p_*\cM^\mu_{x^o})^*$ with the filtration dual to $F_\bbullet p_*\cM^\mu_{x^o}$. In particular, $F_0\bbD p_*\cM^\mu_{x^o}= (p_*\cM^\mu_{0,x^o})^*$.

Assumption~\ref{prop:microcrit}\eqref{prop:microcrit4} is therefore equivalent to the fact that $\ccP^\mu_{x^o}$ strictly shifts the filtrations by $w$ (up to conjugation).

Using flatness of $R_F\cDXDth$ over $R_F\cDXD$ and \eqref{eq:RF}, identify the filtered module $\ov{\bbD p_*\cM^\mu}$ with the filtered module $\cDXDth\otimes_{\cDXD}\bbD\cM$. Then, Assumption~\ref{prop:microcrit}\eqref{prop:microcrit4} is also equivalent to the fact that $\id\otimes\ccP_{x^o}$ strictly shifts the filtrations by $w$.

\item
By Assumption~\ref{prop:microcrit}\eqref{prop:microcrit2}, $S$ is Hermitian of weight $w$ on $(G,G_0)$ if and only if $S^\wedge$ is so on $(G^\wedge,G^\wedge_0)$.

Denote by $\cD$ the sheaf $\cD_X[\theta,\theta^{-1}]\langle t\rangle$, with $\theta^{-1}=\partial_t$, and by $\cD^\wedge$ its tensor product by $\cO_X\lcr\theta\rcr[\theta^{-1}]$. We therefore have $G=\cD\otimes_{\cDXt}M$ and $G^\wedge=\cD^\wedge\otimes_{\cDXt}M$. Let $F_\bbullet\cDXt$ be the filtration by the degree of differential operators, denote by $F_\bbullet\cD,F_\bbullet\cD^\wedge$ the similar filtrations.

We may argue as in the first step, starting from $(M,M_0)$: indeed, $M_0$ also has Poincar\'e rank one and, by Assumption~\ref{prop:microcrit}\eqref{prop:microcrit3} and Proposition~\ref{prop:muG0}, $G^\wedge_{0,x^o}$ is also free. One also uses Remark~\ref{rem:dualconnexionsfiltr}\eqref{rem:dualconnexionsfiltr2}.

Hence, $S$ is Hermitian of weight $w$ on $(G,G_0)$ near $x^o$ if and only if $\id\otimes\ccP:\cD^\wedge\otimes_{\cDXt}\bbD M\to\cD^\wedge\otimes_{\cDXt}M$ strictly shifts the filtrations by $w$.

\item
Last, we conclude by noticing that the analytization morphism $\cD^\wedge\otimes_{\cDXt}\cbbullet\to p_*(\cDXDth\otimes_{\cDXt}\cbbullet)$ is strict with respect to the filtrations.\qedhere
\end{enumerate}
\end{proof}

\section{Gauss-Manin system and Brieskorn lattice}\label{sec:GMB}

\subsection{The setting}\label{subsec:setting}
Let $\AfN$ be the affine space with coordinates $u_1,\dots,u_N$. Let $X$ be a complex manifold and let $\cU\subset\AfN\times X$ be the closed subset defined by an ideal in $\cO_X[u_1,\dots,u_N]$. We assume that the projection $q:\cU^{\an}\to X$ is smooth of relative dimension $n$, with $n\geq0$.

Let $F$ be a section of $\cO_X[u_1,\dots,u_N]$. It defines a function $F:\cU\to \CC$. We put $\Phi=(F,q):\cU^{\an}\to \CC\times X$.

Let $\delta_a$ denote the distance function to some point $a\in\AfN$. For $x\in X$, put $\cU_x=q^{-1}(x)\subset\AfN$. If $x^o\in X$ is fixed, put $\cU^o=q^{-1}(x^o)\subset\AfN$ and $f=F_{\cU^o}:\cU^o\to\CC$. We will assume that $f$ has only \emph{isolated critical points}. We will denote by $\mu=\mu(f)$ the sum of the Milnor numbers $\mu(f,u)$ of $f$ at its critical points $u$.

Following \cite{N-Z90}, we say that $f$ is \emph{M-tame} if, for some choice of $a\in\AfN$, for any $\eta>0$ there exists $R(\eta)>0$ such that, for any $r\geq R(\eta)$, the spheres $\delta_a(x)=r$ are transversal to $f^{-1}(t)$ for any $t$ with $\module{t}\leq\eta$.

Denote by $D(0,\eta)$ the open disc centered at $0$ in $\CC$ and $B(a,R)$ the open ball centered at $a$ in $\AfN$.

If $f$ is M-tame, there exist $\eta>0$, $R>0$, $\epsilon>0$ and a neighbourhood of $x^o$ in $X$, that we still denote by $X$, such that
\begin{enumerate}
\item
all critical points of $f$ are contained in $\cU^o\cap B(a,R-\epsilon)\cap f^{-1}(D(0,\eta))$,
\item
the fibers $\Phi^{-1}(t,x)\subset\AfN$ are transversal to the spheres $\delta(x)=r$ for any $x\in X$, $t\in D(0,\eta)$ and $r\in{}]R-\epsilon,R+\epsilon[$,
\item
moreover, the fibers $f^{-1}(t)=\Phi^{-1}(t,x^o)$ are transversal to spheres $\delta(x)=r$ for any $r>R-\epsilon$ and $t\in D(0,\eta)$.
\end{enumerate}

\begin{figure}[htb]
\centerline{\includegraphics[scale=.9]{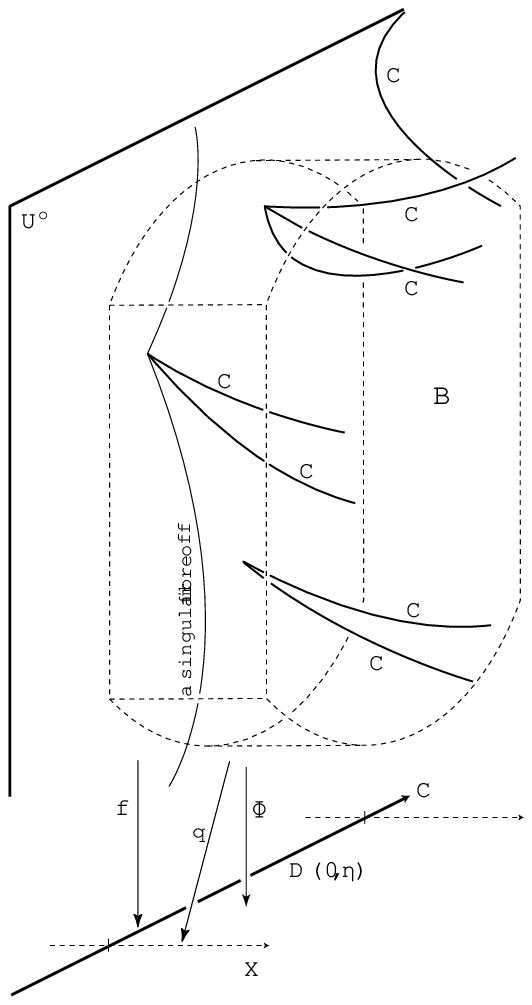}}
\caption{}\label{fig:1}
\end{figure}

In the following, we fix $a,\eta,R, \epsilon$ and $X$ as above. We denote by $\cB$ the open set $\cU\cap B(a,R)\cap F^{-1}(D(0,\eta))$ and we put $D=D(0,\eta)$. We therefore have the following diagram:
\[
\xymatrix@C=1.5cm{
&\cB\ar[ldd]_-{F}\ar[rdd]^-q\ar[d]^-{\Phi}\ar@{^{ (}->}[r]^-{i_F}& D\times\cB\ar[ld]^(.3){\;\id\times q}|!{[l];[dd]}\hole\\
&D\times X\ar[rd]_p\ar[ld]\\
D&&X
}
\]

Let $C$ denote the critical set of~$\Phi_{|\cB}$ (see fig\ptbl\ref{fig:1}). Then $\Phi_{|C}:C\to\Afu\times X$ is finite onto its image and $q_{|C}:C\to X$ is finite and onto (at least if $C$ is not empty). We also denote by $C'$ the other components of the critical set of $\Phi$ in $\cU$, on which $\Phi$ may not be finite. Notice that $C\cup C'$ is locally a complete intersection, being defined by the vanishing of the partial derivatives of $F$ with respect to ``vertical'' local coordinates on $\cU^{\an}$. As $q_{|C}$ is finite, we get:

\begin{lemme}\label{lem:qClocfree}
The sheaf $q_*\cO_C$ is $\cO_X$-locally free of rank $\mu$.\hfill\qed
\end{lemme}

\begin{remarque}\label{rem:omegadF}
Denote by $\Omega^\cbbullet_{\cB/X}$ the sheaves of holomorphic differential forms relative to the smooth morphism $q:\cB\to X$. Then the complex $(\Omega^\cbbullet_{\cB/X},d_{\cB/X}F\wedge)$ has cohomology in degree $n$ only. Moreover, $\Omega^{n}_{\cB/X}\big/(d_{\cB/X}F\wedge\Omega^{n-1}_{\cB/X})$ is supported on $C$ and, locally on $X$, is isomorphic to $\cO_C$ (by choosing a local generator of $\Omega^{n}_{\cB/X}$ near each point of \hbox{$C\cap q^{-1}(x)$}). Hence, $q_*\big[\Omega^{n}_{\cB/X}\big/(d_{\cB/X}F\wedge\Omega^{n-1}_{\cB/X})\big]$ is $q_*\cO_C$-locally free of rank one. Notice also that the complex $(\Phi_*(\Omega^\cbbullet_{\cB/X}),d_{\cB/X}F\wedge)$ has cohomology in degree $n$ only, as $\Phi$ is Stein.
\end{remarque}

Let $\xi$ be a local section of the tangent sheaf $\Theta_X$, \ie a local holomorphic vector field on $X$. Choose a local lifting $\xi'$ of $\xi$ as a vector field on $\cU^{\an}$ and consider the class of $\cL_{\xi'}(F)$ in $\cO_C$. By definition of $C$, this class does not depend on the choice of the lifting, because $\cL_\eta(F)$ belongs to the ideal of $C$ if $\eta$ is ``vertical''. We therefore denote this class by $[\cL_\xi(F)]$.

\begin{definition}\label{def:KS}
The \emph{Kodaira-Spencer map} $\varphi$ attached to $\Phi$ is the $\cO_X$-linear morphism
\[
\Theta_X\to q_*\cO_C,\quad
\xi\mto [\cL_\xi(F)].
\]
\end{definition}

\begin{Remarques}
\begin{enumerate}
\item
Usually, the critical locus of $\Phi$ on $\cU$ contains components other than $C$, on which $\Phi$ may not be finite. In other words, the function $F(\cdot,x)$ on the smooth affine variety $\cU_x$ may have critical points which disappear at infinity when $x\to x^o$. By the choice of $\cB$, we do not care about these critical points. On the other hand, $\cB$ is big enough so that restricting to $\cB$ looses no information concerning $f$ on $\cU^o$.

For the same reason, we \emph{do not} assert that, for $\module{t}<\eta$, the fibres $\{F(\cdot,x)=t\}\subset \cU_x$ and $\{F(\cdot,x)=t\}\cap B(a,R)$ have the same topological type. This only holds \emph{a priori} for the fibres $f=t$ on $\cU^o$.
\item
We may replace the affine manifold $\cU$ with a Stein closed submanifold $\cU^{\an}$ of $B(a,R+\epsilon)\times X$. In such a way, one recovers for instance the classical ``local'' situation of an isolated critical point of a holomorphic function.
\end{enumerate}
\end{Remarques}

\begin{exemples}
In the examples below, we consider a regular function $f:U\to\CC$ on an affine manifold on which global coordinates exist ($U=\CC^n$ or $U=(\CC^*)^n$) and we consider a one-parameter family $F(\cdot,x)=f+xg$ on $\cU=U\times X$, where the class of $g$ in the Jacobian quotient $\cO(U)/(\partial f)$ is nonzero.
\begin{enumerate}
\item
Consider the following situation (\cf \cite[\T3.3.7]{Sabbah96c}): $\cU=\mathbb{A}^2\times X$, $X$ is a disc with coordinate $x$ and $F(u_0,u_1,x)=u_0^5+u_1^5+xu_0^3u_1^3$. Then, for any $x$, $F(\cdot,x)$ has a critical point at $u_0=u_1=0$ with Milnor number equal to $(5-1)^2=16$. For any $x\neq0$, $F(\cdot,x)$ has five other critical points $\frac{5}{3x}(\zeta,1/\zeta)$ with $\zeta^5=-1$. These points disappear at infinity when $x\to 0$. The first Betti number of the smooth fiber of $f=F(\cdot,0)$ is $16$, but that of the smooth fiber of $F(\cdot,x)$ for $x\neq0$ is $21$.

\item
Let $U=(\CC^*)^2$ and $\cU=U\times X$. Consider the function $f(u_1,u_2)=u_1+u_2+1/(u_1u_2)$ (\cf \cite{Barannikov00}) and the perturbation $F(u_1,u_2,x)=f(u_1,u_2)+x/(u_1u_2)^2$. Then $f$ has $3$ critical points and, for $x\neq0$, $F(\cdot, x)$ has $5$ critical points, two of them approaching $u_1=u_2=0$ when $x\to 0$. The corresponding critical values go to infinity. A simple Euler characteristic computation shows that the Betti number of the smooth fibre of $F(\cdot, x)$ is strictly bigger than that of $f$ if $x\neq0$.

\item
In both previous examples, the critical locus of $\Phi$ contains components other than $C$, but $\Phi$ remains finite on all the components. One may vary a little bit the previous example to obtain an example where $\Phi$ is not finite on these other components: take $F=f+xf^2$. Notice that, up to a constant, $f^2$ and $1/(u_1u_2)^2$ have the same class in the Jacobian quotient $\CC[u_1,u_1^{-1},u_2,u_2^{-1}]/(u_1\partial f/\partial u_1,u_2\partial f/\partial u_2)$. The conclusion concerning the Betti numbers remains the same.
\end{enumerate}
\end{exemples}

\subsection{The Gauss-Manin system}
We collect here various known results concerning the Gauss-Manin system (\cf \cite{Brieskorn70,Pham79,Pham81,Pham85,S-S85,Oda87,MSaito89}). These usually apply to a local situation around each critical point of $f$. We will show how to adapt them to the global situation considered here. The results of \cite{Sch-Sch94} will be useful. The presentation follows that of \cite{MSaito89}.

The \emph{Gauss-Manin complex} is the direct image complex $\Phi_+\cO_{\cB}$ in $D^b(\cD_{D\times X})$. Choosing local coordinates $(x_1,\dots,x_m)$ on $X$, so that $q=(q_1,\dots,q_m)$, with associated vector fields $\partial_{x_1},\dots,\partial_{x_m}$, we may express it as
\[
\Phi_+\cO_{\cB}=\bR\Phi_*\Big(\Omega^{n+m+\bbullet}_{\cB}[\partial_t,\partial_{x_1},\dots,\partial_{x_m}],\nabla\Big),
\]
with $\nabla$ defined as
\[
\nabla(\eta\otimes P)=d\eta\otimes P-\sum_j dq_j\wedge\eta\otimes\partial_{x_j}P-dt\wedge\eta\otimes\partial_t P.
\]
The $\cD_{D\times X}$-structure on the complex is defined by
\[
\partial_t(\eta\otimes P)=\eta\otimes(\partial_tP),\quad \partial_{x_i}(\eta\otimes P)=\eta\otimes(\partial_{x_i}P),\quad g(t,x)\cdot (\eta\otimes 1)=g(t,x)\eta\otimes 1.
\]

The Gauss-Manin system is the $\cD_{D\times X}$-module $\cM=\cH^0(\Phi_+\cO_{\cB})$. As we will see below, when the conditions of \T\ref{subsec:setting} are satisfied, the other cohomologies of $\Phi_+\cO_{\cB}$ are $\cO_{D\times X}$-locally free of finite rank, the rank being computed in terms of the Betti numbers of $\cB$.

Consider the filtered complex $F_k\Big(\Omega^{n+m+\bbullet}_{\cB}[\partial_t,\partial_{x_1},\dots,\partial_{x_m}],\nabla\Big)$, with
\begin{equation}\label{eq:filtr}
F_k\Omega^{n+m+\ell}_{\cB}[\partial_t,\partial_{x_1},\dots,\partial_{x_m}]= \sum_{\module{\alpha}\leq k+\ell}\Omega^{n+m+\ell}_{\cB}\partial^\alpha
\end{equation}
and $\partial^\alpha=\partial_t^{\alpha_0}\partial_{x_1}^{\alpha_1}\cdots\partial_{x_m}^{\alpha_m}$.

Recall (see \cite[Lemma 2.2]{MSaito89}) that, taking the realization of $\bR\Phi_*$ using the Godement canonical flabby resolution, we have
\[
\Phi_+{\cO_\cB}=\varinjlim_k\bR\Phi_*F_k\Big(\Omega^{n+m+\bbullet}_{\cB}[\partial_t,\partial_{x_1},\dots,\partial_{x_m}],\nabla\Big).
\]
As $\Phi$ is Stein and as each term in $F_k$ is $\cO_\cB$-coherent, it follows that
$$
\Phi_+\cO_{\cB}=\Big(\Phi_*\Omega^{n+m+\bbullet}_{\cB}[\partial_t,\partial_{x_1},\dots,\partial_{x_m}],\nabla\Big).
$$

If we forget the $\cD_X$-structure on $\Phi_+\cO_\cB$, we may compute this complex as a relative de~Rham complex, \ie as the direct image of $i_+\cO_\cB$ viewed as a $\cD_{D\times\cB/X}$-module, if $i:\cB\hto D\times\cB$ denotes the graph inclusion of $F$ (see \cite[Lemma 2.4]{MSaito89}):

\begin{proposition}\label{prop:GM}
For any $j\in\ZZ$, $\cH^j\Phi_+\cO_{\cB}$ is isomorphic (as a $\cD_{D\times X/X}$-module) to the $j$-th cohomology of the complex
\[
\Big(\Phi_*\Omega^{n+\bbullet}_{\cB/X}[\partial_t], d_{\cB/X}-\partial_td_{\cB/X}F\wedge\Big),
\]
where the $\cO_{D\times X}\langle\partial_t\rangle$-action is defined by the usual formulas:
$$\partial_t\cdot[\eta\partial_t^k]=[\eta\partial_t^{k+1}],\quad g(t,x)\cdot[\eta\partial_t]=[g(F,x)\eta\partial_t]-[g'_t(F,x)\eta].
$$
In particular, $\cH^j\Phi_+\cO_{\cB}=0$ for $j>0$.\hfill\qed
\end{proposition}

\begin{remarque}\label{rem:GMfiltr}
One may filter the previous complex using a formula analogous to \eqref{eq:filtr}. Then the previous proposition holds true with filtration (\cf \loccit.).
\end{remarque}

The main finiteness result for the Gauss-Manin system of M-tame functions is a direct application of general results for ``elliptic pairs'' proved in \cite{Sch-Sch94}.

\begin{theoreme}\label{th:GM}
The Gauss-Manin complex satisfies the following properties:
\begin{enumerate}
\item\label{th:GM1}
The cohomology modules $\cH^j(\Phi_+\cO_{\cB})$ are coherent, holonomic and regular $\cD_{D\times X}$-modules.
\item\label{th:GM2}
$\cH^j(\Phi_+\cO_{\cB})=0$ for $j\geq1$ and $\cH^j(\Phi_+\cO_{\cB})$ is a locally free $\cO_{D\times X}$-module of rank $\dim H^{j+n}(\cU^o)$ if $j<0$.
\item\label{th:GM3}
The direct image $\cM^o\defin \cH^0f_+\cO_{\cB^o}$ is equal to the restriction to $D\times \{x^o\}$ of $M^{\an}\defin (\cH^0f_+\cO_{\cU^o})^{\an}$.
\item\label{th:GM4}
The fibres of the projection $p:D\times X\to X$ are noncharacteristic with respect to the Gauss-Manin complex $\Phi_+\cO_{\cB}$ and we have $i_{x^o}^+\cM=\cH^0(f_+\cO_{\cB^o})$.
\item\label{th:GM5}
The Poincar\'e duality morphism induces a $\cD_{D\times X}$-linear morphism $\bbD\cM\to \cM$, the kernel and cokernel of which are $\cO_{D\times X}$-locally free of finite rank.
\item\label{th:GM6}
The filtration of $\cM$ induced by \eqref{eq:filtr} is good.
\end{enumerate}
\end{theoreme}

\begin{proof}
We apply \cite[Cor\ptbl8.1]{Sch-Sch94} to
\begin{itemize}
\item
the morphism $\Phi:\Phi^{-1}(D(0,\eta)\times X)\to D(0,\eta)\times X$, that we restrict to the family of open sets $B(a,r)\times X$ with $r\in {}]R-\epsilon,R+\epsilon[$,
\item
the constant sheaf restricted to $\cU\cap (B(a,R+\epsilon)\times X)$ and
\item
the $\cD$-module $\cO_{\cU}$ restricted to this open set.
\end{itemize}
We conclude that $\Phi_+\cO_{\cB}$ has coherent $\cD_{D\times X}$-cohomology and that we have a morphism
\[
\bbD\Phi_+\cO_{\cB}\isom \Phi_\dag\cO_{\cB}\to \Phi_+\cO_{\cB}
\]
in $D^b(\cD_{D\times X})$, where $\Phi_\dag$ denotes the direct image of $\cD$-modules with proper support (we have used that $\cO_{\cB}$ is selfdual). Using \cite[Th\ptbl7.5]{Sch-Sch94}, we get the usual Kashiwara's estimate for the characteristic variety of the direct images $\Phi_+$ or $\Phi_\dag$, showing that their cohomologies are holonomic $\cD_{D\times X}$-modules. Moreover, the same argument as in \cite[Cor\ptbl7.6]{Sch-Sch94} shows the $\cO_{D\times X}$-local freeness of $\cH^j(\Phi_+\cO_{\cB})$ for $j\neq0$ (recall that $\cH^j(\Phi_+\cO_{\cB})=0$ for $j>0$, \cf Proposition \ref{prop:GM}). That the cohomology of the cone of $\Phi_\dag\cO_{\cB}\to \Phi_+\cO_{\cB}$ is $\cO_{D\times X}$-locally free may be seen in the same way.

Let us show the regularity of $\Phi_+\cO_{\cB}$ (by the previous arguments, we would only need to show the regularity of $\cM$). Let $Z$ be a hypersurface in $D\times X$. Notice that we also may apply the results of \cite{Sch-Sch94} to the localized module $\cO_{\cB}(*\Phi^{-1}(Z))$, as it also satisfies the ellipticity condition. We may work with $\Phi_\dag\cO_{\cB}$ instead of $\Phi_+\cO_{\cB}$. Therefore, we may apply the same argument as in \cite{Mebkhout89} for the direct image of the irregularity sheaf: we have $\Phi_\dag\cO_{\cB}(*\Phi^{-1}(Z))\isom (\Phi_\dag\cO_{\cB})(*Z)$; the irregularity sheaf of $\Phi_\dag\cO_{\cB}$ along $Z$ is thus equal to $\bR\Phi_!$ of that of $\cO_{\cB}$, \ie is equal to $0$.

For the noncharacteristic property in \ref{th:GM}\eqref{th:GM4}, see \cite[p\ptbl281]{Pham79}. The ``base change'' assertion is proved as in \cite[VI, 8.4]{Borel87} (\cf also \cite[Prop\ptbl1.6]{D-M-S-S99}).

For \ref{th:GM}\eqref{th:GM3}, let us recall the proof given in \cite[\T1]{N-S97}. Remark first that the previous results also apply (with the same arguments) to $f_+\cO_{f^{-1}(D_r)}^{\an}$ for any disc $D_r$, hence to $f_+\cO_{\cU^o}^{\an}$. They also show that the restriction morphism
\[
(f_+\cO_{\cU^o}^{\an})_{|D}=f_+\cO_{f^{-1}(D)}^{\an}\to f_+\cO_{\cB^o}^{\an}
\]
is a quasi-isomorphism. Using a projectivization $\ov f:Z\to\Afu$ of $f$ and denoting by $j:\cU^o\hto Z$ the inclusion, one finds $f_+\cO_{\cU^o}=\ov f_+(j_+\cO_{\cU^o})$, $(f_+\cO_{\cU^o})^\an=\ov f_+(j_+\cO_{\cU^o})^\an$ and $(f_+\cO_{\cU^o}^\an)=\ov f_+(\bR j_*\cO_{\cU^o}^\an)$. The natural morphism $\ov f_+(j_+\cO_{\cU^o})^\an\to\ov f_+(\bR j_*\cO_{\cU^o}^\an)$ is a quasi-isomorphism: indeed, as both $\cD_{\Afuan}$-complexes have regular holonomic cohomology, it is enough to prove this after applying the de~Rham functor; this functor commutes with $\ov f_+$ and $\bR\ov f_*$ (see, \eg \cite[II.5.5]{M-N90}); last, the natural morphism $\DR(j_+\cO_{\cU^o})^\an\to \DR(\bR j_*\cO_{\cU^o}^\an)$ is a quasi-isomorphism, by the Grothendieck comparison theorem.

Let us end with the proof of \ref{th:GM}\eqref{th:GM6}. It is enough to prove the $\cO_{D\times X}$-coherence of the \emph{Brieskorn lattice} $\cM_0$ defined as $F_0\cM$. Indeed, we have
\begin{equation}\label{eq:Br}
\begin{split}
\cM_0=\image\big[\Phi_*\Omega^n_{\cB}\rightarrow\cM\big]&=\image\big[\Phi_*\Omega^n_{\cB/X}\rightarrow\cM\big]\\
&=\Phi_*(\Omega^{n}_{\cB/X})\big/d_{\cB/X}F\wedge \Phi_*(d_{\cB/X}\Omega^{n-2}_{\cB/X}),
\end{split}
\end{equation}
as the complex $(\Omega^{\cbbullet}_{\cB/X},d_{\cB/X}F\wedge)$ has cohomology in degree at most $n$. Now, it is known that the $\cO_{D\times X}$-coherence of $\cM_0$ is reduced to the $\cO_{D\times X}$-coherence of the relative cohomology $\Phi_*\Omega^j_{\cB/D\times X}$, which can be proved as in \cite{B-G80} for instance.
\end{proof}

\begin{Remarques}\label{rem:NCGM}
\begin{enumerate}
\item\label{rem:NCGM1}
We may argue as in \cite[p\ptbl281]{Pham79} to conclude that $\cM$ satisfies the noncharacteristic property \NC.

\item\label{rem:NCGM2}
If $R'$ is a radius $\geq R$, then it also follows from \cite[Cor\ptbl8.1]{Sch-Sch94} that the restriction morphism $\Phi_+\cO_{\cB'}\to\Phi_+{\cB}$ (maybe defined after restriction of $X$ to a smaller neighbourhood $X'$ of $x^o$) is a quasi-isomorphism, which is compatible with the filtrations. In particular, $(\cM',\cM'_0)\to(\cM,\cM_0)$ is an isomorphism on $D\times X'$. Similarly, if one chooses another system of balls on $\cU$, then for sufficiently large radii $R$, the corresponding filtered Gauss-Manin systems are isomorphic to $(\cM,\cM_0)$ on $D\times X''$, if $X''$ is a sufficiently small neighbourhood of $x^o$. Moreover, the Poincar\'e duality morphisms are compatible with these isomorphisms.
\end{enumerate}
\end{Remarques}

\subsection{The Brieskorn lattice (case without parameters)} \label{subsec:algcomp}
In this paragraph, we put $U=\cU^o$ and we omit the exponent~${}^o$, as we do not use parameters (\ie we assume that $X$ is reduce to a point). We assume as above that $U$ is affine. Let $f:U\to \Afu$ be a regular function on $U$. Recall (\cf \cite{Sabbah96b}) that the algebraic Gauss-Manin system and the algebraic Brieskorn lattice are given by the following formulas, where we use algebraic differential forms on the affine manifold~$U$:
\begin{align*}
M&=\Omega^{n}(U)[\partial_t]\big/(d-\partial_tdf\wedge)\Omega^{n-1}(U)[\partial_t],\\
M_0&=\image [\Omega^{n}(U)\rightarrow M]\\[5pt]
G&=\Omega^{n}(U)[\tau,\tau^{-1}]\big/(d-\tau df\wedge) \Omega^{n-1}(U)[\tau,\tau^{-1}]\\
&=\Omega^{n}(U)[\theta,\theta^{-1}]\big/(\theta d-df\wedge) \Omega^{n-1}(U)[\theta,\theta^{-1}]\\
G_0&=\image\big[\Omega^{n}(U)[\theta]\rightarrow G\big].
\end{align*}
That we get such a formula for $G_0$ starting from the definition given in Proposition~\ref{prop:muG0}\eqref{prop:muG02} when $X={\rm pt}$ follows from Remark~\ref{rem:mulattice}\eqref{rem:mulattice1}.

Denote by $f_+$ the algebraic direct image of $\cD_U$-modules, so that $M=\cH^0f_+\cO_U$. Remark that the cone of $f_+\cO_U\MRE{\partial_t}f_+\cO_U$ is the direct image of $\cO_U$ by the constant map. We thus have an exact sequence
\begin{multline*}
\cdots\to\cH^jf_+\cO_U\MRE{\partial_t}\cH^jf_+\cO_U\to H^{j+n}(U,\CC)\to\\
\cdots\to M\MRE{\partial_t}M\to H^{n}(U,\CC)\to0
\end{multline*}
as $\cH^jf_+\cO_U=0$ for $j\not\in[-n+1,0]$. If $f$ is M-tame, or cohomologically tame (\cf \cite{Sabbah96b}), then $\cH^jf_+\cO_U$ is isomorphic to $\cO_{\Afu}^{r_j}$ (for some integer $r_j$) with its usual $\partial_t$-action, if $j<0$. Hence the maps $\partial_t:\cH^jf_+\cO_U\to\cH^jf_+\cO_U$ are onto for $j<0$. Therefore, we get an exact sequence
\[
0\to H^{n-1}(U,\CC)\to M\MRE{\partial_t}M\to H^{n}(U,\CC)\to0.
\]
In particular, $\dim H^{n}(U,\CC)-\dim H^{n-1}(U,\CC)$ is nothing but the Euler characteristics of the algebraic de~Rham complex of $M$. As $M$ is regular (included at infinity), this is also the global Euler characteristics of the analytic de~Rham complex of $M^\an$. This Euler characteristics may be computed as the difference between the generic dimension of the horizontal sections of $M^\an$ (which is equal to the generic dimension $\dim_{\CC(t)}\big(\CC(t)\otimes_{\CC[t]}M\big)$ of $M$) and the sum of dimensions of vanishing cycles of the complex $\DR M^\an$ (which is equal to $\dim_{\CC(\tau)}\big(\CC(\tau)\otimes_{\CC[\tau]}G\big)$, see, \eg \cite[p\ptbl78]{Malgrange91}). One deduces that
\begin{multline*}
\dim_{\CC(t)}\big(\CC(t)\otimes_{\CC[t]}M\big)\\
= \dim_{\CC(\tau)}\big(\CC(\tau)\otimes_{\CC[\tau]}G\big)+\dim H^{n}(U,\CC)-
\dim H^{n-1}(U,\CC).
\end{multline*}
Using that, for any $c\in\CC$, the vector space $\phi_{t-c}\DR M^\an$ is equal to $\bR\Gamma\phi_{f-c}\CC_{U^\an}$ (\cf \cite[Cor\ptbl8.4 and proof of Prop\ptbl9.2]{Sabbah96b} for cohomologically tame functions, and \cite[\T1]{N-S97} for M-tame functions), we also deduce that $\dim_{\CC(\tau)}\CC(\tau)\otimes_{\CC[\tau]}G=\mu$.

That $f$ has only isolated singularities (and assuming $n\geq2$) implies that the complex $(\Omega^{\cbbullet}(U),df\wedge)$ has cohomology in degree $n$ only, hence
\[
M_0=\Omega^n(U)\big/df\wedge d\Omega^{n-2}(U)\quad\text{and}\quad G_0=\Omega^n(U)[\theta]\big/(\theta d-df\wedge)\Omega^{n-1}(U)[\theta]
\]
(see also \cite{Douai02} where this is stated for cohomologically tame polynomials).

Consider $\cM\defin\cH^0f_+\cO_{U^{\an}}$ with its lattice $\cM_0$. M-tameness of $f$ implies, as in Theorem \ref{th:GM}, the regular holonomicity of $\cM$ and the $\cO_{\Afuan}$-coherence of $\cM_0$. We also have $M^{\an}=\cM$ (see \cite[\T1]{N-S97}), and $M_0\subset\Gamma(\Afuan,\cM_0)\cap M$, so $M_0^{\an}\subset\cM_0$ and therefore $M_0=p_*\wt{M_0^\an}\subset p_*\wt\cM_0$. As $p_*\wt\cM_0$ is finitely generated over $\CC[t]$, this implies the finiteness of $M_0$ as a $\CC[t]$-module (see \cite{Sabbah96b} and \cite{N-S97}, see also \cite{D-S01}), hence, by \cite{Sabbah96b}, that of $G_0$ (and its freeness) as a $\CC[\theta]$-module, according to the regularity of $M$.

If $f$ is M-tame, one may also consider the direct image $f_+\cO_\cB$, where $\cB$ is defined in \T\ref{subsec:setting} (taking $X=\text{pt}$), that we now denote by $\cM_\cB$, and therefore define the analytic Brieskorn lattice $\cM_{\cB,0}$ as
the image of $\Omega^n(\cB)$ in $\Omega^n(\cB)[\partial_t]/d_f\Omega^{n-1}(\cB)=\cH^0f_+\cO_\cB$. Arguing as in Theorem \ref{th:GM}, one shows that the natural restriction morphism $\cM_{|D}\to\cM_\cB$ (which sends $\cM_{0|D}$ into $\cM_{\cB,0}$) is an isomorphism.

\begin{lemme}\label{lem:GMalg}
The algebraization of $(\cM_\cB,\cM_{\cB,0})$ is $(M,M_0)$.
\end{lemme}

\begin{proof}[Sketch of proof]
The proof will have 3 steps:
\begin{enumerate}
\item
We will show that $M_0^\an$ can be computed analytically. Let $Z$ be a smooth partial compactification of $U$ so that $f$ extends as a projective morphism $\ov f: Z\to\Afu$ and that $\Delta\defin Z\moins U$ is a divisor [it is not necessary to assume that $Z$ is smooth, it is assumed here to simplify the notation]. We will show that $M_0^{\an}$ is equal to the image of $\Gamma(Z,\Omega^n_{Z^\an}[*\Delta])$ in $\cH^0\ov f_+\cO_{Z^\an}[*\Delta]=\cH^0\big(\Gamma(Z,\Omega^{\cbbullet+n}_{Z^\an}[*\Delta])[\partial_t],d_{\ov f}\big)$. In order to do that, consider a smooth projective compactification $\cZ$ of $Z$ on which $\ov f$ extends as $\ov f:\cZ\to\PP^1$. Denote by $i:\cZ\hto\cZ\times\PP^1$ the graph inclusion of $\ov f$ and put $\cN=i_+\cO_{\cZ}\big(*(\Delta\cup\ov f^{-1}(\infty))\big)$, equipped with a natural lattice $\cN_0$. If $q$ denotes the projection $\cZ\times\PP^1\to\PP^1$, then the $\cO_{\PP^1}(*\infty)$-module $\wt M_0$ associated with $M_0$ is obtained as the image of $q_*(\Omega_{\cZ}^n\otimes_{\cO_\cZ}\cN_0)$ in $\cH^0\big(q_*\big(\Omega_{\cZ}^{\cbbullet}\otimes_{\cO_\cZ}\cN\big)\big)$. By compactness of $\cZ$, $\wt M_0^\an$ is computed using the corresponding analytic objects. Restrict then to $\Afuan$.
\item
We show that the natural morphism $M_0^\an\to\cM_0$ is an isomorphism by a Mayer-Vietoris computation. Fix some big disc $D\subset\Afuan$ and cover $f^{-1}(D)$ by open sets $\cB$, $U^\an\moins\cB'$ for some sets $\cB,\cB'$ as in \T\ref{subsec:setting}, with $\cB'\subset\cB$. One can compute $M_0^\an$ and $\cM_0$ using this covering. One then has to show that the filtered direct image of $\cO_{U^\an}$ on $f^{-1}(D)\moins\cB'$ and that of $\cO_{X^\an}[*\Delta]$ (viewed as a filtered $\cD_{X^\an}[*\Delta]$-module) on $\ov f^{-1}(D)\moins\cB'$ coincide, as the filtered direct images computed on the other open sets clearly coincide. This follows from the regularity along $\Delta$ of $\cO_{X^\an}[*\Delta]$ and from the comparison Theorem (as in the proof of \ref{th:GM}\eqref{th:GM3}): the direct image is $\cO_D$-locally free of finite rank and its filtration is the trivial one.
\item
Last, we show that $\cM_{0|D}\to\cM_{\cB,0}$ is an isomorphism by showing that the corresponding microlocal lattices are isomorphic. By the computation above, this follows from the $\cO_D$-local freeness of the direct images computed on $U^\an\moins\cB'$ and on $\cB\cap (U^\an\moins\cB')$.\qedhere
\end{enumerate}
\end{proof}

\subsection{The Brieskorn lattice}\label{subsec:Br}
We now come back to the original situation, where $X$ is not necessarily reduced to a point.
Recall that, by Proposition \ref{prop:GM}, we have
\[
\cM=\Phi_*(\Omega^{n}_{\cB/X})[\partial_t]\Big/(d_{\cB/X}-\partial_td_{\cB/X}F\wedge)\Phi_*(\Omega^{n-1}_{\cB/X})[\partial_t],
\]
and that $\cM_0$ is defined as the image of $\Phi_*\Omega^{n}_{\cB/X}$ in $\cM$, which can be computed by \eqref{eq:Br}. We know, by \ref{th:GM}\eqref{th:GM6}, that the sheaf $\cM_0$ is $\cO_{D\times X}$-coherent and is a lattice of $\cM$. It follows from Remark~\ref{rem:GMfiltr} that $\cM_0$ has Poincar\'e rank one (\cf Definition~\ref{def:lattices}).

We may associate to $(\cM,\cM_0)$ a pair $(G,G_0)$ as in \T\ref{subsec:FTlattice}. We call $G_0$ the \emph{Brieskorn lattice} associated with $F$.

\begin{proposition}\label{prop:locfree}
The Brieskorn lattice is $\cO_X[\theta]$-locally free of rank $\mu$ and we have $G_0/\theta G_0\simeq q_*\big[\Omega^{n}_{\cB/X}\big/d_{\cB/X}F\wedge\Omega^{n-1}_{\cB/X}\big]$. Moreover, the restriction $G_0^o$ of $G_0$ at $x^o$ is equal to the Brieskorn lattice of $f_{|\cB^o}$, which is nothing but the algebraic Brieskorn lattice of $f:\cU^o\to\Afu$ as defined in \T\ref{subsec:algcomp}.
\end{proposition}

Therefore, by Lemma \ref{lem:qClocfree}, $G_0/\theta G_0$ is $\cO_X$-locally free of rank $\mu$.

\begin{proof}[Proof of Proposition \ref{prop:locfree}]
For the first part, according to Proposition \ref{prop:muG0}, it is enough to prove that $p_*\cM_0^\mu$ is $\cO_X\lcr\theta\rcr$-locally free. One then may argue as in \cite[p\ptbl276--284]{Pham79}, using the sheaves of relative microdifferential operators $\cE_{D\times\cB/X}$ and $\cE_{D\times\cB/X\to D\times X/X}$ instead of the absolute ones, as we do not care here about the $\cD_X$-structure.

For the second part, remark that, as the tensor product is right exact, Formula~\eqref{eq:Br} shows that the restriction of $\cM_0$ at $x^o$ coincides with $\cM^o_0$. Therefore, the same property holds for the algebraization $M_0$, and then for $G_0$.
\end{proof}

\begin{theoreme}\label{th:PD}
The Poincar\'e duality morphism \ref{th:GM}\eqref{th:GM5} induces a nondegenerate Hermitian pairing $S$ of weight $n=\dim\cB/X$ on the Gauss-Manin system $(G,G_0)$.
\end{theoreme}

\begin{proof}
By Proposition~\ref{prop:microcrit}, we are reduced to a local statement on $D\times X$. Let $(c,x^o)\in\Sigma$. Choose a very small neighbourhood $D_c\times V$ of $(c,x^o)$ such that $\Phi^{-1}(D_c\times V)$ intersects small balls around each critical point above $(c,x^o)$ in a transversal way. This gives a covering of $\Phi^{-1}(D_c\times V)$ by these small balls and the complement in $\Phi^{-1}(D_c\times V)$ of smaller balls. The theorem applies to~$\Phi$ restricted to each of these open sets: for the small balls, this is the local situation at a critical point, \cf \cite{MSaito89}; for the complement of the smaller balls, $\Phi$ is then smooth and the corresponding $\cM^\mu_{x^o}$ is zero. Therefore, it applies to~$\Phi$ itself.
\end{proof}

\begin{remarque}
The microlocal sesquilinear pairing attached to $S$ is equal to the direct sum of the microlocal sesquilinear pairings attached to each critical point of $f$. It is known (\cf \cite{MSaito89}) that these pairings are equal, up to a constant, to the microlocal pairing obtained from the higher residue pairings of K.~Saito \cite{KSaito83}.
\end{remarque}

\subsection{The Malgrange-Kashiwara filtration and the spectrum} \label{subsec:MKspectre}
We keep notation of \T\ref{subsec:algcomp} and we assume that $X$ is reduced to a point. For the convenience of the reader, we briefly recall the basic definitions (see \T\ref{subsubsec:nearby} and \eg \cite[\T1]{Sabbah96b} for more details). Let $V_\bbullet\CC[\tau]\langle\partial_\tau\rangle$ be the increasing filtration of $\CC[\tau]\langle\partial_\tau\rangle$ defined by
\begin{align*}
V_{-k}\CC[\tau]\langle\partial_\tau\rangle&=\tau^k\CC[\tau]\langle\tau\partial_\tau\rangle\quad\mbox{for }k\geq 0\\
V_{k}\CC[\tau]\langle\partial_\tau\rangle&=V_{k-1}\CC[\tau]\langle\partial_\tau\rangle+\partial_\tau V_{k-1}\CC[\tau]\langle\partial_\tau\rangle\quad\mbox{for }k\geq 1.
\end{align*}
There exists a unique increasing exhaustive filtration $V_\bbullet G$ of $G$, indexed by the union of a finite number of subsets $\alpha+\ZZ\subset\QQ$, satisfying the following properties:
\begin{enumerate}
\item
For every $\alpha$, the filtration $V_{\alpha+\ZZ}G$ is good relatively to $V_\bbullet\CC[\tau]\langle\partial_\tau\rangle$;
\item
For every $\beta\in\QQ$, $\tau\partial_\tau+\beta$ is nilpotent on $\gr_{\beta}^{V}G\defin V_\beta G/V_{<\beta}G$.
\end{enumerate}
By assumption, each $V_\beta G$ is a finite type module over $\CC[\tau]\langle\tau\partial_\tau\rangle$. Because $\tau$ is invertible on $G$ the map induced by $\tau$
\[
\tau:V_\beta G\longrightarrow V_{\beta-1}G
\]
is bijective. Consequently, for every $\beta\in\QQ$ we have $\CC[\tau,\tau^{-1}]\otimes_{\CC[\tau]}V_\beta G=G$.

Consider also the filtration $G_\bbullet$ of $G$ by free $\CC[\theta]$-modules of rank $\mu$ defined by $G_k=\theta^{-k} G_0=\tau^kG_0$ for $k\in\ZZ$.

For $\beta\in\QQ$, set $$V_\beta G\cap G_0\big/\big( V_\beta\cap G_{-1}+V_{<\beta}\cap G_0\big)=\gr_\beta^V(G_0/G_{-1})$$ and let $\nu_\beta=\dim \gr_\beta^V(G_0/G_{-1})$. Notice that, for any $p$, we have an isomorphism
\[
\theta^p=\tau^{-p}:\gr_\beta^V(G_0/G_{-1})\isom\gr_{\beta-p}^V(G_p/G_{p-1}).
\]
The set of pairs $\{\beta,\nu_\beta\}$ for which $\nu_\beta\neq0$ is called the {\em spectrum} of $(G,G_0)$. The spectral polynomial of $(G,G_0)$ is
\[
\SP_f(S)\defin\prod_{\beta\in\QQ}(S+\beta)^{\nu_\beta}.
\]
Its degree is equal to $\mu$.

\section{The Frobenius structure}\label{sec:Frobenius}
We keep notation of \T\ref{subsec:setting}. We say that the family $\Phi$ is a \emph{universal unfolding} of~$f$ if the Kodaira-Spencer map $\varphi$ attached to $\Phi$ (\cf Def\ptbl\ref{def:KS}) is an isomorphism in a neighbourhood of $x^o$. We therefore restrict $X$ to such a neighbourhood. We then have $\dim X=\mu$. Throughout this section we assume that $X$ is a universal unfolding of~$f$. We will recall how one may construct on the germ at $x^o$ of $X$ a natural Frobenius structure on~$X$, following the method of K.~Saito \cite{KSaito83b} and M.~Saito \cite{MSaito89,MSaito91}, as adapted to this situation in \cite{Sabbah96c}.

\subsection{F-manifold structure}\label{subsec:weakF}
The Kodaira-Spencer isomorphism $\varphi:\Theta_X\isom q_*\cO_C$ pulls back the algebra structure of $q_*\cO_C$ to $\Theta_X$. This defines a product~$\star$ with unit on $\Theta_X$. By definition, the \emph{Euler vector field} $\gE$ of the structure is $\varphi^{-1}([F])$, where $[F]$ is the class of $F$ in $q_*\cO_C$.

Notice that the algebra $q_*\cO_C$ is the direct sum of the corresponding algebras attached to each critical point of~$f$. At this point, the structure is nothing but the direct sum of the local structures, so that we may apply the known results in this case: we have the structure of a F-manifold on $X$ (see \cite{H-M99,Hertling00}).

\subsection{Frobenius structure}\label{subsec:F}
We recall in Appendix \ref{app:B} the method of M.~Saito (\cf \cite{MSaito89}) to find a solution to Birkhoff's problem for the algebraic Brieskorn lattice $G_0^o\subset G^o$ attached to $f:\cU^o\to\Afu$, starting from any filtration of $\oplus_{\alpha\in[0,1[}\gr_\alpha^VG^o$ opposite to the filtration induced by $G^o_\bbullet$ and satisfying nice properties. Any nice solution (called $V^+$-solution in Appendix \ref{app:B}) gives a $\CC[\theta]$-basis $\varepsilong^o$ of $G_0^o$ in which the $\partial_\theta$-action takes the form
\begin{equation}\label{eq:birk}
\theta^2\partial_\theta\varepsilong^o=\varepsilong^o\cdot \big(A_0^o+\theta A_\infty\big),
\end{equation}
where $A_0^o,A_\infty$ are two constant $\mu\times\mu$ matrices, and $A_\infty$ is semisimple. The characteristic polynomial of $A_\infty$ is equal to the spectral polynomial $\SP_f(S)$ of the Brieskorn lattice of $f$ (\cf \T\ref{subsec:MKspectre}), which is also equal to the polynomial associated with the ``Hodge spectrum of $f$ at infinity'', that is, the spectrum (as defined by Steenbrink) of the limit mixed Hodge structure attached (as in \cite{Elzein86,S-Z85,MSaito87}) to $H^{n}(\cU^o,f^{-1}(t),\QQ)$ when $t\to\infty$ (\cf \cite{Sabbah96b}).

We denote by $R_\infty$ the endomorphism of $G_0^o/\theta G_0^o$ having $-A_\infty$ as matrix in the basis induced by $\varepsilong^o$.

On the other hand, recall that $G_0^o/\theta G_0^o=\Omega^{n}(\cU^o)/df\wedge\Omega^{n-1}(\cU^o)$ is a free rank-one module on the Artin algebra $q_*\cO_{C^o}$.
\begin{definition}
Let $\omega^o$ be a element of $\Omega^{n}(\cU^o)/df\wedge\Omega^{n-1}(\cU^o)$. We say that
\begin{enumerate}
\item
$\omega^o$ is \emph{primitive} if it generates $\Omega^{n}(\cU^o)/df\wedge\Omega^{n-1}(\cU^o)$ as a $q_*\cO_{C^o}$-module;
\item
$\omega^o$ is \emph{homogeneous} if it is an eigenvector of $R_\infty$.
\end{enumerate}
\end{definition}

Assume that a primitive homogeneous element $\omega^o$ exists and consider the locally free $\cO_X[\theta]$-module $G_0$: it is equipped with an integrable connection $\nabla$ with a pole of type $1$ at $\theta=0$ and a regular singularity at infinity (and no other pole). It is also equipped with a Hermitian form $S$ of weight $w$. The solution of Birkhoff's problem given by the basis $\varepsilong^o$ extends, according to theorems of B.~Malgrange \cite{Malgrange83c,Malgrange86}, to a solution in a neighbourhood of~$x^o$. This implies (see also \cite{Sabbah96c}) the existence of a flat connection $\nablaf$ on the locally free $\cO_X$-module $G_0/\theta G_0$ in some open neighbourhood of $x^o$. The pairing $S$ induces a symmetric $\nablaf$-flat nondegenerate pairing $g$ on $G_0/\theta G_0$ and $\omega^o$ extends by $\nablaf$-parallel transport to a $\nablaf$-horizontal local section $\omega$ of $G_0/\theta G_0$, called a \emph{primitive section}.

The infinitesimal period mapping $\Theta_X\isom G_0/\theta G_0$ defined by K.~Saito \cite{KSaito83b} (see also \cite{Sabbah96c}) that $\omega$ induces, is an isomorphism near $x^o$ and, transporting the structures existing on $G_0/\theta G_0$, gives rise to a Frobenius structure on $X$, compatible with the weak Frobenius structure of \T\ref{subsec:weakF}.

In conclusion, according to the previous results, in order to exhibit a Frobenius structure on $X$ compatible with the natural weak Frobenius structure, it is enough to find a $V^+$-solution to Birkhoff's problem and to construct such a primitive homogeneous element $\omega^o\in\Omega^{n}(\cU^o)/df\wedge\Omega^{n-1}(\cU^o)$.

It follows from Remark~\ref{rem:NCGM}\eqref{rem:NCGM2} that the germ of this Frobenius structure at $(X,x^o)$ does not depend on the choice of the sufficiently big ball $\cB$ defining the filtered Gauss-Manin system $(G,G_0)$.

Notice that there may exist many different Frobenius structures, attached to different choices of solutions of the Birkhoff problem or to the choice of a primitive homogeneous element $\omega^o$. We will be mainly interested to the most canonical one.

\subsection{Canonical Frobenius structure}\label{subsec:canoFrob}
We wish to fix a natural choice of the Frobenius structure. This has to be done at two levels.

A natural (or canonical) choice of the solution to Birkhoff's problem, responsible for a canonical choice of a flat connection and metric on the vector bundle $G_0/\theta G_0$. Hodge Theory furnishes such a choice, according to M.~Saito \cite[Lemme~2.8]{MSaito89}, through a natural filtration opposite to the Hodge filtration (the necessary Hodge Theory in the affine case is done in \cite{Sabbah96a,Sabbah96b}).

We will say that a primitive homogeneous element is \emph{canonical} if 
\begin{enumerate}
\item
it is an eigenvector of $R_\infty$ corresponding to the minimal element $\alpha_{\min}$ of the spectrum of $f$ (\ie the minimal exponent),
\item
up to a constant, it is the only such eigenvector, \ie $\alpha_{\min}$ has multiplicity one in the spectrum.
\end{enumerate}

That a canonical primitive homogeneous element does exist is proved by M.~Saito in \cite[Remark 3.11]{MSaito91} in the singularity case. We prove this below (\cf \T\ref{subsec:Frobnondeg}) for convenient nondegenerate Laurent polynomials on $(\CC^*)^n$. As we remarked at the end of the introduction, this also holds for convenient nondegenerate polynomials on $\CC^n$ satisfying the supplementary assumption that the coefficients of the linear forms defining the Newton boundary are nonnegative. In all these cases, the class of a volume form gives such a canonical primitive homogeneous element. The following examples illustrate phenomena which can occur when this supplementary assumption is not satisfied.

\begin{Exemples}
\begin{enumerate}
\item $f(x,y)=x+xy^{2}+y$ on $\CC^2$. There is only one spectral number $\alpha=\alpha_{\min}=1$, which therefore has multiplicity $\mu=2$. However, the class of volume form $dx\wedge dy$ is homogeneous, and furnishes a (noncanonical) primitive homogeneous element.
\item
$f(x,y)=x+y+x^2y^4$ on $\CC^2$. The class of the volume form $dx\wedge dy$ is \emph{not} homogeneous. Its order with respect to the $V$-filtration is $3/4\in{}]\alpha_{\min},\alpha_{\min}+1[$ with $\alpha_{\min}=1/2$. However, $\alpha_{\min}$ has multiplicity one in the spectrum and the corresponding eigenvector (class of $ydx\wedge dy$) is in fact primitive, hence is a canonical primitive homogeneous element.
\end{enumerate}
\end{Exemples}

\begin{remarque}
Assume that we found a canonical primitive homogeneous element $\omega^o$ as above.
Given a universal unfolding of $f$ parametrized by $(X,x^o)$, we can then construct on the base space $(X,x^o)$ a canonical Frobenius structure. When $f$ is the germ of an analytic function, two universal unfoldings are analytically isomorphic, hence the Gauss-Manin systems, with their Brieskorn lattice and sesquilinear pairings, are isomorphic, giving rise to isomorphic Frobenius structures. For functions $f$ as in \T\ref{subsec:setting}, we cannot assert that two universal unfoldings are isomorphic. Therefore, it lacks here a proof of the independence up to isomorphism of the Frobenius structure with respect to the particular choice of the unfolding.

For the example treated in the second part of this paper \cite{D-S02b}, and more generally if all the critical points of $f$ are simple and all the critical values are distinct, the Frobenius structure is semi-simple, hence completely characterized by its initial value at $x^o$ (\cf \cite[Main Theorem p\ptbl188]{Dubrovin96}, see also \cite[\T II.3]{Manin96}, \cite[Th\ptbl5.1.2]{Sabbah96c}). We therefore get the independence (up to isomorphism) with respect to the choice of the unfolding in such a case.

In general, one can expect that the Gauss-Manin system with its Brieskorn lattice is completely determined by the Gauss-Manin systems with Brieskorn lattices of the universal unfolding of each critical point plus the Stokes structure at $x^o$. Would this happen to be true, we would get the desired independence, as in the previous case.
\end{remarque}

\subsection{Basic recipe}\label{subsec:recipe}
It is possible to recover some information on the Frobenius structure on $X$ by an algebraic computation on the Gauss-Manin system of $f$. Let us indicate the recipe to get it.

{\renewcommand{\theenumi}{\alph{enumi}}
\begin{enumerate}
\item
Compute the Gauss-Manin system $G^o$ and its Brieskorn lattice $G_0^o$.
\item
Compute the ``good basis'' $\epsilong^o$ of $G_0^o$ as a $\CC[\theta]$-module.
\item
Find a primitive homogeneous element $\omega^o$, which should be part of the basis $\epsilong^o$, and denote by $\alpha$ the corresponding eigenvalue of $A_\infty$. Denote $\omega^o=\epsilon_0^o$, and let $\alpha(k)$ be the eigenvalue of $A_\infty$ corresponding to $\epsilon^o_k$, $k=0,\dots,\mu-1$ (so that $\alpha(0)=\alpha$)
\end{enumerate}}

Then we get (\cf for instance \cite{Sabbah96c} or \cite[Chap\ptbl VII]{Sabbah00}):
\begin{enumerate}
\item
there exist flat coordinates $(t_0,\dots,t_{\mu-1})$ on $X$ centered at $x^o$ such that
\begin{itemize}
\item
the basis $\partial_i ^o\in T_{x^o}X$ corresponds, under the Kodaira-Spencer map $\varphi$ and identification $q_*\cO_{C^o}\isom G^o_0/\theta G_0^o$ induced by the multiplication by $\omega^o$, to the basis $\epsilong^o$;
\item
the Euler vector field is given by
\[
\gE=\sum_{k=0}^{\mu-1}\big[(1+\alpha-\alpha(k))t_k+c_k\big]\partial_{t_k},
\]
where the $c_k$'s are the coefficient of $[f\omega^o]\in G_0^o/\theta G_0^o$ in the basis induced by $\epsilong^o$;
\end{itemize}
\item
the homogeneity constant $D$ of the structure (\ie such that $\cL_\gE(g)=Dg$, where $g$ is the ``metric'') is given by
\[
D=2\alpha+2-\dim\cU^o.
\]
\end{enumerate}

\section{The case of Laurent polynomials}\label{sec:laurent}

Let $f\in\CC[u_1,u_1^{-1},\dots,u_n,u_n^{-1}]\defin\CC[u,u^{-1}]$ be a Laurent polynomial in $n$ variables. Write $f(u)=\sum_{k\in\ZZ^n}a_ku^k$ and put $\Supp f=\{k\in\ZZ^n\mid a_k\neq0\}$. Denote by $\Gamma(f)$ the convex hull in $\RR^n$ of the set $\Supp(f)\moins\{0\}$. We will assume from now on that $f$ is \emph{nondegenerate with respect to its Newton polyhedron} and \emph{convenient} (\cf \cite{Kouchnirenko76}). In particular, $0$ belongs to the interior of $\Gamma(f)$. It is known that such an $f$ is M-tame, by applying the same reasoning as in \cite{Broughton88}.

\subsection{The Newton filtration}
For any face $\sigma$ of dimension $n-1$ of the boundary $\partial\Gamma(f)$, denote by $L_\sigma$ the linear form with coefficients in $\QQ$ such that $L_\sigma\equiv1$ on $\sigma$. For $g\in\CC[u,u^{-1}]$, put $\phi_\sigma(g)=\max_aL_\sigma(a)$, where the max is taken on the exponents of monomials appearing in $g$, and set $\phi(g)=\max_\sigma\phi_\sigma(g)$.

\begin{Remarques}\label{rem:Newton}
\begin{enumerate}
\item\label{rem:Newton1}
Of course, if $L_\sigma$ has integral coefficients, then $\phi_\sigma(g)$ is an integer for any $g\in\CC[u,u^{-1}]$.
\item\label{rem:Newton2}
For $g,h\in\CC[u,u^{-1}]$, we have
\[
\phi(gh)\leq\phi(g)+\phi(h)
\]
with equality if and only if there exists a face $\sigma$ such that $\phi(g)=\phi_\sigma(g)$ and $\phi(h)=\phi_\sigma(h)$.
\item\label{rem:Newton3}
As $0$ belongs to the interior of $\Gamma(f)$, we have $\phi(g)\geq0$ for any $g\in\CC[u,u^{-1}]$ and $\phi(g)=0$ if and only if $g\in\CC$. This would not remain true without this convenient assumption.
\item
Put $\frac{du}{u}=\frac{du_1}{u_1}\wedge\cdots\wedge\frac{du_n}{u_n}$ and let $U$ denote the torus $(\CC^*)^n$ with coordinates $u_1,\dots,u_n$. If $\omega\in\Omega^n(U)$, write $\omega=gdu/u$ and define $\phi(\omega)\defin\phi(g)$.
\end{enumerate}
\end{Remarques}

Consider the Newton increasing filtration $\cN_\bbullet\Omega^n(U)$ indexed by $\QQ$, defined by
\[
\cN_\alpha\Omega^n(U)\defin\{gdu/u\in\Omega^n(U)\mid\phi(g)\leq\alpha\}.
\]
The previous remark shows that $\cN_\alpha\Omega^n(U)=0$ for $\alpha<0$ and $\cN_0\Omega^n(U)=\CC\cdot du/u$.

Extend this filtration to $\Omega^n(U)[\theta]$ by putting
\[
\cN_\alpha\Omega^n(U)[\theta]\defin\cN_\alpha\Omega^n(U)+\theta\cN_{\alpha-1}\Omega^n(U)+\cdots+\theta^k\cN_{\alpha-k}\Omega^n(U)+\cdots
\]
and induce this filtration on $G_0$:

\begin{definition}[Newton filtration of the Brieskorn lattice]
The Newton filtration of the Brieskorn lattice is defined by
\begin{align*}
\cN_\alpha G_0&\defin \image\big[\cN_\alpha\Omega^n(U)[\theta]\hookrightarrow\Omega^n(U)[\theta] \rightarrow G_0\big]\\
&=\cN_\alpha\Omega^n(U)[\theta]\big/\big(\theta d_f\Omega^{n-1}(U)[\theta]\cap \cN_\alpha\Omega^n(U)[\theta]\big),
\end{align*}
where we have put $\theta d_f=\theta d-df\wedge$.
\end{definition}

\begin{lemme}\label{lem:Newton}
The Newton filtration on $G_0$ satisfies the following properties:
\begin{enumerate}
\item\label{lem:Newton1}
$\theta\cN_\alpha G_0\subset\cN_{\alpha+1}G_0$,
\item\label{lem:Newton2}
$\cup_\alpha\cN_\alpha G_0=G_0$,
\item\label{lem:Newton3}
$\cN_\alpha G_0=0$ if $\alpha<0$ and $\dim\cN_0G_0=1$.
\end{enumerate}
\end{lemme}

\begin{proof}
(1) is clear, (2) follows from $\cup_\alpha\cN_\alpha\Omega^n(U)=\Omega^n(U)$, and (3) follows from the similar statement for $\Omega^n(U)$.
\end{proof}

\begin{definition}[The Newton filtration on the Gauss-Manin system]
For any $\alpha\in\QQ$, we put
\[
\cN_\alpha G\defin\cN_\alpha G_0 +\tau\cN_{\alpha+1}G_0+\cdots+\tau^k\cN_{\alpha+k}G_0+\cdots.
\]
\end{definition}

For instance, we have
\[
\cN_0G=\image\big[(\cN_0\Omega^n(U)+\tau\cN_1\Omega^n(U)+\cdots+\tau^k\cN_k\Omega^n(U)+\cdots)\rightarrow G\big].
\]
From the definition, we clearly get $\tau\cN_\alpha G\subset\cN_{\alpha-1} G$, which implies that $\cN_\alpha G$ is a $\CC[\tau]$-module, and the Newton filtration on $G$ is exhaustive. Notice however that we do not have $\cN_\alpha G=0$ for $\alpha\ll0$. Nevertheless, it follows from the regularity of $G$ at $\tau=0$ and from the identification of the Newton filtration with the Malgrange-Kashiwara filtration (\cf \T\ref{subsec:MKspectre}) given by Lemma \ref{lem:VNfil} below, that we have $\cap_\alpha\cN_\alpha G_{|\tau=0}=0$ and $\cN_\alpha G=G$ out of $\tau=0$, \ie $$\CC[\tau,\tau^{-1}]\ootimes_{\CC[\tau]}\cN_\alpha G=G.$$

\begin{theoreme}\label{th:Newton}
Assume that $f$ is convenient and nondegenerate with respect to its Newton polyhedron. Then the Newton filtration $\cN_\bbullet G_0$ on the Brieskorn lattice coincides with the filtration $G_0\cap V_\bbullet G$ induced on $G_0$ by the Malgrange-Kashiwara filtration $V_\bbullet G$.
\end{theoreme}

The case of convenient nondegenerate polynomials on the affine space $\mathbb{A}^{\!n}$ has been treated in \cite{Sabbah96b}, and we will follow the proof given there. The case of germs of analytic functions which are convenient and nondegenerate goes back to \cite{K-V85,MSaito88}. However, the proof that we give here is somewhat simpler than that of \loccit., as it does not make reference to duality. One can easily adapt the simpler proof below to the case of \loccit. The proof of this theorem will be given in \T\ref{subsec:proofthNewton}.

\subsection{Division}\label{subsec:division}
Denote by $\cJ(f)$ the ideal $(u_1\partial f/\partial u_1,\dots,u_n\partial f/\partial u_n)$ of $\CC[u,u^{-1}]$.

\begin{lemme}\label{lem:Kouch}
Let $g\in\cJ(f)$. Then there exist $g_1,\dots,g_n\in \CC[u,u^{-1}]$ such that
$g=\sum_ig_iu_i\partial f/\partial u_i$ with, for any $i=1,\dots,n$,
\begin{enumerate}
\item\label{lem:Kouch2}
$\phi(g_i)\leq\phi(g)-1$,
\item\label{lem:Kouch3}
$\phi(g_iu_i\partial f/\partial u_i)\leq\phi(g)$,
\item\label{lem:Kouch4}
$\phi\big(\partial(u_ig_i)/\partial u_i\big)\leq\phi(g)-1$, and
\item\label{lem:Kouch5}
$\phi(u_ig_i)\leq\phi(g)-1+\phi(u_i)$.
\end{enumerate}
\end{lemme}

\begin{proof}
Consider with Kouchnirenko (\cf \cite[Th\ptbl 4.1]{Kouchnirenko76}) the surjective map
\begin{align*}
\partial_1:\CC[u,u^{-1}]^n&\to\cJ(f)\\
(g_1,\dots,g_n)&\mto\sum_ig_iu_i\partial f/\partial u_i.
\end{align*}
Put $\cN_\alpha(\CC[u,u^{-1}]^n)=(\cN_\alpha\CC[u,u^{-1}])^n$. As $\phi(g_iu_i\partial f/\partial u_i)\leq\phi(g_i)+1$ (\cf Remark \ref{rem:Newton}\eqref{rem:Newton2}), we have $\partial_1\cN_{\alpha-1}(\CC[u,u^{-1}]^n)\subset\cN_\alpha\CC[u,u^{-1}]\cap \cJ(f)$. According to \loccit., this inclusion is an \emph{equality}, \ie $\partial_1$ is strict for $\cN_\bbullet$. This gives \ref{lem:Kouch}\eqref{lem:Kouch2} and \eqref{lem:Kouch3}. Now, \ref{lem:Kouch}\eqref{lem:Kouch4} follows from \ref{lem:Kouch}\eqref{lem:Kouch2}, as $\phi(u_i\partial g_i/\partial u_i)\leq\phi(g_i)$, and \ref{lem:Kouch}\eqref{lem:Kouch5} follows from \ref{lem:Kouch}\eqref{lem:Kouch3}, as $\phi(u_ig_i)\leq\phi(u_i)+\phi(g_i)$.
\end{proof}

Put $E_\alpha=\gr_\alpha^\cN\big(\Omega^n(U)/df\wedge\Omega^{n-1}(U)\big)$.

\begin{proposition}\label{prop:omega}
Assume that $\omega_1,\dots,\omega_\mu\in\Omega^n(U)$ are such that their principal parts project onto a basis of $\oplus_\alpha E_\alpha$. Then, any element $\omega\in\cN_\alpha\Omega^n(U)$ may be written as
\[
\omega=\sum_ia_i\omega_i+df\wedge\eta
\]
with $a_i\in\CC$ such that $a_i=0$ if $\phi(\omega_i)>\alpha$, and with $\phi(df\wedge\eta)\leq\alpha$, $\phi(d\eta)\leq\alpha-1$.
\end{proposition}

\begin{proof}
There exist a unique family $(a_i)_{i=1,\dots,\mu}$ of complex numbers such that $\omega-\sum_ia_i\omega_i\defin gdu/u$ belongs to $\cN_\alpha\Omega^n(U)\cap (df\wedge\Omega^{n-1}(U))$; then $g$ belongs to $\cJ(f)$ and satisfies $\phi(g)\leq\alpha$. This family $(a_i)$ clearly satisfies $a_i=0$ if $\phi(\omega_i)>\alpha$. Write $g=\sum_ig_iu_i\partial f/\partial u_i$ as in Lemma \ref{lem:Kouch}, so that $\phi(u_i\partial g_i/\partial u_i)\leq\alpha-1$. Put then $\eta=\sum_i\epsilon_ig_i\wh{\frac{du_i}{u_i}}$, where $\epsilon_i$ is a suitable sign $\pm1$, so that $\omega-\sum_ia_i\omega_i=df\wedge\eta$. Remark now that $d\eta=\big(\sum_i\epsilon_iu_i\partial g_i/\partial u_i\big)\frac{du}{u}$.
\end{proof}

\begin{remarque}\label{rem:free}
It follows from this proposition that $G_0$ is $\CC[\theta]$-free (a property that we know to hold in a more general situation, \cf Proposition \ref{prop:locfree}) and, more precisely, that the classes in $G_0$ of $\omega_1,\dots,\omega_\mu$ form a $\CC[\theta]$-basis of $G_0$ adapted to the Newton filtration. Indeed, this proposition implies that any $[\omega]\in \cN_\alpha G_0$ may be written as a finite sum
\begin{equation}\label{eq:omegabis}
[\omega]=\sum_{k=0}^{[\alpha]}\ \sum_{i\mid\phi(\omega_i)\leq\alpha-k}a_i^{(k)}\theta^k[\omega_i]=\sum_{i=1}^\mu\Big(\sum_{k=0}^{\alpha-\phi(\omega_i)}a_i^{(k)}\theta^k\Big)[\omega_i].
\end{equation}
Therefore, we have found $\mu$ generators of $G_0$ over $\CC[\theta]$, giving rise to a surjective morphism $\varpi:\CC[\theta]^\mu\to G_0$. This family also generates $G$ as a $\CC[\theta,\theta^{-1}]$-module and, as $G$ is $\CC[\theta,\theta^{-1}]$-free of rank $\mu$, this family is a $\CC[\theta,\theta^{-1}]$-basis of $G$. The kernel of $\varpi$ is thus a torsion submodule, hence is $0$. In particular, the previous decomposition \eqref{eq:omegabis} is unique.

Notice also that each $\cN_\alpha G_0$ is a finite dimensional vector space.

Last, observe that if $[\omega]\in G_0$ is written as $\sum_ia_i(\theta)[\omega_i]$, then the order of $[\omega]$ with respect to the Newton filtration $\cN_\bbullet G_0$ is given by the formula
\[
\ord_\cN[\omega]=\max_i(\deg a_i+\phi(\omega_i)).
\]
\end{remarque}

We immediately deduce from this remark:
\begin{corollaire}\label{cor:omega}
We have $\cN_\alpha G_0\cap \theta G_0=\theta\cN_{\alpha-1}G_0$.\hfill\qed
\end{corollaire}

\Subsection{Proof of Theorem \ref{th:Newton}}\label{subsec:proofthNewton}

\begin{lemme}\label{lem:VNfil}
The filtration $\cN_\bbullet G$ is equal to the Malgrange-Kashiwara filtration $V_\bbullet G$.
\end{lemme}

\begin{proof}
By uniqueness of the Malgrange-Kashiwara filtration, it is enough to prove the following properties for the Newton filtration on $G$: for any $\alpha\in\QQ$,
\begin{enumerate}
\item\label{lem:VNfil1}
$\cN_\alpha G$ has finite type over $\CC[\tau]$,
\item\label{lem:VNfil2}
$\tau\cN_\alpha G\subset \cN_{\alpha-1}G$, with equality if $\alpha<1$, and $\partial_\tau\cN_{\alpha}G\subset\cN_{\alpha+1}G$,
\item\label{lem:VNfil3}
on $\gr_\alpha^\cN G$, $\tau\partial_\tau+\alpha$ is nilpotent.
\end{enumerate}
\noindent
[Notice that the ``goodness'' property $\cN_{\alpha+1}G=\partial_\tau\cN_{\alpha}G+\cN_{<\alpha+1}G$ for $\alpha>0$ follows from \eqref{lem:VNfil2} and \eqref{lem:VNfil3}.]

\smallskip
\eqref{lem:VNfil1}
As each $\cN_\beta G_0$ is a finite dimensional vector space, it is enough to show that, for a given $\alpha$, there exists $k_0$ such that, for any $k\geq k_0$, we have
\[
\tau^k\cN_{\alpha+k}G_0\subset\CC[\tau]\big(\cN_\alpha G_0+\cdots+\tau^{k_0}\cN_{\alpha+k_0}G_0\big).
\]
Choose $k_0$ such that $\cN_{\alpha+k_0}G_0+\theta G_0=G_0$ and let $k\geq k_0$. We have, according to Corollary \ref{cor:omega},
\begin{align*}
\cN_{\alpha+k}G_0=\cN_{\alpha+k_0}G_0+(\theta G_0\cap\cN_{\alpha+k}G_0)
&=\cN_{\alpha+k_0}G_0+\theta\cN_{\alpha+k-1}G_0\\
&\dots\\
&=\cN_{\alpha+k_0}G_0+\cdots+\theta^{k-k_0}\cN_{\alpha+k_0}G_0,
\end{align*}
which gives the desired inclusion.

\eqref{lem:VNfil2} We have yet seen that $\tau\cN_\alpha G\subset\cN_{\alpha-1}G$, and the inclusion $\partial_\tau\cN_\alpha G\subset\cN_{\alpha+1}G$ follows from $\phi(f)=1$, by definition of the action of $\partial_\tau$. In particular we get $\tau\partial_\tau\cN_\alpha G\subset\cN_\alpha G$. As $\cN_\alpha G_0=0$ for $\alpha<0$, we get $\cN_\alpha G=\tau\cN_{\alpha+1}G$ for $\alpha<0$.

\eqref{lem:VNfil3} Let $\sigma$ be a face of dimension $n-1$ of $\partial\Gamma(f)$. Denote by $\xi_\sigma$ the vector field $L_\sigma(u\partial_u)$. Let $g\in\CC[u,u^{-1}]$. Then we have in $G$ the relation $[\xi_\sigma(g)du/u]=\tau[\xi_\sigma(f)gdu/u]$, and therefore
\begin{multline}\label{eq:phig}
(\tau\partial_\tau+\phi_\sigma(g))[gdu/u]\\
=-\Big(\tau\big[(f-\xi_\sigma(f))gdu/u\big] +\big[(\xi_\sigma(g)-\phi_\sigma(g)g)du/u\big]\Big).
\end{multline}
(In the local case, a formula of this kind can be found in \cite{B-G-M-M89}.)
Now, for any~$g$, the support of $\xi_{\sigma}(g)-\phi_{\sigma }(g)g$ is contained in that of $g$ (it is obtained from that of $g$ by taking out
the points corresponding to monomials $u^a$ such that $\phi_{\sigma}(u^a) =\phi_{\sigma}(g)$, because, if $g$ is a monomial, we have
$\xi_{\sigma }(g)-\phi_{\sigma }(g)g=0$). Hence, for any face $\sigma'$, we have $\phi_{\sigma'}(\xi_\sigma(g)-\phi_\sigma(g)g)\leq\phi_{\sigma'}(g)$ and, similarly, $\phi_{\sigma'}\big((f-\xi_\sigma(f))g\big)\leq\phi_{\sigma'}(g)+1$ (it follows that
$\phi(\xi_{\sigma}(g)-\phi_{\sigma}(g)g)\leq \phi(g)$ and $\phi((f-\xi_{\sigma }(f))g)\leq \phi(g)+1$); moreover, these inequalities are strict if $\sigma'=\sigma$. Denote by $N(g)$ the number of faces $\sigma$ of $\partial\Gamma(f)$ such that $\phi_\sigma(g)=\phi(g)$. Applying the previous relation to any such face successively and to any initial monomial of $g$ (\ie a monomial $u$ such that $\phi(u^a)=\phi(g)$), we eventually get
\[
(\tau\partial_\tau+\phi(g))^{N(g)}[gdu/u]\in\cN_{<\phi(g)}G.\qedhere
\]
\end{proof}

\subsubsection*{End of the proof of Theorem \ref{th:Newton}}
According to Lemma \ref{lem:VNfil}, we are reduced to proving that $\cN_\alpha G_0=\cN_\alpha G\cap G_0$ for any $\alpha\in\QQ$. Let $[\omega]=v_\alpha+\tau v_{\alpha+1}+\cdots+\tau^r v_{\alpha+r}$ be in $\cN_\alpha G\cap G_0$, with each $v_{\alpha+j}$ in $\cN_{\alpha+j}G_0$. Multiplying by $\theta^r$, we find that $v_{\alpha+r}\in \cN_{\alpha+r}G_0\cap\theta G_0=\theta\cN_{\alpha+r-1}G_0$, after Corollary \ref{cor:omega}. By decreasing induction on $r$, we find that $\omega$ belongs to $\cN_\alpha G_0$.\hfill\qed

\begin{corollaire}[The spectrum]
The spectrum (or the spectral polynomial) of the Newton filtration on $G_0$ is equal to the spectrum (or the spectral polynomial) of the Malgrange-Kashiwara filtration (which is also the spectrum at infinity in the sense of Steenbrink of the Laurent polynomial $f$).\hfill\qed
\end{corollaire}

\noindent
From \cite{Sabbah96b}, we conclude:
\begin{corollaire}
The spectrum of the Newton filtration is contained in $[0,n]$ and is symmetric with respect to $n/2$.\hfill\qed
\end{corollaire}

\begin{remarque}
For a cohomologically tame or M-tame function on a smooth affine complex manifold, it can be conjectured, after C.~Hertling (for analytic germs) and A.~Dimca (for tame polynomials) that the spectrum of $(G,G_{0})$, written as \hbox{$\alpha_{1}\leq\cdots\leq\alpha_{\mu}$}, satisfies the following inequality:
\[
\frac{1}{\mu }\sum_{i=1}^{\mu }\Big(\alpha_{i}-\frac{n}{2}\Big)^{2}\geq \frac{\alpha_{\mu }-\alpha_{1}}{12}
\]
(if it is true, this inequality is the best possible). See also \cite{Brelivet02} for the case of two-variable polynomials. Notice that, for convenient and nondegenerate Laurent polynomials, one has $\alpha_1=0$ and $\alpha_\mu=n$, so that the inequality is written as
\[
\frac{1}{\mu }\sum_{i=1}^{\mu }\Big(\alpha_{i}-\frac{n}{2}\Big)^{2}\geq \frac{n}{12}.
\]
\end{remarque}

Recall that, for $\alpha\in\RR$, one denotes by $\lceil\alpha\rceil$ the smallest integer larger than or equal to $\alpha$. We also obtain:

\begin{corollaire}[A {$\CC[\tau]$}-basis of $V_0G$]
Let $\omega_1,\dots,\omega_\mu\in\Omega^n(U)$ be as in Proposition \ref{prop:omega}. For each $i=1,\dots\mu$, put $k_i=\lceil\phi(\omega_i)\rceil$. Then the classes of $\tau^{k_i}\omega_i$ in $G$ form a $\CC[\tau]$-basis of $V_0G$.
\end{corollaire}

\begin{proof}
We know that the classes of $\omega_1,\dots,\omega_\mu$, hence of $\tau^{k_1}\omega_1,\dots,\tau^{k_\mu}\omega_\mu$, form a $\CC[\tau,\tau^{-1}]$-basis of $G=\CC[\tau,\tau^{-1}]\otimes_{\CC[\tau]}V_0G$. By Nakayama, it is therefore enough to prove that the classes of $\tau^{k_i}\omega_i$ in $\gr_0^VG$ form a $\CC$-basis. This is even true in the graded space
\begin{align*}
\ooplus_k\ooplus_{\alpha\in{}]-1,0]}\gr_k^G\gr_\alpha^VG&=\ooplus_{\alpha\in{}]-1,0]}\ooplus_k\gr_\alpha^V(G_k/G_{k-1})\notag\\
&=\ooplus_{\alpha\in{}]-1,0]}\ooplus_k \gr_\alpha^\cN(G_k/G_{k-1})\quad \text{(after Th\ptbl\ref{th:Newton})}\\
&\simeq \ooplus_{\alpha}\gr_\alpha^\cN(G_0/G_{-1}).\qedhere
\end{align*}
\end{proof}

\begin{exemple}[The case of two variables]
We assume that $n=2$. It is then easy to compute the spectrum, according to the symmetry property.

First, remark that $\alpha\in\QQ\cap[0,1[$ belongs to the spectrum if and only if there exists a monomial $g\in\CC[u,u^{-1}]$ with $\phi(g)=\alpha$: indeed, assume that
\[
g=h+a_1u_1\frac{\partial f}{\partial u_1}+a_2u_2\frac{\partial f}{\partial u_2},
\]
with $\phi(h)<\alpha$ and $\phi(a_1u_1\sfrac{\partial f}{\partial u_1})\leq\alpha$ ($i=1,2$). According to Lemma \ref{lem:Kouch}, one may assume that $\phi(a_i)\leq\phi(g)-1$, hence $\phi(a_i)<0$ as $\phi(g)<1$. Therefore, $a_i=0$ and $g=h$, a contradiction.

Now, we have determined the part of the spectrum contained in $\QQ\cap [0,1[{}$: it is enough to compute $\phi$ of every monomial corresponding to a point in the interior of the polyhedron $\Gamma(f)$. By symmetry, we get the part in $\QQ\cap{}]1,2]$. The total multiplicity being equal to $\mu$, we also get the multiplicity of $1$ in the spectrum. Notice that this is the analogue of Arnold's ``butterfly''.
\end{exemple}

\subsection{The Frobenius structure}\label{subsec:Frobnondeg}
By Lemma \ref{lem:Newton}\eqref{lem:Newton3}, the smallest index of the Newton filtration is $0$, and it has multiplicity one, an eigenvector being the class of $du/u$. By Theorem \ref{th:Newton}, this also holds for the $V$-filtration, and therefore $\omega^o\defin du/u$ is the canonical homogeneous primitive element in $G_0$, as defined in \T\ref{subsec:canoFrob}.

Saito's method (\cf \T\ref{sec:Frobenius} and Appendix \ref{app:B}) gives the existence of a canonical Frobenius structure on any germ of universal unfolding of the convenient nondegenerate Laurent polynomial~$f$. The homogeneity constant is $D=2-n$.

\renewcommand{\thesection}{\Alph{section}}
\renewcommand{\thesubsubsection}{\thesubsection.\arabic{subsubsection}}
\setcounter{section}{0}
\refstepcounter{section}\label{app:A}
\section*{Appendix \ref{app:A}}

In this appendix we give the missing proofs of the results in Section \ref{sec:Fourier}. These are due to B.~Abdel-Gadir and B.~Malgrange and are adapted from \cite{AbdelGadir97a} and \cite{AbdelGadir97b}. We give them for the convenience of the reader, as these references are not published and hardly accessible (however, see also \cite{Malgrange95b}). We keep notation of Section \ref{sec:Fourier}.

\subsection{Algebraization of $\cD$-modules}\label{subsec:algebr}
Let us begin with a preliminary remark. Let $Z$ be a complex manifold and let $\Sigma$ be a divisor in $Z$. The sheaf $\cO_Z(*\Sigma)$ of meromorphic functions on $Z$ with poles on $\Sigma$ at most is a coherent sheaf of rings. The order of the pole defines a filtration by $\cO_Z$-coherent submodules. Coherent $\cO_Z(*\Sigma)$-modules have locally good filtrations. As a consequence, Cartan-Oka Theorem applies to sufficiently small compact polycylinders and a $\cO_Z(*\Sigma)$-module $\cF$ is coherent iff, for any sufficiently small compact polycylinder $K$ in $Z$, $\cF(K)$ has finite type over $\Gamma(K,\cO_Z(*\Sigma))$ and, for any $x\in K$, $\cO_{Z,x}\otimes_{\cO_Z(K)}\cF(K)\to\cF_x$ is an isomorphism.

Let now $\PP$ be a projective space of arbitrary dimension (in the next subsection, we specialize to $\PP=\PP^1$). We fix a hyperplane at infinity $H_\infty$ and still denote by $\infty$ the divisor $H_\infty\times X\subset\PP\times X$. We also fix coordinates $t$ on the affine space $\bbA=\PP\moins H_\infty$. We still denote by $p$ the projection $\PP\times X\to X$.

We review here the relationship between coherent $\cDXt$-modules (on $X$) and good $\cD_{\PP\times X}(*\infty)$-modules (on $\PP\times X$). Recall that the sheaves $\cO_X[t]$, $\cO_X[t]\langle\partial_t\rangle$ and $\cDXt$ are coherent. Similarly, the sheaves $\cO_{\PP\times X}(*\infty)$ and $\cD_{\PP\times X}(*\infty)$ are coherent. (See \cite{Malgrange75,G-M90,Kaup87}.)

Consider the following categories:
\begin{itemize}
\item
$\Mod_{\coh}(\cDXt)$ is the category of coherent left $\cDXt$-modules. Let $M$ be a $\cDXt$-module. It is coherent if and only if for any sufficiently small compact polycylinder $K\subset X$,
\begin{itemize}
\item
$\Gamma(K,M)$ is finitely generated as a $\Gamma(K,\cD_X)[t]\langle\partial_t\rangle$-module,
\item
for any $x\in K$, the natural morphism $\cO_{X,x}\otimes_{\cO(K)}\Gamma(K,M)\to M_x$ is an isomorphism.
\end{itemize}

\item
$\Mod_{\coh}\big(\cD_{\PP\times X}(*\infty)\big)$ is the category of coherent left $\cD_{\PP\times X}(*\infty)$-modules and $\Mod_{\pgood}\big(\cD_{\PP\times X}(*\infty)\big)$ is the full subcategory of objects $\cM$ such that, locally on $X$, there exists a coherent $\cO_{\PP\times X}$-module $\cF$ and a surjective morphism
\[
\cD_{\PP\times X}(*\infty)\ootimes_{\cO_{\PP\times X}}\cF\to\cM\to0.
\]
Given a $p$-good $\cD_{\PP\times X}(*\infty)$-module $\cM$, there exists, locally on $X$, a coherent $\cD_{\PP\times X}$-module $\cN$ which admits a good filtration such that $\cM$ is isomorphic to $\cD_{\PP\times X}(*\infty)\otimes_{\cD_{\PP\times X}}\cN$.
\end{itemize}

We have $p_*\big(\cD_{\PP\times X}(*\infty)\big)= \cDXt$, making $\cD_{\PP\times X}(*\infty)$ a (left and right) $p^{-1}\cDXt$-module. If $M$ is a $\cDXt$-module, we put $$p^*M=\cD_{\PP\times X}(*\infty)\ootimes_{p^{-1}\cDXt}p^{-1}M.$$
We also denote by $M^{\an}$ the restriction of $p^*M$ to the open set $\Afuan\times X$.

\begin{theoreme}[\cite{AbdelGadir97a}]
The direct image $p_*$ (partial algebraization) induces an equivalence of categories $\Mod_{\pgood}\big(\cD_{\PP\times X}(*\infty)\big)\isom\Mod_{\coh}(\cDXt)$. A quasi-inverse functor is given by $p^*$ (partial analytization).
\end{theoreme}

\begin{proof}[Sketch of proof]
One shows that $p_*$ takes values in $\Mod_{\coh}$ and that the two natural transformations $\id\to p_*p^*$ and $p^*p_*\to\id$ are isomorphisms. These assertions are local on $X$, so that we may assume that the $\pgood$ objects are generated by a $\cO_{\PP\times X}$-coherent submodule.
\begin{enumerate}
\item
One first defines analogous categories $\Mod_{\pgood}\big(\cO_{\PP\times X}(*\infty)\big)$ and $\Mod_{\coh}(\cO_X[t])$, and proves the analogous statement for these categories. The $\pgood$ objects we consider take, locally on $X$, the form $\cF(*\infty)$ where $\cF$ is $\cO_{\PP\times X}$-coherent. The theorem is a consequence of Grauert-Remmert Theorems A and B which say the following:
\begin{enonce}{}\label{th:G-R}
Let $\cF$ be a coherent $\cO_{\PP\times X}$-module. Then $\bR^k p_*\cF$ are $\cO_X$-coherent for any $k$ and for any relatively compact open subset $W$ of $X$, there exists $q_0\in\NN$ such that, for any $q\geq q_0$ and restricting to $\PP\times W$,
\begin{itemize}
\item
the natural morphism $p^*p_*\cF(q)\to\cF(q)$ is onto,
\item
$\bR^kp_*\cF(q)=0$ for $k\geq1$.
\end{itemize}
\end{enonce}
\noindent
[As usual, we denote by $\cO_{\PP\times X}(q)$ the inverse image of $\cO_{\PP}(q)$, which is the rank one bundle on $\PP$ with first Chern number equal to $q$, and we put $\cF(q)=\cO_{\PP\times X}(q)\otimes_{\cO_{\PP\times X}}\cF$.]

Notice that, by taking inductive limits, we have, as $p$ is proper, $\bR^kp_*\cF(*\infty)=0$ for $k\geq1$ and any such $\cF$.

\medskip
Fix a compact polycylinder $K$ in $X$ and take $q_0$ as in \ref{th:G-R} associated to $\cF$ and some neighbourhood of $K$. As $p_*\cF(q_0)$ is generated on $K$ by its sections (Cartan Theorem~A), we get a surjective morphism, composition of two surjective morphisms:
\[
i_{\PP\times K}^{-1}\cO_{\PP\times X}^\ell\longrightarrow i_{\PP\times K}^{-1}p^*p_*\cF(q_0)\longrightarrow i_{\PP\times K}^{-1}\cF(q_0),
\]
the kernel $\cG$ of which is coherent in some neighbourhood of $\PP\times K$. We therefore have an exact sequence on such a neighbourhood:
\[
0\to\cG(*\infty)\to\cO_{\PP\times X}(*\infty)^\ell\to\cF(*\infty)\to0.
\]
For $q$ big enough, we also have an exact sequence in some neighbourhood of $K$:
\[
0\to p_*\cG(q)\to p_*\cO_{\PP\times X}(q-q_0)^\ell\to p_*\cF(q)\to0
\]
and by taking direct limits, we get an exact sequence
\[
0\to p_*\cG(*\infty)\to \cO_X[t]^\ell\to p_*\cF(*\infty)\to 0.
\]
Arguing similarly for $\cG$, we find that, in some neighbourhood of $K$, $\cF(*\infty)$ has a free presentation
\[
\cO_{\PP\times X}(*\infty)^{\ell_1}\to\cO_{\PP\times X}(*\infty)^{\ell_0}\to\cF(*\infty)\to0
\]
so that $p_*$ of it is a presentation of $p_*\cF(*\infty)$ by free $\cO_X[t]$-modules. We conclude that $p_*\cF(*\infty)$ is $\cO_X[t]$-coherent and that the natural morphism
$p^*p_*\cF(*\infty)\to\cF(*\infty)$ is an isomorphism, as this is true for $\cF=\cO_{\PP\times X}$. That the natural morphism $N\to p_*p^*N$ is an isomorphism for a coherent $\cO_X[t]$-module follows from the projection formula and the fact that $\bR^kp_*\cO_{\PP\times X}(*\infty)=0$ for $k\geq1$.

\item
The proof for $\cD$-modules is very similar. Having fixed an object $\cM$ of $\Mod_{\pgood}(\cD_{\PP\times X}(*\infty))$, we work locally on $X$ so that we may assume that $\cM=\cN(*\infty)$ where $\cN$ is $\cD_{\PP\times X}$-coherent and has a $\cO_{\PP\times X}$-coherent submodule $\cF$ with a surjective morphism $\cD_{\PP\times X}\otimes_{\cO_{\PP\times X}}\cF\to\cN$. We therefore have a surjective morphism $\cD_{\PP\times X}\otimes_{\cO_{\PP\times X}}\cO_{\PP\times X}(-q_0)^\ell\to\cN\to0$, the kernel $\cN'$ of which is a coherent $\cD_{\PP\times X}$-module having a good filtration, by Artin-Rees lemma, hence having a $\cO_{\PP\times X}$-coherent submodule $\cG$ with a surjective morphism $\cD_{\PP\times X}\otimes_{\cO_{\PP\times X}}\cG\to\cN'$. One may end the argument as above.\qedhere
\end{enumerate}
\end{proof}

\begin{remarque}\label{rem:lattice}
Consider the categories for which objects are pairs made of a $p$-good $\cD_{\PP\times X}(*\infty)$-module and a $p$-good $\cO_{\PP\times X}(*\infty)$-submodule (\resp a coherent $\cDXt$-module with a coherent $\cO_X[t]$-submodule). The morphisms are morphisms of $\cD$-modules which send the $\cO$-submodule of the source into the $\cO$-submodule of the target.
It follows from the proof that $(p_*,p^*)$ also give an equivalence between these categories.

A similar argument applies to the category of $p$-good $\cD_{\PP\times X/X}(*\infty)$-modules and coherent $\cO_X[t]\langle\partial_t\rangle$-modules, or the same category with a $p$-good $\cO(*\Sigma)$-submodule.
\end{remarque}

\subsubsection*{Behaviour with respect to direct images, inverse images or duality}
Let $f:X\to Y$ be a holomorphic mapping and put $\wt f=(\id_{\PP}\times f):\PP\times X\to \PP\times Y$. There are direct image functors with proper support $f_\dag$ from $D^+(\cDXt)$ to $D^+(\cDYt)$, and $\wt f_\dag$ from $D^+(\cD_{\PP\times X}(*\infty))$ to $D^+(\cD_{\PP\times Y}(*\infty))$. Then, standard arguments show that they correspond each other through $p_*,p^*$. A similar result holds for inverse images of $\cD$-modules.

Let $\cM$ be in $\Mod_{\coh}\big(\cD_{\PP\times X}(*\infty)\big)$ (\resp in $\Mod_{\pgood}\big(\cD_{\PP\times X}(*\infty)\big)$), that we write locally on $\PP\times X$ (\resp locally on $X$) as $\cM=\cN(*\infty)$, where $\cN$ is an object of $\Mod_{\coh}\big(\cD_{\PP\times X}\big)$ (\resp $\Mod_{\pgood}\big(\cD_{\PP\times X}\big)$). As $\cO_{\PP\times X}(*\infty)$ is $\cO_{\PP\times X}$-flat, we have locally
\begin{align*}
\cExt_{\cD_{\PP\times X}(*\infty)}^k\big(\cM,\cD_{\PP\times X}(*\infty)\big)&=
\cExt_{\cD_{\PP\times X}}^k\big(\cN,\cD_{\PP\times X}\big)\ootimes_{\cO_{\PP\times X}} \cO_{\PP\times X}(*\infty)\\
&=\cExt_{\cD_{\PP\times X}}^k\big(\cN,\cD_{\PP\times X}\big)\ootimes_{\cD_{\PP\times X}} \cD_{\PP\times X}(*\infty),
\end{align*}
therefore $\cExt_{\cD_{\PP\times X}(*\infty)}^k\big(\cM,\cD_{\PP\times X}(*\infty)\big)$ vanishes for $k<\dim(\PP\times X)$ and otherwise belongs (after a $\text{right}\to\text{left}$ transformation) to $\Mod_{\coh}\big(\cD_{\PP\times X}(*\infty)\big)$ (\resp $\Mod_{\pgood}(\cD_{\PP\times X}(*\infty))$). Moreover, if $\cM$ is $p$-good, by considering a resolution of $p_*\cM$ by flat $\cDXt$-modules, we have isomorphisms:
\[
p_*\cExt_{\cD_{\PP\times X}(*\infty)}^k\big(\cM,\cD_{\PP\times X}(*\infty)\big)\isom \cExt_{\cDXt}^k\big(p_*\cM,\cDXt\big).
\]

\subsection{Holonomic modules}\label{subsec:hol}
Say that $\cM$ (\resp $M$) is \emph{holonomic} if, for any $k\neq\dim(\PP\times X)$, one has
\[
\cExt_{\cD_{\PP\times X}(*\infty)}^k\big(\cM,\cD_{\PP\times X}(*\infty)\big)=0\quad (\text{\resp } \cExt_{\cDXt}^k(M,\cDXt)=0).
\]

\begin{remarque}
If $\cN$ is a holonomic $\cD_{\PP\times X}$-module, it has a good filtration (\cf \cite{Malgrange92b}), hence is $p$-good. Then $\cM\defin\cN(*\infty)$ is $p$-good and holonomic as a $\cD_{\PP\times X}(*\infty)$-module. Moreover, by \cite{Kashiwara78}, $\cM$ itself is $\cD_{\PP\times X}$-holonomic (hence coherent as a $\cD_{\PP\times X}$-module and $p$-good as a $\cD_{\PP\times X}$ or $\cD_{\PP\times X}(*\infty)$-module).

Notice also that, if $\PP=\PP^1$, it is easy to show, without using \cite{Malgrange92b}, that $\cN$ is $p$-good: locally on $X$, one gets a good filtration of $\cN$ by gluing local good filtrations.
\end{remarque}

\noindent
Applying \cite{Kashiwara78}, one also gets:

\begin{proposition}\label{prop:holo}
Let $\cM$ be an object of $\Mod_{\coh}\big(\cD_{\PP\times X}(*\infty)\big)$. The following conditions are equivalent:
\begin{enumerate}
\item
$\cM$ is holonomic;
\item
considered as a $\cD_{\PP\times X}$-module, $\cM$ is holonomic (hence coherent and $p$-good);
\item
$\cM$ is $p$-good and $M=p_*\cM$ is a holonomic $\cDXt$-module;
\item
$\cM_{|\bbA^{\an}\times X}$ is a holonomic $\cD_{\bbA^{\an}\times X}$-module.
\end{enumerate}
Moreover, when these conditions are satisfied, one has
\[
\cExt_{\cD_{\PP\times X}(*\infty)}^d\big(\cM,\cD_{\PP\times X}(*\infty)\big)=
\cExt_{\cD_{\PP\times X}}^d\big(\cM,\cD_{\PP\times X}\big)(*\infty)
\]
with $d=\dim(\PP\times X)$, \ie the dual of $\cM$ as a $\cD_{\PP\times X}(*\infty)$-module is the localization at infinity of the dual of $\cM$ as a $\cD_{\PP\times X}$-module.\hfill\qed
\end{proposition}

\begin{remarque}\label{rem:holo}
As an immediate corollary of this proposition, one obtains that local analytic properties of holonomic $\cD$-modules hold for holonomic $\cDXt$-modules. In particular, if $M$ is a holonomic $\cDXt$-module,
\begin{itemize}
\item
if $Y$ is a analytic hypersurface of $X$, then $M(*Y)$ is holonomic,
\item
there exist Bernstein functional equations with respect to any holomorphic function on $X$,
\item
for any $x\in X$, the germ $M_x$ has a finite Jordan-H\"older filtration.
\end{itemize}
\end{remarque}

\begin{definition}\label{def:singetreg}
Let $M$ be a holonomic $\cDXt$-module.
\begin{enumerate}
\item\label{def:singetreg1}
Let $\ov \Sigma$ be a closed analytic hypersurface in $\PP\times X$ and put $\Sigma=\ov\Sigma\cap(\bbA\times X)$. We say that the singular locus of $M$ is contained in $\Sigma$ if $p^*\cM$ is $\cO_{\bbA^{\an}\times X}$-locally free of finite rank near any point of $(\bbA\times X)\moins\Sigma$. We denote by $\cHol_\Sigma(\cDXt)$ the category of holonomic $\cDXt$-modules with singular locus contained in~$\Sigma$.
\item\label{def:singetreg2}
We say that $M$ is regular included at infinity if the holonomic $\cD_{\PP\times X}$-module $p^*M$ is regular.
\end{enumerate}
\end{definition}

Given a holonomic $\cDXt$-module $M$, there exists, locally on $X$, a closed hypersurface $\Sigma$ containing the singular locus of $M$: indeed, locally on $X$, there exists a finite number of closed analytic subsets $Z_i$ of $\PP\times X$ such that $p^*\cM$ is $\cO$-locally free of finite rank near any point of $(\PP\times X)\moins\cup_iZ_i$ (take the family $Z_i$ so that the characteristic variety of $p^*M$ as a $\cD_{\PP\times X}$-module is equal to $T^*_{\PP\times X}(\PP\times X)\cup\bigcup_iT^*_{Z_i}(\PP\times X)$); now, each $Z_i$, being locally (on $X$) defined by homogeneous polynomials with $X$-holomorphic coefficients, is contained in a hypersurface $\ov \Sigma$ as soon as it is not equal to $\PP\times X$.

\subsubsection{Proof of Proposition \ref{prop:extension}} \label{proof:extension}
First, remark that the restriction functor, from the category $\cHol_\Sigma(\cD_{\Afuan\times X})$ to the category $\cHol_\Sigma(\cD_{D\times X})$ is an equivalence: indeed, on $D\times X\moins\Sigma$, $\cM$ is a locally free $\cO$-module of finite rank with a flat connection, hence determined by a linear representation of $\pi_1(D\times X\moins\Sigma)$, and $\pi_1(\Afuan\times X\moins\Sigma)\to \pi_1(D\times X\moins\Sigma)$ is an isomorphism. This functor also induces an equivalence between the regular holonomic objects. In what follows, we implicitly replace $D$ with $\Afuan$.

We first work on a relatively compact open set $U$ of $X$. Then there exists a disc $D'\subset \PP^1$ centered at $\infty$ such that $D'\times U\cap\Sigma=\emptyset$. On $D^{\prime*}\times U$, $\cM$ is a holomorphic bundle with a flat connection. Recall that the restriction from $D'\times U$ to $D^{\prime*}\times U$ induces an equivalence from the category of meromorphic bundles with a flat connection on $D'\times U$ with poles along $\{0'\}\times U$ to that of bundles with a flat connection on $D^{\prime*}\times U$: an inverse functor is given by Deligne's canonical extension (\cf \cite{Deligne70}, see also \cite{Malgrange87}); any morphism of bundles with flat connection extends in a unique way to a morphism of the corresponding Deligne's extensions, as horizontal sections of a flat bundle extend as meromorphic sections of the Deligne extension. Using Deligne's extension, one may extend $\cM_{\Afuan\times U}$ to some $\wt\cM_{\PP^1\times U}$. Now we may glue various extensions existing on various relatively compact sets, according to the unique lifting of morphisms.

The compatibility with duality or with base change is a direct consequence of the unique lifting of morphisms.\hfill\qed

\subsubsection{Proof of Corollary \ref{cor:locfreeextension}} \label{proof:locfreeextension}

Keep notation of \T\ref{proof:extension}. Consider the free $\cO_{D\times X}$-submodule $\cM_0$ of $\cM$. For $x^o\in X$, fix a free $\cO_{D'\times{x^o}}$-submodule $\cM^{\prime o}_0\subset\wt\cM^o_{|D'}$ in such a way that the restriction $\cM_0^o$ and $\cM^{\prime o}_0$ glue together as a trivial bundle on $\PP^1$ (it is possible to choose such a $\cM^{\prime o}_0$ according to the classification of bundles on $\PP^1$). By construction, in some neighbourhood $U$ of $x^o$, $\wt\cM_{D'\times U}$ is the inverse image by the projection $D'\times U\to D'$ of $\wt\cM^o_{|D'}$. Denote by $\cM'_0$ the inverse image of $\cM^{\prime o}_0$ by this projection. This is a locally free $\cO_{D'\times U}$-submodule of $\wt\cM_{D'\times U}$.

Restrict the situation to $\PP^1\times U$. Now, $\cM_0$ and $\cM'_0$ glue together as a vector bundle $\cN_0$ on $\PP^1\times U$, the restriction of which to $\PP^1\times \{x^o\}$ is trivial. Therefore, by rigidity of trivial bundles on $\PP^1$, $\cN_0$ is a free $\cO_{\PP^1\times V}$-module, if $V$ is some (maybe) smaller neighbourhood of $x^o$. We have $\wt\cM_0=\cO_{\PP^1\times V}(*\infty)\otimes_{\cO_{\PP^1\times V}}\cN_0$, \resp $\wt\cM=\cO_{\PP^1\times V}(*\Sigma\cup\infty)\otimes_{\cO_{\PP^1\times V}}\cN_0$, which is therefore free as a $\cO_{\PP^1\times V}(*\infty)$ (\resp $\cO_{\PP^1\times V}(*\Sigma\cup\infty)$) module, and $p_*\wt\cM_0=M_0$ (\resp $p_*\wt\cM=M$) is a free $\cO_V[t]$ (\resp $\cO_V[t](*\Sigma)$)-module.\hfill\qed

\subsubsection{Behaviour with respect to moderate nearby/vanishing cycles}\label{subsubsec:nearby}
Let $Y\subset X$ be a smooth hypersurface. Consider the $V$-filtration (increasing, indexed by $\ZZ$) of $\cDXt$ relative to $Y$ (a section $P\in \Gamma(U,\cDXt)$ is a section of $V_k\cDXt$ on $U$ iff for any $j\geq0$, $P\cdot \Gamma(U,\cI_Y^j)\subset\Gamma(U,\cI_Y^{j-k})$, where $\cI_Y$ denotes the ideal sheaf of $Y$ in $X$). There is a natural notion of a good $V$-filtration on a coherent $\cDXt$-module. This notion also exists for coherent $\cD_{\PP\times X}$ or $\cD_{\PP\times X}(*\infty)$-modules.

Any holonomic $\cD_{\PP\times X}$-module has a canonical good filtration (increasing, indexed by $\ZZ$), called Malgrange-Kashiwara filtration with respect to $Y$ (see, \eg \cite{M-S86} for more details on the $V$-filtration): it is the unique good $V$-filtration (\emph{a priori} only locally defined, but globally defined by uniqueness) such that, on each graded piece $\gr_k^V$, the operator $x\partial_x+k$ has a minimal polynomial with roots in $]-1,0]$, if $x$ is a local equation of $Y$. We then say that any holonomic $\cD_{\PP\times X}$-module is \emph{specializable} along $Y$. We define two functors $\psi_Y=\gr_0^V$ and $\phi_Y=\gr_1^V$ from holonomic $\cD_{\PP\times X}$-modules to holonomic $\cD_{\PP\times Y}$-modules (Bernstein-Kashiwara, see, \eg \cite[Cor\ptbl4.6.3]{M-S86}). The functor $\psi_Y$ is called the (moderate) \emph{nearby cycles functor} and $\phi_Y$ is called the (moderate) \emph{vanishing cycles functor}.

If the roots of the minimal polynomial of $x\partial_x+k$ on $\gr_k^V$ are real or rational, it is natural to extend the definition of the $V$-filtration so that it is indexed by $\RR$ or $\QQ$ with a (locally) finite set of jumps modulo $\ZZ$. For $\alpha\in\RR$, each $V_\alpha\cM$ is $V_0(\cD_{X\times\PP^1})$-coherent and $x\partial_x+\alpha$ is nilpotent on $\gr_\alpha^V\cM\defin V_\alpha\cM\big/V_{<\alpha}\cM$.

One may similarly define the notion of specializable $\cD_{\PP\times X}(*\infty)$ or $\cDXt$-module, and similarly define the nearby or vanishing cycles functors. Using Proposition \ref{prop:holo} and the previous results on holonomic $\cD_{\PP\times X}$-modules, as well as the uniqueness of the Malgrange-Kashiwara filtration, one gets:

\begin{corollaire}\label{cor:Vfil}
Holonomic $\cD_{\PP\times X}(*\infty)$ or $\cDXt$-modules are specializable along~$Y$. The Malgrange-Kashiwara filtrations and the nearby/vani\-shing cycles functors correspond each other under $p_*$ and $p^*$.\hfill\qed
\end{corollaire}

\subsection{Partial Fourier transform}
We now fix $\PP=\PP^1$, we denote by $t$ the coordinate on $\Afu$ and by $\tau$ the Fourier coordinate on the affine line $\Afuh$. Proposition \ref{prop:Fourierhol} clearly follows from the characterization of holonomic $\cDXt$-modules \emph{via} the vanishing of $\cExt^k$ for $k\neq\dim X+1$, as this condition is invariant by Fourier transform.

Let $M$ be $\cDXt$-holonomic. Notice that $\wh M$ is supported in $\{0\}\times X$ if and only if $M=p^+N$ for some holonomic $\cD_X$-module $N$ (this follows from Kashiwara's equivalence applied to holonomic $\cDXT$-modules supported on $\tau=0$, after partial Fourier transform), and then $\wh M$ is the direct image of $N$ by the inclusion $\{0\}\times X\hookrightarrow\Afuh\times X$. In such a case, $(\wh M)^{\an}$ is regular iff $M$ is so.

\subsubsection{Proof of Theorem \ref{th:Fourierreg}}\label{proof:Fourierreg}
Let us begin with the regularity statement. The assertion is local with respect to $X$, so we will work in some neighbourhood of a point $x^o\in X$. We will argue by induction on the dimension of $X$, the case $\dim X=0$ being treated for instance in \cite{Malgrange91}.

By induction on the dimension, we may reduce to the case where $M_{x^o}$ has no $\cDXt$-submodule supported in a strict analytic germ $(Z,x^o)\subset (X,x^o)$ (if $M$ is supported on such a set, we may assume that this set is smooth by considering a finite map $f:(Z,x^o)\to(\CC^k,0)$ and by replacing $M$ with $f_+M$).

Denote by $\Sigma_i$ ($i\in I$) the irreducible components of $\Sigma$ and by $m_i$ the multiplicity of the conormal space $T_{\Sigma_i}(\Afuan\times X)$ in the characteristic variety of $M^{\an}$. Then it is known (see, \eg \cite[Prop\ptbl1.5]{Malgrange91}) that $\wh M$ is a holonomic $\cDXT$-module having generic rank (as a $\cO$-module) equal to $\wh r(M)=\sum_im_i\deg p_{|\Sigma_i}$.

There exists an analytic hypersurface $Y\subset X$ such that, for any $x\not\in Y$, the fibre $p^{-1}(x)$ is noncharacteristic with respect to $p^*M$, and similarly there exists $\wh Y\subset X$ corresponding to $\wh M$. If $x\not\in Y$, $i_x^+M$ has only one cohomology module, hence so has $i_x^+\wh M^{\an}$, where $i_x$ corresponds to the inclusion $\{x\}\hookrightarrow X$. If moreover $x\not\in \wh Y$, then the singular locus of $i_x^+\wh M^{\an}$ is the intersection of the singular locus of $\wh M^{\an}$ with $\Afuhan\times\{x\}$. As $i_x^+M$ is regular at infinity, the singular locus of $i_x^+\wh M^{\an}$ is reduced to $\{(0,x)\}$, hence the singular locus of $\wh M^{\an}_{|X\moins Y\cup\wh Y}$ is reduced to $\{0\}\times(X\moins Y\cup\wh Y)$. The singular locus of $\wh M^{\an}$ is thus contained in $[\Afuhan\times (Y\cup \wh Y)]\cup[\{0\}\times X]$. In order to prove the regularity of $\wh M^{\an}$, it is enough to prove the regularity of the meromorphic bundle $\wh M^\an_{\rm loc}$ of rank $\wh r$ with connection obtained by localizing $\wh M^{\an}$ along its singular locus. Indeed, $\wh M^\an$ has no submodule supported in $\Afuhan\times (Y\cup \wh Y)$, and any submodule supported in $\{0\}\times X$ is regular, as we noticed above.

We assume first that $\dim X=1$. Denote by $x$ a local coordinate on $X$ centered at $x^o$. Denote by $V_\bbullet\wh M^{\an}$ the Malgrange-Kashiwara filtration of $\wh M^{\an}$ with respect to the hypersurface $x=0$. Denote simply by $\psi$ and $\phi$ the corresponding nearby and vanishing cycles functors. By uniqueness of such a filtration, we have $V_\bbullet(\wh M^{\an})=(V_\bbullet\wh M)^{\an}$, where the right-hand filtration is taken in the category of $\cDXT$-modules.

\begin{lemme}
The module $\wh M^{\an}$ is regular in a neighbourhood of $\Afuhan\times\{x^o\}$ if and only if $\psi\wh M^{\an}$ (which is holonomic on $\Afuhan\times\{x=0\}$) is generically $\cO$-locally free of rank $\wh r$.
\end{lemme}

\begin{proof}[Sketch of proof]
The ``nearby cycles'' module $\psi\wh M^{\an}$ only depends on the formal module $\cO_{\Afuan}\lcr x\rcr\otimes_{\cO_{\Afuan\times X}}\wh M^{\an}$, where we denote (somewhat inaccurately) $\cO_{\Afuan}\lcr x\rcr=\varprojlim_n\cO_{\Afuan\times X}/(x^n)$. Generically along $\Afuhan\times\{x^o\}$ one may decompose this formal module as the direct sum of a regular one and of a purely irregular one. Then $\psi$ of the purely irregular one vanishes identically. Therefore, the equality of ranks (as in the lemma) is equivalent to the vanishing of the purely irregular component, or equivalently to the fact that the formal module is regular. This, in turn, is equivalent to the regularity of $\wh M^{\an}$ itself.
\end{proof}

Remark then that
\[
\psi(\wh M^{\an})=(\psi\wh M)^{\an}=(\wh{\psi M})^{\an}.
\]
Now, the finiteness of $p$ on $\Sigma$ insures that $M$ has no ``vanishing cycles at $t=\infty$'' with respect to the function $x$, hence $\wh r(M)=\wh r(\psi M)$ (the characteristic cycle of $\psi M^{\an}$ is easily computed from that of $M^{\an}$, as $M^{\an}$ is regular).

\medskip
Let us now consider the case $\dim X\geq1$. It is a matter of verifying that the restriction of the meromorphic connection $\wh M^\an_{\rm loc}$ considered above to any disc of a holomorphic family of discs transverse to the generic part of the singular locus is regular. For the component $\{0\}\times X$, this follows from the regularity at the origin of $i_x^+\wh M^{\an}$ for any $x\not\in Y\cup\wh Y$, as $i_x^+M$ is regular at infinity (see, \eg \cite{Malgrange91}). Fix a germ of disc $\Delta\subset X$ transverse to $(Y\cup \wh Y)$ at a generic smooth point. Then $\Afuhan\times\Delta$ is noncharacteristic with respect to $\wh M^{\an}$. We apply the first part of the proof to the restriction of $M$ to $\Afu\times\Delta$.

\medskip
Now to the second part of the theorem (the noncharacteristic case). We reduce to the case where $\dim X=1$ as above. Then,
knowing that $\wh M^{\an}$ is regular, to show that $\Afuhan\times\{x^o\}$ is not contained in the singular locus of $\wh M^{\an}$ amounts to showing that the ``vanishing cycles module'' $\phi\wh M^{\an}$ vanishes identically out of $\tau=0$. Going back to $M$, this is equivalent to showing that $\phi(M)$ is isomorphic to $\CC[t]^d$ for some $d$. By the noncharacteristic assumption, we even have $\phi(M)=0$, hence the result.\hfill\qed

\begin{remarque}
In \cite{AbdelGadir97b}, the proof of Theorem \ref{th:Fourierreg} is given in the noncharacteristic case only. The proof proceeds differently: the partial Fourier transform is realized by taking the direct image of $\CC[\tau]\otimes_\CC M\otimes e^{t\tau}$ with respect to the projection $(t,\tau)\mto\tau$. The second point of the theorem is obtained by applying Kashiwara's estimate for the characteristic variety.
\end{remarque}

\subsection{A remark on duality of meromorphic connections} \label{subsec:dualmerom}
Let $Z$ be a complex manifold and let $\Sigma$ be a divisor in $Z$. A \emph{meromorphic connection} $\cM$ on $Z$ with poles on $\Sigma$ will mean a coherent $\cO_Z(*\Sigma)$-module equipped with a flat connection. By \cite{Malgrange95}, we know that such a sheaf is locally stably free on $Z$. On the other hand, viewed as a $\cD_Z$-module, we know from \cite{Kashiwara78} that $\cM$ is $\cD_Z$-holonomic. Denote by $\cM^*$ the dual meromorphic connection $\cHom_{\cO_Z(*\Sigma)}(\cM,\cO_Z(*\Sigma))$ with its natural connection, and by $\bbD\cM$ the dual left holonomic $\cD_Z$-module (may be not localized along $Z$). Then, after \cite{Kashiwara78}, the localized module $(\bbD\cM)(*\Sigma)$ is still $\cD_Z$-holonomic, and is a meromorphic connection. In the following, we put $\cO=\cO_Z(*\Sigma)$, $\Theta=\Theta_Z(*\Sigma)$, $\Omega^k=\Omega^k_Z(*\Sigma)$ and $\cD=\cD_Z(*\Sigma)$. Notice that $(\bbD\cM)(*\Sigma)$ is the dual of $\cM$ in the category of left $\cD$-modules, that we denote $\bbD_\cD\cM$.

\begin{lemme}\label{lem:dualconnexions}
There is a canonical and functorial isomorphism of meromorphic connections $\bbD_\cD\cM\simeq\cM^*$.
\end{lemme}

Notice that the lemma also gives another proof of the vanishing of $\cExt_\cD^j(\cM,\cD)$ for $j\neq\dim Z$.

\begin{proof}
Put $n=\dim Z$. The Spencer complex $\Sp^\cbbullet(\cD)$ is a left $\cD$-resolution of $\cO$ by $\cD$-locally free modules and the de~Rham complex $\Omega^{n+\cbbullet}(\cD)$ is a right $\cD$-resolution of $\omega\defin\Omega^n(*\Sigma)$. Consider the complex $\cK^\cbbullet\defin\Omega^{\cbbullet}(\cD)\otimes_{\cO}\cM$, where each term is equipped with its natural structure of right $\cD$-module. It is a complex of right $\cD$-modules. As $\cM$ is locally $\cO$-stably free, this complex is (up to a shift by $n$) a right $\cD$-resolution of $\omega\otimes_{\cO}\cM=\cM^r$ by locally stably free $\cD$-modules. Hence, the complex $\cL^{\cbbullet}=\cHom_{\cD}(\cK^{-\cbbullet},\cD)$ can be used to compute $\bbD_\cD$: we have $\bbD_\cD=\cH^0\cL^{\cbbullet}$.

Notice that the right $\cD$-module $(\Omega^k\otimes_{\cO}\cD)\otimes_{\cO}\cM$ is isomorphic to \hbox{$(\Omega^k\otimes_{\cO}\cM)\otimes_{\cO}\cD$}, where now the right structure is the trivial one.

Recall now that, if $\cF$ is any coherent $\cO$-module, then the natural morphism
\[
\cD\otimes_{\cO}\cF^*\to\cHom_\cO(\cF,\cD)=\cHom_{\cD}(\cF\otimes_\cO\cD,\cD)
\]
of left $\cD$-modules is an isomorphism (there may be an ambiguity in the middle term: let us specify that the $\cO$-module structure of $\cD$ is the right one).

Therefore, we have $\cL^k\simeq\cD\otimes_{\cO}(\wedge^{-k}\Theta\otimes_{\cO}\cM^*)$. Use now the left $\cD$ structure on $\cM^*$ to identify $\cL ^k$ to $\cM^*\otimes_{\cO}\Sp^k(\cD)$ as left $\cD$-modules. Computing the differentials gives $\cL^\cbbullet\simeq\cM^*\otimes_{\cO}\Sp^\cbbullet(\cD)$. As $\cM^*$ is locally $\cO$-stably free, this is a resolution of $\cM^*$ as a left $\cD$-module.
\end{proof}

There is a filtered variant of this lemma. We now assume that $Z=D'\times X$, where $D'$ is a disc with coordinate $\theta$ and $\Sigma=\{\theta=0\}$.

On $\cD$, consider the increasing filtration $F_\bbullet\cD$ such that $\cO_Z$ has degree $0$, each vector field on $X$ and the logarithmic vector field $\theta\partial_\theta$ have degree one, and $\theta^{-1}$ has degree one. Hence the vector field $t\defin\theta^2\partial_\theta$ has degree zero and $F_0\cD=\cO_Z\langle t\rangle$. Consider the sheaf of graded rings $R_F\cD=\oplus_kF_k\cD\cdot z^k$ on $Z$, where $z$ is a new (commuting) variable. We consider the category of graded $R_F\cD$-modules, the morphisms being graded. The sheaf $R_F\cD$ is coherent, and $\cD$-modules with good filtrations correspond to coherent graded $R_F\cD$-modules without $z$-torsion. The shift of filtrations corresponds to the multiplication by the corresponding power of $z$.

A $0$-graded $\cHom$-acyclic resolution of a $R_F\cD$-module is a resolution by $\cHom$-acyclic graded $R_F\cD$-modules for which the differentials have degree $0$. Any coherent $R_F\cD$-module has locally such resolutions, which are used to define the dual $R_F\cD$-complex, that we denote by $\bR\cHom_{R_F\cD}^0(\cbbullet,R_F\cD)$, where the exponent $0$ indicates that the differentials have degree $0$.

Denote by $\Theta'$ the sheaf of linear combinations over $R_F\cO=\cO_Z[z\theta^{-1}]$ of fields $z\xi$ ($\xi$ a vector field on $X$) and $t$. This is a locally free $R_F\cO$-module. Denote also by $\Omega^{\prime1}$ its $R_F\cO$-dual. We have a connection $\nabla:R_F\cD\to\Omega^{\prime1}\otimes_{R_F\cO}R_F\cD$, which is $z$-graded of degree $0$.

We may therefore define the de~Rham complex and the Spencer complex of $R_F\cD$, which are respectively $0$-graded $\cHom$-acyclic resolutions of $\omega'\defin\wedge^n\Omega^{\prime1}$ and $R_F\cO$. Moreover, the rules for going from left to right are similar to those for $\cD_Z$ or $\cD$.

\begin{lemme}\label{lem:dualconnexionsfiltr}
Let $R\cM$ be a coherent graded left $R_F\cD$-module which is $R_F\cO$-locally free. Then,
\begin{enumerate}
\item\label{lem:dualconnexionsfiltr1}
there exists a coherent $\cD$-module $\cM$ with a good filtration $F_\bbullet\cM$ such that $R\cM=R_F\cM$; moreover, $\cM$ is $\cO$-locally free;
\item\label{lem:dualconnexionsfiltr2}
the dual complex $\bR\cHom_{R_F\cD}^0(R\cM,R_F\cD)$ has cohomology in degree $n$ only and $\bbD_{R_F\cD}R\cM$ is strict, of the form $R_F\bbD_\cD\cM$ for some good filtration $F_\bbullet\bbD\cM$;
\item\label{lem:dualconnexionsfiltr3}
there is a canonical isomorphism $\bbD_{R_F\cD}R\cM\simeq\cHom_{R_F\cO}(R\cM,R_F\cO)$ which induces, by restriction to $z=1$, the isomorphism of Lemma~\ref{lem:dualconnexions}.
\end{enumerate}
\end{lemme}

\begin{proof}
To get \ref{lem:dualconnexionsfiltr}\eqref{lem:dualconnexionsfiltr3}, we use the same proof as in Lemma~\ref{lem:dualconnexions} with $\Omega^{\prime\cbbullet}(R_F\cD)$ and $\Sp^\cbbullet(R_F\cD)$ introduced above. By restricting to $z=1$, we get the isomorphism of Lemma~\ref{lem:dualconnexions}. Then \ref{lem:dualconnexionsfiltr}\eqref{lem:dualconnexionsfiltr1} is clear, as a locally free $R_F\cO$-module has no $z$-torsion, and \ref{lem:dualconnexionsfiltr}\eqref{lem:dualconnexionsfiltr2} follows similarly from \ref{lem:dualconnexionsfiltr}\eqref{lem:dualconnexionsfiltr3}.
\end{proof}

\begin{remarque}\label{rem:dualconnexionsfiltr}
We mainly use formal variants of these lemmas.
\begin{enumerate}
\item\label{rem:dualconnexionsfiltr1}
Let $D$ be a disc in $\CC$ and consider the sheaf $\cDXDth$ introduced in \eqref{eq:cDXDth}. For any $(c,x^o)\in D\times X$, we consider the germ $\cDXDth_{(c,x^o)}$. Introduce the sheaf $\Theta$ of vector fields, generated over $\cO_{X,x^o}\lcr\theta\rcr[\theta^{-1}]$ by $t-c$ and $\Theta_{X,x^o}$. We have $$
\cO_{X,x^o}\lcr\theta\rcr[\theta^{-1}]= \cDXDth_{(c,x^o)}/\cDXDth_{(c,x^o)}\cdot\Theta.
$$ Introduce the corresponding Spencer complex and de~Rham complex, to get the analogue of Lemma~\ref{lem:dualconnexions} for a finitely generated $\cDXDth_{(c,x^o)}$-module which is $\cO_{X,x^o}\lcr\theta\rcr[\theta^{-1}]$ free.

Filter $\cDXDth_{(c,x^o)}$ as in Lemma~\ref{lem:dualconnexionsfiltr}. Then, the analogue of this lemma holds for filtered $\cDXDth_{(c,x^o)}$-modules.

\item\label{rem:dualconnexionsfiltr2}
We may take for $\cD$ the sheaf $\cD_X[\theta,\theta^{-1}]\langle t\rangle$, with $\theta^{-1}=\partial_t$, or also its tensor product $\cD^\wedge$ by $\cO_X\lcr\theta\rcr[\theta^{-1}]$.
\end{enumerate}
\end{remarque}

\refstepcounter{section}\label{app:B}
\section*{Appendix \ref{app:B}}

We quickly review M. Saito's solution to Birkhoff's problem for the Brieskorn lattice \cite{MSaito89}. We keep notation of \T\ref{subsec:MKspectre}, which applies to any free $\CC[\tau,\tau^{-1}]$-module $G$ of rank $\mu$ equipped with a connection having singularities at $0$ and $\infty$ only, $\tau=0$ being regular, and with a lattice $G_0$, which is a free $\CC[\theta]$-submodule of rank $\mu$ stable by $t=\theta^2\partial_\theta$. The only supplementary assumption made here is that the indices of the Malgrange-Kashiwara filtration are rational numbers (the same reasoning would apply with real numbers). We denote by $N$ the nilpotent endomorphism of $H_{\alpha}\defin\gr_{\alpha}^VG$ induced by $\tau\partial_\tau+\alpha$.

\subsection{}
The following assertions are equivalent (\cf \eg \cite[\T IV.5]{Sabbah00}):
\begin{enumerate}
\item
there exists a basis $\varepsilong$ of $G_0$ satisfying \eqref{eq:birk},
\item
there exists a $\CC[\tau]$-lattice $G^{\prime 0}\subset G$ which is logarithmic (\ie stable by $\tau\partial_\tau$) such that $G_0=G_0\cap G^{\prime 0}\oplus\theta G_0$.
\end{enumerate}

On the other hand, according to the existence of a Levelt normal form (see \eg \cite[lemme II.1.2]{Sabbah00}), there is a bijective correspondence between logarithmic lattices $G^{\prime 0}$ and exhaustive decreasing filtrations $H^{\cbbullet}=\oplus_{\alpha\in[0,1[}H_{\alpha}^\bbullet$ of $H=\oplus_{\alpha\in[0,1[}H_{\alpha}$, which are stable by $N$. Put $G^{\prime k}\defin\tau^kG^{\prime 0}$. To $G^{\prime 0}$ is associated the filtration $H_{\alpha}^k=(G^{\prime k}\cap V_{\alpha}G)/(G^{\prime k}\cap V_{<\alpha}G)$.

\begin{definition}\label{def:Vsol}
A solution $G^{\prime 0}$ to Birkhoff's problem for $G_0$ is \emph{compatible with $V_\bbullet G$} (or is a \emph{$V$-solution}) if, for any $\alpha\in\QQ$, we have
\begin{equation}\tag*{(\protect\ref{def:Vsol})($*$)}\label{def:Vsol*}
G_0\cap V_{\alpha}G=(G_0\cap G^{\prime 0}\cap V_{\alpha}G)\oplus(\theta G_0\cap V_{\alpha}G).
\end{equation}
\end{definition}

\subsection{}
Recall that, for any $\alpha$, the space $G_0\cap V_{\alpha}G$ is finite dimensional (see \eg \cite{Sabbah96b}). We will repeatedly use that for $\beta\ll0$, we have $G_0\cap V_\beta=\{0\}$ and $V_\beta\subset G^{\prime 0}$. Then, by induction, \ref{def:Vsol*} is equivalent to
\begin{equation}\tag*{(\protect\ref{def:Vsol})($**$)}\label{def:Vsol**}
G_0\cap V_{\alpha}G=\ooplus_{k\geq0}(\theta^kG_0\cap \theta^kG^{\prime 0}\cap V_{\alpha}G).
\end{equation}
Denote by $\alpha_{\min}$ the minimal spectral number: $\alpha_{\min}=\min\{\alpha\mid G_0\cap V_{\alpha}G\neq0\}$. Then, if a $V$-solution $G^{\prime 0}$ exists, we have, for any $\alpha\in[\alpha_{\min},\alpha_{\min}+1[$,
\[
(G_0\cap G^{\prime 0}\cap V_{\alpha}G)=G_0\cap V_{\alpha}G\hto G_0/\theta G_0,
\]
as $(\theta G_0\cap V_{\alpha}G)=\theta(G_0\cap V_{\alpha-1}G)=0$. Hence, for such an $\alpha$, $G_0\cap G^{\prime 0}\cap V_{\alpha}G$ is independent of the choice of the $V$-solution.

Given a $V$-solution $G^{\prime 0}$, we then have
\[
(\tau\partial_\tau+\alpha)(G_0\cap G^{\prime 0}\cap V_{\alpha}G)\subset (G_0\cap G^{\prime 0}\cap V_{\alpha}G)\oplus \tau(G_0\cap G^{\prime 0}\cap V_{\alpha+1}G),
\]
which is equivalent to
\[
(\theta\partial_\theta-\alpha)(G_0\cap G^{\prime 0}\cap V_{\alpha}G)\subset (G_0\cap G^{\prime 0}\cap V_{\alpha}G)\oplus \tau(G_0\cap G^{\prime 0}\cap V_{\alpha+1}G).
\]
We then say that $G^{\prime 0}$ is a $V^+$-solution to Birkhoff's problem for $G_0$ if, for any~$\alpha$, we moreover have
\[
(\theta\partial_\theta-\alpha)(G_0\cap G^{\prime 0}\cap V_{\alpha}G)\subset (G_0\cap G^{\prime 0}\cap V_{<\alpha}G)\oplus \tau(G_0\cap G^{\prime 0}\cap V_{\alpha+1}G).
\]

\subsection{}
Assume that $G^{\prime 0}$ is a $V^+$-solution. Then, the corresponding matrix $A_\infty$, which is the first component of $\theta\partial$ in the previous decomposition, is \emph{semisimple}, and its spectrum is equal to the spectrum of $(G,G_0)$. Moreover, $G_0\cap G^{\prime 0}\cap V_{\alpha}G$ is nothing but the sum of the eigenspaces of $A_\infty$ acting on $G_0\cap G^{\prime 0}$ corresponding to the eigenvalues $\leq\alpha$. In particular, any element of $G_0\cap V_{\alpha_{\min}}G= G_0\cap G^{\prime 0}\cap V_{\alpha_{\min}}G$ is an eigenvector (with eigenvalue $\alpha_{\min}$) of~$A_\infty$.

\subsection{}
Let $S$ be a Hermitian nondegenerate sesquilinear pairing of weight $w$ on $(G,G_0)$ (\cf Definition~\ref{def:Hemitianw}). It induces nondegenerate linear pairings
\begin{equation}\label{eq:sesqH}
\begin{split}
H_{\alpha}\otimes_\CC H_{1-\alpha}&\to\CC,\qquad \alpha\in{}]0,1[,\\
H_0\otimes_\CC H_0&\to \CC.
\end{split}
\end{equation}
On the other hand, we say that a solution $G^{\prime0}$ to Birkhoff's problem for $G_0$ is \emph{compatible with $S$} (or a $S$-solution) if the induced pairing
\[
(G_0\cap G^{\prime0})\otimes_\CC(G_0\cap G^{\prime0})\to\CC[\tau,\tau^{-1}]
\]
takes values in $\CC\cdot\tau^{-w}$.

\begin{proposition}\label{prop:oppose}
Let $G^{\prime 0}$ be a logarithmic $\CC[\tau]$-lattice of $G$.
\begin{enumerate}
\item\label{prop:oppose1}
$G^{\prime 0}$ is a $V$-solution to Birkhoff's problem for $G_0$ iff for any $\alpha\in[0,1[$, its associated filtration $G^{\prime \cbbullet}\gr_{\alpha}^VG$ is opposite to $G_\bbullet\gr_{\alpha}^VG$.
\item\label{prop:oppose2}
$G^{\prime 0}$ is a $V^+$-solution iff moreover both filtrations are (B)-opposed in the sense of \cite{MSaito89}, \ie for any $\alpha\in[0,1[$ and any $k$,
\[
N(G^{\prime k}\gr_{\alpha}^VG)\subset G^{\prime k+1}\gr_{\alpha}^VG.
\]
\item
\label{prop:oppose3}
A $V$-solution $G^{\prime 0}$ is a $S$-solution iff
\[
\begin{cases}
(G^{\prime k}\gr_{\alpha}^VG)^\perp=G^{\prime w-k}\gr_{1-\alpha}^VG&\text{if } \alpha\in{}]0,1[\\
(G^{\prime k}\gr_0^VG)^\perp=G^{\prime w-k+1}\gr_0^VG,
\end{cases}
\]
where orthogonality is taken with respect to the pairing \eqref{eq:sesqH}.
\end{enumerate}
\end{proposition}

\begin{Remarques}
\begin{enumerate}
\item
A proof in the microlocal situation was given by M.~Saito in \cite[Th\ptbl3.6]{MSaito89}. An adaptation to the affine situation was given in \cite{Sabbah96b}. We sketch here a simpler proof of \eqref{prop:oppose1}. The proof of \eqref{prop:oppose2} and \eqref{prop:oppose3} is then easy, see \loccit.\ or \cite[Rem\ptbl IV.5.13]{Sabbah00}. Notice also that another proof of \eqref{prop:oppose2} can be obtained from Theorem~3.2.1 in \cite{Douai99}.
\item
According to \cite[Lemma~2.8]{MSaito89}, there is a natural choice of an opposite filtration to $G_\bbullet(\oplus_{\alpha\in[0,1[}H_\alpha)$ which gives rise to a $V^+,S$-solution to Birkhoff's problem in case $G_\bbullet H$ is the Hodge filtration of a mixed Hodge structure for which the weight filtration is obtained from the monodromy filtration of $N$ and such that the pairing \ref{eq:sesqH} is a morphism of mixed Hodge structures. That these properties are satisfied for the Gauss-Manin system, its Brieskorn lattice and the duality pairing $S$ of Theorem~\ref{th:PD}, follows from \cite{MSaito89} in the singularity case, and from \cite{Sabbah96a,Sabbah96b} in the affine case.
\end{enumerate}
\end{Remarques}

\begin{proof}
Let $H$ be any finite dimensional $\CC$-vector space equipped with an exhaustive increasing filtration $H_\bbullet$ (indexed by $\ZZ$, say) and with two exhaustive filtrations $F_\bbullet H$ (increasing) and $F^{\prime\cbbullet}H$ (decreasing). We denote by $\gr H$ the graded space associated to $H_\bbullet$. Then $F_\bbullet H$ and $F^{\prime\cbbullet}H$ naturally induce filtrations $F_\bbullet\gr H$ and $F^{\prime\cbbullet}\gr H$.

Let $u$ be a new variable. The $\CC$-vector space $\mathbb{F}=\oplus_ku^kF_kH$ is naturally equipped with the structure of a $\CC[u]$-module, as $F_\bbullet H$ is increasing, and similarly $\mathbb{F}'=\oplus_ku^kF^{\prime k}H$ is a $\CC[u^{-1}]$-module. They glue together along $\CC[u,u^{-1}]\otimes_\CC H$ as a vector bundle $\cF(F_{\bbullet}H,F^{\prime\bbullet}H)$ on $\PP^1$. The following properties are equivalent:
\begin{enumerate}
\item
there is a decomposition $H=\oplus_k(F_kH\cap F^{\prime k}H)$,
\item
for all $k\neq\ell$, one has $F_kH\cap F^{\prime \ell}H\subset (F_{k-1}H\cap F^{\prime \ell}H) + (F_kH\cap F^{\prime \ell+1}H)$,
\item
for any $k$ one has $F_{k-1}H\cap F^{\prime k}H=\{0\}$ and $F_k=(F_kH\cap F^{\prime k}H)+F_{k-1}H$,
\item
the bundle $\cF(F_\bbullet H{},F^{\prime\bbullet}H)$ is isomorphic to the trivial bundle (of rank $\dim H$).
\end{enumerate}
When any of these properties is satisfied, we say that $F_\bbullet H$ and $F^{\prime\cbbullet}H$ are \emph{opposite}.

\begin{lemme}\label{lem:2}
Assume that $F_\bbullet\gr H$ et $F^{\prime\bbullet}\gr H$ are opposite. Then so are $F_\bbullet H$ et $F^{\prime\bbullet}H$ and, for any $k$, one has $\gr(F_kH\cap F^{\prime k}H)=F_k\gr H\cap F^{\prime k}\gr H$.
\end{lemme}

\begin{proof}[Sketch of proof]
By induction, one reduces to a filtration $H^\cbbullet$ of length two, and one uses that an extension of trivial bundles on $\PP^1$ is trivial. This gives the first part. The second part reduces to a dimension count.
\end{proof}

Let us sketch the proof of \ref{prop:oppose}\eqref{prop:oppose1}. Assume that the filtrations $G_\bbullet H_{\alpha}$ and $G^{\prime\cbbullet} H_{\alpha}$ are opposite for any $\alpha\in[0,1[$. Then, $G_\bbullet\gr_\alpha^VG$ and $G^{\prime\cbbullet}\gr_\alpha^VG$ are opposite for any $\alpha$. According to Lemma \ref{lem:2}, for any $\beta<\alpha$, the filtrations $G_\bbullet(V_{\alpha}G/V_{<\beta}G)$ and $G^{\prime\cbbullet}(V_{\alpha}G/V_{<\beta}G)$ remain opposite and, for all $k$,
\[
\dim\big[G_k(V_{\alpha}G/V_{<\beta}G)\cap G^{\prime k}(V_{\alpha}G/V_{<\beta}G)\big]= \sum_{\beta\leq\gamma\leq\alpha}\dim (G_k\gr_\gamma^VG\cap G^{\prime k}\gr_\gamma^VG).
\]
For $k$ and $\alpha$ fixed, let $\beta\ll0$ such that $G_k\cap V_{<\beta}G=\{0\}$ and $V_{<\beta}G\subset G^{\prime k}$. The left-hand term above is the dimension of
\[
\big[(G_k\cap V_{\alpha}G+V_{<\beta}G)\cap G^{\prime k}+V_{<\beta}G\big]\big/V_{<\beta}G\stackrel{\sim}\longleftarrow G_k\cap G^{\prime k}\cap V_{\alpha}G.
\]
Therefore,
\[
\dim (G_k\cap G^{\prime k}\cap V_{\alpha}G) =\sum_{\gamma\leq\alpha}\dim (G_k\gr_\gamma^VG\cap G^{\prime k}\gr_\gamma^VG).
\]
Hence
\begin{align*}
\sum_{j\geq0}\dim (G_{k-j}\cap G^{\prime k-j}\cap V_{\alpha}G)&=\sum_{j\geq0}\sum_{\gamma\leq\alpha}\dim (G_{k-j}\gr_\gamma^VG\cap G^{\prime k-j}\gr_\gamma^VG)\\
&=\sum_{\gamma\leq\alpha}\sum_{j\geq0}\dim (G_{k-j}\gr_\gamma^VG\cap G^{\prime k-j}\gr_\gamma^VG)\\&=\sum_{\gamma\leq\alpha}\dim G_k\gr_\gamma^VG=\dim (G_k\cap V_{\alpha}G),
\end{align*}
because $G_\bbullet\gr_\gamma^VG$ and $G^{\prime\cbbullet}\gr_\gamma^VG$ are opposite for any $\gamma$. So, \ref{def:Vsol**} is true at the level of dimensions. It remains to show the vanishing of the two by two intersections of the terms in \ref{def:Vsol**}, and it is enough to show that, for any $k$ and any $\alpha$, $G_{k-1}\cap G^{\prime k}\cap V_{\alpha}G=0$. But the image of this term in $\gr_{\alpha}^VG$ vanishes, by oppositeness. Hence
\[
G_{k-1}\cap G^{\prime k}\cap V_{\alpha}G\subset G_{k-1}\cap G^{\prime k}\cap V_{<\alpha}G\subset\cdots\subset G_{k-1}\cap G^{\prime k}\cap V_{<\beta}G=\{0\},
\]
if $\beta\ll0$ is chosen so that $G_{k-1}\cap V_{<\beta}G=\{0\}$.

Conversely, assume that $G^{\prime 0}$ is a $V$-solution. We have to show that, for any $\alpha$ and~$k$,
\begin{gather}
\label{eq:a} (G_{k-1}\cap V_{\alpha}G+V_{<\alpha}G)\cap G^{\prime k}\subset V_{<\alpha}G\\
\label{eq:b} 
G_k\cap V_{\alpha}G\subset (G_k\cap V_{\alpha}G+V_{<\alpha}G)\cap G^{\prime k}+(G_{k-1}\cap V_{\alpha}G)+V_{<\alpha}G.
\end{gather}
Clearly, \eqref{eq:b} is a consequence of \ref{def:Vsol*}. Fix $\alpha$ and $k$ and choose $\beta\ll0$ such that $V_{<\beta}G\subset G^{\prime k}$, then choose $\ell\gg0$ such that $V_{<\alpha}G=(G_\ell\cap V_{<\alpha}G)+V_{<\beta}G$. Then
\begin{align*}
(G_{k-1}\cap V_{\alpha}G+V_{<\alpha}G)\cap G^{\prime k}&=(G_{k-1}\cap V_{\alpha}G+G_\ell\cap V_{<\alpha}G+V_{<\beta}G)\cap G^{\prime k}\\
&=(G_{k-1}\cap V_{\alpha}G+G_\ell\cap V_{<\alpha}G)\cap G^{\prime k}.
\end{align*}
On the other hand, \ref{def:Vsol**} implies  
\[
(G_{k-1}\cap V_{\alpha}+V_{<\alpha}\cap G_\ell)=\big[\oplus_{j\geq1}G_{k-j}\cap G^{\prime k-j}\cap V_{\alpha}\big]+ \big[\oplus_{i=0}^{\ell-k}G_{k+i}\cap G^{\prime k+i}\cap V_{<\alpha}\big].
\]
In the sum, the right-hand term is contained in $G^{\prime k}$, and the left-hand term, contained in $G_{k-1}\cap V_{\alpha}G$, meets $G^{\prime k}$ at $0$ only. Therefore,
\[
(G_{k-1}\cap V_{\alpha}G+G_\ell\cap V_{<\alpha}G)\cap G^{\prime k} =\big[\oplus_{i=0}^{\ell-k}G_{k+i}\cap G^{\prime k+i}\cap V_{<\alpha}G\big]\subset V_{<\alpha}G,
\]
hence \eqref{eq:a}.
\end{proof}

\backmatter
\hyphenation{Kyo-to Grund-lehren}

\providecommand{\bysame}{\leavevmode ---\ }
\providecommand{\og}{``}
\providecommand{\fg}{''}
\providecommand{\smfandname}{\&}
\providecommand{\smfedsname}{\'eds.}
\providecommand{\smfedname}{\'ed.}
\providecommand{\smfmastersthesisname}{M\'emoire}
\providecommand{\smfphdthesisname}{Th\`ese}

\end{document}